\documentclass[a4paper]{article}
\usepackage{amssymb,amscd,amsfonts,mathrsfs,amsmath}
\usepackage[all]{xy}
\input amssym.def

\def\double{\mathbb}

\def\cc{{\double C}}
\def\nn{{\double N}}
\def\zz{{\double Z}}

\def\rr{{\double R}}

\newtheorem{theorem}{Theorem}[section]
\newtheorem{lemma}[theorem]{Lemma}
\newtheorem{corollary}[theorem]{Corollary}
\newtheorem{definition}[theorem]{Definition}

\newtheorem{example}[theorem]{Example}

\def\Kt{K^{\mathrm{top}}}

\def\res{\mathop{\mathrm{Res}}\limits_{z=0}}
\def\Pf{\mathop{\mathrm{Pf}}}

\def\cp{\rtimes}
\def\si{\sigma}

\def\cinf{C^{\infty}}
\def\cinfc{C^{\infty}_c}

\newcommand{\be}{\begin{equation}}
\newcommand{\ee}{\end{equation}}
\newcommand{\beq}{\begin{eqnarray}}
\newcommand{\eeq}{\end{eqnarray}}
\newcommand{\om}{\omega}
\newcommand{\Om}{\Omega}
\newcommand{\al}{\alpha}
\def\nat{\natural}

\newcommand{\la}{\lambda}
\newcommand{\Ec}{{\mathscr E}}

\newcommand{\Lc}{{\mathscr L}}
\newcommand{\non}{\nonumber}
\newcommand{\eps}{\varepsilon}
\newcommand{\Sc}{{\mathscr S}}

\newcommand{\Rc}{{\mathscr R}}
\newcommand{\Mc}{{\mathscr M}}
\newcommand{\Nc}{{\mathscr N}}
\newcommand{\Ic}{{\mathscr I}}
\newcommand{\Jc}{{\mathscr J}}

\newcommand{\Ind}{{\mathop{\mathrm{Ind}}}}
\def\ch{\mathrm{ch}}

\def\tch{\mathrm{\,\slash\!\!\!\! \ch}}
\def\cs{\mathrm{cs}}
\def\re{\mathrm{Re}}

\newcommand{\Tr}{{\mathop{\mathrm{Tr}}}}

\newcommand{\tr}{{\mathop{\mathrm{tr}}}}
\newcommand{\Ac}{{\mathscr A}}

\newcommand{\te}{\theta}

\newcommand{\Det}{\textup{Det}}
\newcommand{\cqfd}{\hfill\rule{1ex}{1ex}}

\def\Id{\mathrm{Id}}

\def\d{\partial}
\def\dd{\mathrm{\bf d}}

\def\Bc{{\mathscr B}}
\def\Cc{{\mathscr C}}
\def\Jc{{\mathscr J}}
\def\Kc{{\mathscr K}}
\def\Fc{{\mathscr F}}
\def\Pc{{\mathscr P}}

\def\ker{\mathop{\mathrm{Ker}}}

\def\bb{\overline{b}}

\def\psib{\overline{\psi}}

\def\hom{{\mathop{\mathrm{Hom}}}}
\def\dom{{\mathop{\mathrm{Dom}}}}

\def\End{{\mathop{\mathrm{End}}}}

\def\hotimes{\hat{\otimes}}

\def\gt{\widetilde{g}}

\def\Omh{\widehat{\Omega}}
\def\Th{\widehat{T}}
\def\chih{\widehat{\chi}}
\def\etah{\widehat{\eta}}

\def\Dh{\widehat{D}}
\def\Rch{\widehat{\mathscr R}}
\def\Mch{\widehat{\mathscr M}}
\def\Jch{\widehat{\mathscr J}}
\def\Sch{\widehat{\mathscr S}}

\def\Xh{\widehat{X}}
\def\Jh{\widehat{J}}

\def\Kc{{\mathscr K}}

\def\mod{\ \mathrm{mod}\ }

\def\eh{\hat{e}}
\def\gh{\hat{g}}
\def\uh{\hat{u}}

\def\hh{\hat{h}}

\begin{document}

\begin{center}

{\bf QUASIHOMOMORPHISMS AND THE \\
RESIDUE CHERN CHARACTER}
\vskip 1cm
{\bf Denis PERROT}
\vskip 0.5cm
Universit\'e de Lyon, Universit\'e Lyon 1,\\
CNRS, UMR 5208 Institut Camille Jordan,\\
43, bd du 11 novembre 1918, 69622 Villeurbanne Cedex, France \\[2mm]
{\tt perrot@math.univ-lyon1.fr}\\[2mm]
\today
\end{center}
\vskip 0.5cm
\begin{abstract}
We develop a general procedure, based on the renormalized eta-cochain, which allows to find ``local'' representatives of the bivariant Chern character of finitely summable quasihomomorphisms. In particular, using zeta-function renormalization we obtain a bivariant generalization of the Connes-Moscovici residue formula, and explain the link with chiral and multiplicative anomalies in quantum field theory.
\end{abstract}

\vskip 0.5cm


\section{Introduction}

This article is the second part of a survey about pushforward maps between non-commutative spaces. In the first part \cite{P5} we stated a Riemann-Roch-Grothendieck theorem computing direct images of primary and secondary invariants between two Fr\'echet $m$-algebras induced by finitely summable quasihomomorphisms. By primary invariants of a non-commutative algebra we mean its usual homotopy invariants: topological $K$-theory $\Kt_n$ and periodic cyclic homology $HP_n$. The secondary invariants are multiplicative $K$-theory \cite{K1, K2} and the unstable versions of cyclic homology. Recall that a $p$-summable quasihomomorphism between two Fr\'echet $m$-algebras $\Ac$ and $\Bc$ is a continuous homomorphism $\rho: \Ac \to \Ec^s\triangleright \Ic^s\hotimes \Bc$, where $\Ec^s$ is a certain $\zz_2$-graded Fr\'echet $m$-algebra containing $\Ic^s\hotimes\Bc$ as a (not necessarily closed) two-sided ideal, and $\Ic^s$ is provided with a continuous supertrace on its $p$-th power. Under adequate assumptions about the algebras involved (existence of extensions etc...) one shows that $\rho$ induces a transformation of degree $-p$ between the long exact sequences relating topological and multiplicative $K$-theories to cyclic homology:
$$
\vcenter{\xymatrix{
\Kt_{n+1}(\Ic\hotimes\Ac) \ar[r] \ar[d]^{\rho_!} & HC_{n-1}(\Ac) \ar[r] \ar[d]^{\ch^p(\rho)} & MK^{\Ic}_n(\Ac)  \ar[r] \ar[d]^{\rho_!}  & \Kt_n(\Ic\hotimes\Ac)  \ar[d]^{\rho_!}  \\
\Kt_{n+1-p}(\Ic\hotimes\Bc) \ar[r]  & HC_{n-1-p}(\Bc) \ar[r]  & MK^{\Ic}_{n-p}(\Bc)  \ar[r]  & \Kt_{n-p}(\Ic\hotimes\Bc) }} 
$$
Here $MK^{\Ic}_n(\Ac)$ is multiplicative $K$-theory associated to the finitely summable algebra $\Ic$, see \cite{P5}. The map $\rho_!$ in multiplicative $K$-theory is a kind of non-commutative analogue of ``integrating Deligne cohomology classes'' along the fibers of a submersion. The pushforward map in cyclic homology $HC_n(\Ac)\to HC_{n-p}(\Bc)$ is induced by a bivariant non-periodic Chern character $\ch^p(\rho)\in HC^p(\Ac,\Bc)$ of degree $p$. Although $\ch^p(\rho)$ is given by a rather explicit formula, the latter is non-local and therefore of little computational use. The aim of this paper is to establish \emph{local formulas}. \\

The principle is based on the existence of a hierarchy of bivariant Chern characters $\ch^n(\rho)\in HC^n(\Ac,\Bc)$ in degrees $n=p+2k$, $k\in \nn$ associated to the quasihomomorphism, all representing the same bivariant periodic cyclic cohomology class. The representative of lowest degree $\ch^p(\rho)$ carries the information about secondary characteristic classes, whereas the higher degrees are stabilized with respect to homotopy. The difference $\ch^n(\rho)-\ch^{n+2}(\rho)=[\d, \tch^{n+1}(\rho)]$ is calculated as the boundary of a bivariant \emph{eta-cochain}. We let $\Ac_0\subset\Ac$ be a dense subalgebra and suppose that the components of the eta-cochain restricted to $\Ac_0$ can be renormalized (in a sense to be precised) in low degrees. Then a local representative of the bivariant Chern character 
\be
\ch_R(\rho) \in HC^p(\Ac_0,\Bc)
\ee
is obtained as an \emph{anomaly formula}, that is, as the boundary of the renormalized eta-cochain. For example if one chooses zeta-function renormalization, the periodic cyclic cohomology class of $\ch_R(\rho)$ is given by a residue formula generalizing the one of Connes and Moscovici \cite{CM95} to the bivariant setting (see Theorem \ref{tres}). The choice of zeta-function renormalization is by no means unique, although we use it extensively here. It is worth insisting that \emph{any} kind of renormalization would work as well, and the best choice is dictated by the geometric situation at hand. We present other types of renormalization elsewhere \cite{P6}, the important fact being that the local bivariant Chern character always appears as a boundary. This method was inspired to us by the completely analogous phenomenon of chiral anomalies occuring in quantum field theory \cite{P4}, and it is one of our goals to explain in some details the equivalence between noncommutative index theory and chiral anomalies.\\ 

As a simple illustration of the bivariant residue formula we establish the equivariant index theorem of \cite{P3}, in a simplified form. Let $G$ be a countable discrete group acting smoothly and properly by isometries on a complete Riemannian manifold $M$ without boundary, with compact quotient. Suppose given a $G$-invariant elliptic differential operator $D$ of order one acting on smooth section of a $G$-equivariant hermitian vector bundle $E\to M$. If $G$ is provided with a right-invariant distance we can form the Fr\'echet $m$-algebra $\Bc$ of functions with rapid decay over $G$, endowed with the convolution product. To these data one associates a $(p+1)$-summable quasihomomorphism $\rho:\Ac\to \Ec^s\triangleright \Ic^s\hotimes\Bc$, where $p= \dim M$ and $\Ac$ is a certain Fr\'echet $m$-algebra completion of the algebraic crossed product $\Ac_0=\cinfc(M)\cp G$. Disregarding the secondary invariants, the topological side of the Riemann-Roch-Grothendieck theorem yields a commutative diagram
$$
\vcenter{\xymatrix{
\Kt_0(\Ic\hotimes\Ac) \ar[r]^{\rho_!} \ar[d] & \Kt_{-p}(\Ic\hotimes\Bc) \ar[d] \\
HP_0(\Ac) \ar[r]^{\ch^p(\rho)} & HP_{-p}(\Bc) }}
$$
The equivariant index of $D$ is defined as the direct image $\Ind_G(D)=\rho_!(e)$ of the canonical line bundle $[e]\in \Kt_0(\Ac)$. In fact $[e]$ is represented by an idempotent $e\in\Ac_0$ and zeta-unction renormalization applies. When $E$ is a Clifford module and $D$ a generalized Dirac operator, the residue formula for the bivariant Chern character $\ch_R(\rho)$ leads to an explicit computation of the periodic Chern character $\ch(\Ind_G(D))\in HP_{-p}(\Bc)$. Its components $\ch_n(\Ind_G(D))\in \Om^n\Bc$ in all degrees $n$ are functions over the set $G^n\cup G^{n+1}$, whose evaluation at a point $\gt$ is given by a localization formula of Lefschetz type
\be
\ch_n(\Ind_G(D))(\gt) = \sum_{M_g}\frac{(-)^{q/2}}{(2\pi i)^{d/2}}\int_{M_g} \widehat{A}(M_g)\, \frac{\ch(E/S,g)}{\ch(S_N,g)}\,\ch_n(e)(\gt)\ ,
\ee
where the sum runs over all the submanifolds $M_g$ fixed by the element $g=g_n\ldots g_1\in G$ when $\gt=(g_1,\ldots,g_n)$ or $g=g_n\ldots g_0$ when $\gt=(g_0,\ldots,g_n)$. In \cite{P3} we actually proved a stronger version of the equivariant index theorem, valid for any locally compact group $G$ with proper and cocompact action on a complete manifold $M$, which yields a representative of the \emph{entire} cyclic homology class of $\ch(\Ind_G(D))$. In principle the method of proof in \cite{P3} should be adaptable to our present framework with minor modifications.\\

The above equivariant index theorem deals only with topological invariants: $\Kt_*$ and $HP_*$. In the last part of the paper we study secondary invariants for spectral triples. One thus gets quasihomomorphisms from $\Ac$ to $\Bc=\cc$ and the secondary invariants correspond to the regulator maps $\rho_!:MK^{\Ic}_{p+1}(\Ac)\to \cc^{\times}$ with values in the multiplicative group of non-zero complex numbers, already introduced by Connes and Karoubi \cite{CK} in the context of algebraic $K$-theory. However our approach by renormalization of the eta-cochain sheds some light on the link between these regulators and anomalies in quantum field theory \cite{MSZ}. In fact when $p$ is even we exhibit an explicit correspondence between the renormalized eta-cochain of a regular spectral triple and the partition function of a noncommutative chiral gauge theory. The local representative $\ch_R(\rho)$ computed as the boundary of the renormalized eta-cochain then exactly corresponds to the appearance of the chiral anomaly as the variation of the partition function under chiral gauge transformations. This allows to find an expression for the regulator map $\rho_!$ on certain multiplicative $K$-theory classes of degree ${p+1}$, in terms of renormalized determinants. \\
More precisely, let $0\to J\Ac_0\to T\Ac_0 \to \Ac_0\to 0$ be the universal free extension associated to the dense subalgebra $\Ac_0\subset \Ac$ and let $\Th\Ac_0$ be the $J\Ac_0$-adic completion of the tensor algebra. In the Cuntz-Quillen formalism for cyclic cohomology \cite{CQ1} the Chern character $\ch_R(\rho)$ is represented by a trace over the pro-algebra $\Th\Ac$. Choose a unitary element (gauge transformation) $g\in M_{\infty}(\Ac_0)^+$ homotopic to the identity in the unitalized matrix algebra over $\Ac_0$ and let $\gh\in M_{\infty}(\Th\Ac_0)^+$ be its lifting under the canonical linear inclusion $\Ac_0\hookrightarrow \Th\Ac$. Choose a transgression $\te$ of the Chern character of $g$ in cyclic homology of degree $\leq p$. Then according to \cite{P5} the pair $(\gh,\te)$ represents a class in $MK^{\Ic}_{p+1}(\Ac)$. In Corollary \ref{cdet} we calculate its pushforward 
\be
\rho_!(\gh,\te) = \exp(\sqrt{2\pi i}\,\ch_R(\rho)\cdot \te)\,\Det_R^{-1}(g) \ \in \cc^{\times}\ ,
\ee
where the determinant is renormalized in a way consistent with the eta-cochain. Hence any change of renormalization for $\Det_R$ is automatically compensated by the corresponding change in $\ch_R(\rho)$, and the result above is renormalization-independent. This interplay between $\Det_R$ and $\ch_R$ also allows to compute explicitly the \emph{multiplicative anomaly}, i.e. the lack of multiplicativity for the renormalized determinant. If $g$ and $h$ are two gauge transformations, denote by $\gh$ and $\hh$ their canonical liftings in $M_{\infty}(\Th\Ac_0)^+$ and $\widehat{gh}$ the lifting of $gh$. Then $\hh^{-1}\gh^{-1} \widehat{gh} -1$ is an element of the pro-nilpotent ideal $M_{\infty}(\Jh\Ac_0)$ and the logarithm $\ln(\hh^{-1}\gh^{-1} \widehat{gh})$ is defined by the usual formal power series. Corollary \ref{cma} gives
\be
\Det_R(gh) = \exp\big( \ch_R(\rho)\cdot \Tr\ln ( \hh^{-1}\gh^{-1} \widehat{gh} ) \big) \, \Det_R(g) \Det_R(h)
\ee
as an immediate consequence of the fact that the regulator $\rho_!:MK^{\Ic}_{p+1}(\Ac)\to \cc^{\times}$ is a morphism of abelian groups. It would be interesting to generalize these computations to target algebras $\Bc  \neq \cc$, e.g. to establish noncommutative analogues of higher analytic torsion.\\

The paper is organized as follows. In section \ref{sbiv} we recall the construction of the bivariant Chern character for quasihomomorphisms and the Riemann-Roch-Grothendieck theorem of \cite{P5}, then we explain how to obtain local representatives of the bivariant Chern character by taking the boundary of a renormalized eta-cochain. We illustrate this general procedure in section \ref{sren} with zeta-function renormalization leading to residue formulas, and treat as an example the equivariant index theorem for proper actions of discrete groups on manifolds. The link with anomalies in quantum field theory is explained in section \ref{sreg}, where we exhibit the explicit correspondence between the renormalized eta-cochain and the partition function of a noncommutative chiral gauge theory. As a consequence we interpret the index map and the Connes-Karoubi regulator of a spectral triple in the light of chiral and multiplicative anomalies.

\section{Bivariant Chern character and eta-cochain}\label{sbiv}

In this section we recall the Riemann-Roch-Grothendieck theorem of \cite{P5} and establish the general philosophy of the anomaly formula giving a ``local'' representative of the bivariant Chern character. As in \cite{P5} we work with Fr\'echet $m$-algebras, i.e. metrizable complete locally convex algebras whose topology is given by a family of submultiplicative seminorms. Equivalently, a Fr\'echet $m$-algebra is always the projective limit of a sequence of Banach algebras \cite{Ph}. \\

We say that a Fr\'echet $m$-algebra $\Ic$ is $p$-summable (with $p\geq 1$ an integer), if there is a continuous trace $\Ic^p\to \cc$ on the $p$th power of $\Ic$. By convention, $\Ic^p$ is the image of the completed projective tensor product $\Ic\hotimes\ldots \hotimes\Ic$ ($p$ times) by the multiplication map in $\Ic$. If $\Bc$ is any other Fr\'echet $m$-algebra, the completed projective tensor product $\Ic\hotimes\Bc$ is again a Fr\'echet $m$-algebra. In \cite{P5} we consider algebras $\Ec$ containing $\Ic\hotimes\Bc$ as a (not necessarily closed) two-sided ideal, subject to the condition that the semi-direct sum $\Ec\ltimes\Ic\hotimes\Bc$ is a Fr\'echet $m$-algebra. This is depicted as
\be
\Ec\triangleright \Ic\hotimes\Bc\ .
\ee
From $\Ec$ and its ideal $\Ic\hotimes\Bc$ we associate a $\zz_2$-graded algebra $\Ec^s=\Ec^s_+\oplus \Ec^s_-$ with ideal $\Ic^s\hotimes\Bc$, using $2\times 2$ matrices. In particular $\Ec^s$ comes equipped with an odd multiplier $F$ such that $F^2=1$ and $[F,\Ec^s_+]\subset \Ic^s_-\hotimes\Bc$. This construction can be done in two different ways, depending on a choice of parity:\\

\noindent Even case: let $\Ec^s_+=\Ec \ltimes \Ic\hotimes\Bc$ be the trivially graded algebra represented in $M_2(\Ec)$ by diagonal matrices $\bigl( \begin{smallmatrix} a+b & 0 \\ 0 & a \end{smallmatrix} \bigr)$ with $a\in \Ec$ and $b\in \Ic\hotimes\Bc$. Then $\Ec^s$ is the direct sum of $\Ec^s_+$ and the odd subspace $\Ec^s_-=F\Ec^s_+$ represented by off-diagonal matrices, with $F=\bigl( \begin{smallmatrix} 0 & 1 \\ 1 & 0 \end{smallmatrix} \bigr)$. The $\zz_2$-graded algebra $\Ic^s=M_2(\Ic)$ is the direct sum of the subalgebra $\Ic^s_+$ of diagonal matrices and the off-diagonal subspace $\Ic^s_-$.\\

\noindent Odd case: let $\Ec^s_+$ be the trivially graded matrix algebra $\bigl( \begin{smallmatrix} \Ec & \Ic\hotimes\Bc \\ \Ic\hotimes\Bc & \Ec \end{smallmatrix} \bigr)$. The first Clifford algebra $C_1=\cc\oplus \eps\cc$ is the $\zz_2$-graded algebra generated by the unit in degree zero and $\eps$ in degree one, with $\eps^2=1$. Then $\Ec^s$ is the tensor product $C_1\hotimes\Ec^s_+=\Ec^s_+\oplus\Ec^s_-$, with odd subspace $\Ec^s_-=\eps\Ec^s_+$. The multiplier $F$ is given by the matrix $\eps \bigl( \begin{smallmatrix} 1 & 0 \\ 0 & -1 \end{smallmatrix} \bigr)$, and $\Ic^s=C_1\hotimes M_2(\Ic)$ is the direct sum of the trivially graded subalgebra $\Ic^s_+=M_2(\Ic)$ and the odd subspace $\Ic^s_-=\eps M_2(\Ic)$.\\

\begin{definition}
Let $\Ac$, $\Bc$, $\Ic$, $\Ec$ be Fr\'echet $m$-algebras. Assume that $\Ic$ is $p$-summable and $\Ec\triangleright \Ic\hotimes\Bc$. A quasihomomorphism from $\Ac$ to $\Bc$ is a continuous homomorphism
\be
\rho: \Ac\to \Ec^s \triangleright \Ic^s\hotimes\Bc 
\ee
sending $\Ac$ to the even degree subalgebra $\Ec^s_+$. The quasihomomorphism is provided with a parity (even or odd) depending on the parity chosen for the construction of $\Ec^s$. In particular, the linear map $a\in\Ac \mapsto [F,\rho(a)]\in \Ic^s_-\hotimes\Bc$ is continuous.
\end{definition}
This definition is adapted from \cite{Cu}. Roughly speaking, a quasihomomorphism of even degree is described by a pair of homomorphisms $(\rho_+,\rho_-):\Ac\rightrightarrows \Ec$ such that the difference $\rho_+(a)-\rho_-(a)$ lies in the ideal $\Ic\hotimes\Bc$ for any $a\in\Ac$. A quasihomomorphism of odd degree is a homomorphism $\rho:\Ac\to M_2(\Ec)$ such that the off-diagonal elements land in $\Ic\hotimes\Bc$. \\

Recall that the cyclic homology of a Fr\'echet $m$-algebra $\Bc$ is computed by the cyclic bicomplex of non-commutative differential forms $\Om\Bc=\bigoplus_{n\geq 0}\Om^n\Bc$, where $\Om^n\Bc = \Bc^+\hotimes\Bc^{\hotimes n}$ is the space of $n$-forms ($\Bc^+$ denotes the unitalization of $\Bc$). The Hochschild operator $b:\Om^n\Bc\to\Om^{n-1}\Bc$ and the boundary map $B:\Om^n\Bc\to \Om^{n+1}\Bc$ are defined in the usual way \cite{C0}. One has $b^2=B^2=bB+Bb=0$ and the $\zz_2$-graded space of \emph{completed} differential forms $\Omh\Bc=\prod_{n\geq 0}\Om^n\Bc$ endowed with the boundary map $b+B$ calculates the periodic cyclic homology $HP_*(\Bc)$. The various filtrations of this complex by degree leads to the non-periodic versions of cyclic homology. See \cite{P5} for a review. An alternative description of cyclic homology is provided by the $X$-complex \cite{CQ1}. If $\Rc$ is a Fr\'echet $m$-algebra one consider the $\zz_2$-graded complex
\be
X(\Rc) \ :\ \Rc\ \mathop{\rightleftarrows}^{\nat\dd}_{\bb}\ \Om^1\Rc_{\nat}\ ,
\ee
where the odd part $\Om^1\Rc_{\nat}$ is the quotient of $\Om^1\Rc$ by the subspace of commutators $[\Rc,\Om^1\Rc]=b \Om^2\Rc$. The boundary map $\nat\dd: \Rc\to \Om^1\Rc_{\nat}$ sends an element $x\in\Rc$ to the class of the one-form $\nat\dd x$, and the boundary map $\bb: \Om^1\Rc_{\nat}\to\Rc$ comes from the Hochschild operator $\bb(\nat x\dd y)=[x,y]$. The main result of \cite{CQ1} adapted to Fr\'echet $m$-algebras states that if $\Rc$ is a \emph{quasi-free} extension of $\Bc$ with continuous linear splitting
$$
0\to \Jc \to \Rc \to \Bc \to 0\ ,
$$
the $X$-complex of the pro-algebra $\Rch=\varprojlim_n \Rc/\Jc^n$ computes the periodic cyclic homology $HP_*(\Bc)$. Note that each power $\Jc^n$ is defined as the image of $\Jc\hotimes\ldots \hotimes\Jc$ ($n$ times) by the multiplication map, and as topological vector space $\Jc^n$ is a direct summand in the quasi-free algebra $\Rc$ (see \cite{Cu1}). One gets the non-periodic versions of cyclic homology by considering the filtrations of $X(\Rc)$ by the subcomplexes
\beq 
F_{\Jc}^{2n}X(\Rc) &:& \Jc^{n+1}+[\Jc^n,\Rc] \ \rightleftarrows \ \nat \Jc^n\dd\, \Rc \label{filtration}\\
F_{\Jc}^{2n+1}X(\Rc) &:& \Jc^{n+1}\ \rightleftarrows\ \nat(\Jc^{n+1}\dd\, \Rc + \Jc^n\dd\, \Jc)\ ,\non
\eeq
where the commutator $[\Jc^n,\Rc]$ is by definition the image of $\Jc^n\dd\Rc$ under the Hochschild operator $b$, and $\Jc^n$ is defined as the unitalized algebra $\Rc^+$ for $n\leq 0$. This is a decreasing filtration and $X(\Rch)$ coincides with the projective limit of the quotient complexes
\be
X_n(\Rc,\Jc)=X(\Rc)/F^n_{\Jc}X(\Rc)\ .
\ee
The non-periodic cyclic homology of $HC_n(\Bc)$ is the homology in degree $n\mod 2$ of $X_n(\Rc,\Jc)$. Note that $X(\Rch)$ is filtered by the subcomplexes $F^n\Xh(\Rc,\Jc)=\ker(X(\Rch)\to X_n(\Rc,\Jc))$. If $\Ac$ is another algebra with quasi-free extension $0\to \Kc\to \Sc\to \Ac\to 0$, the bivariant periodic cyclic cohomology $HP^n(\Ac,\Bc)$ is the homology in degree $n\mod 2$ of the $\zz_2$-graded complex
$$
\hom(X(\Sch), X(\Rch)) = \varprojlim_m\Big(\varinjlim_k \hom(X_k(\Sc,\Kc),X_m(\Rc,\Jc))\Big) \ .
$$
The non-periodic bivariant cyclic cohomology $HC^n(\Ac,\Bc)$ is the homology in degree $n\mod 2$ of the subcomplex of bivariant cochains with order $n$:
$$
\hom^n(X(\Sch), X(\Rch)) = \{f\in \hom\ |\ \forall k, f(F^{k+n}\Xh(\Sc,\Kc)) \subset F^k\Xh(\Rc,\Jc) \}\ .
$$
Finally, one can build a chain map $\gamma: X(\Rch)\to\Omh\Rch$ realizing a homotopy equivalence between the $X$-complex and the $(b+B)$-complex of differential forms over $\Rch$. \\

Using this formalism, a Chern character in bivariant cyclic cohomology $HC^n(\Ac,\Bc)$ is constructed in \cite{P5} for any $p$-summable quasihomomorphism $\rho: \Ac\to \Ec^s \triangleright \Ic^s\hotimes\Bc $ of parity $p\mod 2$, in all degrees $n=p+2k$, $k\geq 0$, provided the algebra $\Ec\triangleright \Ic\hotimes\Bc$ satisfies some admissibility conditions. Let us choose as in \cite{P5} two quasi-free extensions of $\Ac$ and $\Bc$ admitting continuous linear splittings
$$
0\to J\Ac \to T\Ac \to \Ac \to 0\ ,\qquad 0\to \Jc \to \Rc \to \Bc \to 0\ ,
$$
where $T\Ac$ is the non-unital tensor $m$-algebra over $\Ac$, and $J\Ac$ is the kernel of the multiplication map $T\Ac\to\Ac$. The existence of a bivariant Chern character in $HC^n(\Ac,\Bc)$ is guaranteed if the algebra $\Ec\triangleright \Ic\hotimes\Bc$ admits a $\Rc$-admissible extension in the following sense: there exists two algebras $\Mc\triangleright \Ic\hotimes\Rc$ and $\Nc\triangleright \Ic\hotimes\Jc$ and a commutative diagram of extensions (with continuous linear splitting)
\be
\vcenter{\xymatrix{
0\ar[r] & \Nc \ar[r] & \Mc \ar[r] & \Ec \ar[r] & 0 \\
0 \ar[r] & \Ic\hotimes \Jc \ar[r] \ar[u] & \Ic\hotimes \Rc \ar[r] \ar[u] & \Ic\hotimes\Bc \ar[r] \ar[u] & 0}}
\label{adm}
\ee
with adequate properties listed in \cite{P5}, Definition 3.2. Note that the algebra $\Mc$ is not necessarily quasi-free. The bivariant Chern character in $HC^n(\Ac,\Bc)$ is represented by a chain map between the $X$-complexes $X(\Th\Ac)$ and $X(\Rch)$, compatible with the adic filtrations induced by the ideals $J\Ac$ and $\Jc$. The first step is to lift the quasihomomorphism to the quasi-free algebras $T\Ac$ and $\Rc$, using the diagram of extensions (\ref{adm}). To $\Mc\triangleright \Ic\hotimes\Rc$ we associate the $\zz_2$-graded algebra $\Mc^s\triangleright \Ic^s\hotimes\Rc$ with multiplier $F$ such that $[F,\Mc^s_+]\subset \Ic^s_-\hotimes\Rc$. In the same way, we have $\Nc^s\triangleright\Ic^s\hotimes\Jc$ with $[F,\Nc^s_+]\subset \Ic^s_-\hotimes\Jc$. Then we get an extension of $\zz_2$-graded algebras
$$
0\to \Nc^s \to \Mc^s \to \Ec^s \to 0\ ,
$$
which are moreover differential algebras for the differential operator of odd degree induced by the graded commutator $[F,\ ]$. The restriction to the even-degree subalgebras yields an extension of trivially graded algebras $0\to \Nc^s_+ \to \Mc^s_+ \to \Ec^s_+ \to 0$, split by a continuous linear map $\si:\Ec^s_+\to \Mc^s_+$ by hypothesis. The universal property of the tensor algebra $T\Ac$ then allows to extend the homomorphism $\rho: \Ac\to \Ec^s_+$ to a continuous homomorphism $\rho_*:T\Ac\to \Mc^s_+$ by setting $\rho_*(a_1\otimes\ldots\otimes a_n)=\si\rho(a_1)\otimes\ldots\otimes \si\rho(a_n)$:
\be
\vcenter{\xymatrix{
0 \ar[r]  & J\Ac \ar[r] \ar[d]_{\rho_*} & T\Ac  \ar[r] \ar[d]_{\rho_*} & \Ac \ar[r]  \ar[d]^{\rho} & 0  \\
0 \ar[r] & \Nc^s_+ \ar[r] & \Mc^s_+ \ar[r] & \Ec^s_+ \ar[r] \ar@/_/@{.>}[l]_{\si} & 0 }} \label{uni}
\ee
A priori $\rho_*$ depends on the choice of linear splitting $\si$, but in a way which will not affect the cohomology class of the bivariant Chern character. This construction may be depicted in terms of a $p$-summable quasihomomorphism $\rho_*: T\Ac \to \Mc^s\triangleright \Ic^s\hotimes \Rc$, compatible with the adic filtration by the ideals in the sense that $J\Ac$ is mapped to $\Nc^s\triangleright \Ic^s\hotimes\Jc$. Hence, $\rho_*$ extends to a quasihomomorphism of pro-algebras
\be
\rho_* : \Th\Ac \to \Mch^s\triangleright \Ic^s\hotimes \Rch\ ,
\ee
where $\Th\Ac$, $\Mch^s$ and $\Rch$ are the adic completions with respect to the ideals $J\Ac$, $\Nc^s$ and $\Jc$ respectively. The bivariant Chern character $\ch^n(\rho)\in HC^n(\Ac,\Bc)$ then exists in any degree $n=p+2k$, $k\geq 0$ and is represented by the composition of chain maps 
\be
\ch^n(\rho): X(\Th\Ac) \stackrel{\gamma}{\longrightarrow} \Omh \Th\Ac \stackrel{\rho_*}{\longrightarrow} \Omh\Mch^s_+ \stackrel{\chih^n}{\longrightarrow}  X(\Rch)\ ,
\ee
where $\gamma$ is the homotopy equivalence and the homomorphism $\rho_*$ extends to non-commutative differential forms in any degree. The chain map $\chih^n$ has only two non-zero components $\chih^n_0:\Om^n \Mch^s_+\to \Rch$ and $\chih^n_1:\Om^{n+1}\Mch^s_+\to \Om^1\Rch_{\nat}$ defined by
$$
\chih^n_0(x_0\dd x_1\ldots\dd x_n) = (-)^n\frac{\Gamma(1+\frac{n}{2})}{(n+1)!} \sum_{\la\in S_{n+1}} \eps(\la)\, \tau(x_{\la(0)}[F,x_{\la(1)}]\ldots [F,x_{\la(n)}])
$$
$$
\chih^n_1(x_0\dd x_1\ldots\dd x_{n+1}) = (-)^n\frac{\Gamma(1+\frac{n}{2})}{(n+1)!} \sum_{i=1}^{n+1}  \tau\nat(x_0[F,x_1]\ldots\dd x_i \ldots [F,x_{n+1}])
$$
where $\tau$ is induced by the supertrace on the $n$-th power of $\Ic^s$ when $n\geq p$. $S_{n+1}$ is the cyclic permutation group over $n+1$ elements. For any $x\in\Mch^s_+$, the commutator $[F,x]$ lies in the ideal $\Ic^s\hotimes\Rch$ hence $\chih^n$ is well-defined whenever $n\geq p$. See \cite{P5} for details. In fact it is possible to weaken the summability degree and require $\Ic$ to be only $(p+1)$-summable, while replacing $\tau$ by the supertrace $\tau'=\frac{1}{2}\tau(F[F,\ ])$. Note that the parity of these supertraces, as well as the parity of the cocycle $\chih^n\in \hom(\Omh\Mch^s_+,X(\Rch))$ coincide with the parity of the quasihomomorphism. As shown in \cite{P5}, the composite map $\ch^n(\rho)=\chih^n\rho_*\gamma$ is compatible with the adic filtrations induced by the ideals $J\Ac$, $\Jc$ and defines a cocycle of order $n$ in the $\zz_2$-graded complex $\hom(X(\Th\Ac),X(\Rch))$ representing a bivariant cyclic cohomology class in $HC^n(\Ac,\Bc)$. The cocycles corresponding to successive values of $n\geq p$ are related by the transgression formula involving the \emph{eta-cochain}
\be
\chih^n-\chih^{n+2}= [\partial , \etah^{n+1}]\quad \mbox{in}\quad \hom(\Omh \Mch^s_+,X(\Rch))\ ,\label{tra}
\ee
where $\partial$ denotes either the boundary operator on $\Omh \Mch^s_+$ or on $X(\Rch)$, and $[\ ,\ ]$ is the graded commutator. The eta-cochain $\etah^{n+1}$ is a cochain in the complex $\hom(\Omh\Mch^s_+,X(\Rch))$ of parity opposite to the quasihomomorphism, whose only non-zero components $\etah^{n+1}_0:\Om^{n+1}\Mch^s_+\to \Rch$ and $\etah^{n+1}_1:\Om^{n+2}\Mch^s_+\to \Om^1\Rch_{\nat}$ are given by
\beq
\lefteqn{\etah^{n+1}_0(x_0\dd x_1\ldots\dd x_{n+1}) =  \frac{\Gamma(\frac{n}{2}+1)}{(n+2)!} \, \frac{1}{2}\tau\Big (F x_0[F,x_1]\ldots [F,x_{n+1}]+ } \non \\
&&\qquad \qquad \qquad  \sum_{i=1}^{n+1}(-)^{(n+1)i} [F,x_i]\ldots [F,x_{n+1}] Fx_0 [F,x_1]\ldots [F,x_{i-1}] \Big) \non
\eeq
\beq
\lefteqn{\etah^{n+1}_1(x_0\dd x_1\ldots\dd x_{n+2}) =} \non\\
&& \frac{\Gamma(\frac{n}{2}+1)}{(n+3)!} \sum_{i=1}^{n+2}  \frac{1}{2}\tau\nat(ix_0 F + (n+3-i)Fx_0)[F,x_1]\ldots\dd x_i \ldots [F,x_{n+2}]\ .\non
\eeq
Again the composite map $\tch^{n+1}(\rho)=\etah^{n+1}\rho_*\gamma$ has good adic properties and the transgression formula $\ch^n(\rho)-\ch^{n+2}(\rho)=[\partial , \tch^{n+1}(\rho)]$ holds as cochains of order $n+2$ in the $\zz_2$-graded complex $\hom(X(\Th\Ac),X(\Rch))$. This implies the equality $S\ch^n(\rho)=\ch^{n+2}(\rho)$ in $HC^{n+2}(\Ac,\Bc)$. In particular the Chern characters associated to different values of $n$ represent the same bivariant periodic cyclic cohomology class $\ch(\rho)\in HP^i(\Ac,\Bc)$, $i\equiv p\mod 2$.\\
Another important property of the bivariant Chern character is its invariance under two types of equivalence relations among quasihomomorphisms, namely conjugation and homotopy (\cite{P5}, Proposition 3.11). If two quasihomomorphisms $\rho$ and $\rho'$ are conjugate, then $\ch^n(\rho)=\ch^n(\rho')$ in $HC^n(\Ac,\Bc)$ in any degree $n$. If however $\rho$ and $\rho'$ are homotopic, one has only $S\ch^n(\rho)=S\ch^n(\rho')$ in $HC^{n+2}(\Ac,\Bc)$. Thus in any case the periodic cyclic cohomology class is invariant.\\

The bivariant Chern character was used in \cite{P5} to define direct images of secondary characteristic classes for Fr\'echet $m$-algebras. Recall that for any such algebra $\Ac$ and finitely summable $\Ic$, one has a long exact sequence of abelian groups
$$
\Kt_{n+1}(\Ic\hotimes\Ac)\to HC_{n-1}(\Ac)\to MK^{\Ic}_n(\Ac)\to \Kt_n(\Ic\hotimes\Ac)\to HC_{n-2}(\Ac)
$$
where $\Kt_*$ is the topological $K$-theory of Phillips \cite{Ph}, $HC_*$ is non-periodic cyclic homology, and $MK^{\Ic}_*$ is multiplicative $K$-theory. The Riemann-Roch-Grothendieck theorem of \cite{P5} establishes the way a quasihomomorphism $\rho: \Ac\to \Ec^s \triangleright \Ic^s\hotimes\Bc $ maps the long exact sequence for $\Ac$ to the long exact sequence for $\Bc$. We say that the finitely summable algebra $\Ic$ is multiplicative if there exists an external product $\boxtimes:\Ic\hotimes\Ic\to \Ic$ compatible with the traces (\cite{P5} Definition 6.1). Then one has

\begin{theorem}[\cite{P5}]\label{trr}
Let $\Ac$, $\Bc$ be Fr\'echet $m$-algebras, and choose a quasi-free extension $0 \to \Jc \to \Rc \to \Bc \to 0$. Let $\rho:\Ac\to \Ec^s\triangleright \Ic^s\hotimes\Bc$ be a quasihomomorphism of parity $p \mod 2$, where $\Ic$ is multiplicative and $(p+1)$-summable in the even case, $p$-summable in the odd case. The algebra $\Ec\triangleright \Ic\hotimes\Bc$ is supposed $\Rc$-admissible. Then $\rho$ induces a graded-commutative diagram 
$$
\vcenter{\xymatrix{
\Kt_{n+1}(\Ic\hotimes\Ac) \ar[r] \ar[d]^{\rho_!} & HC_{n-1}(\Ac) \ar[r] \ar[d]^{\ch^p(\rho)} & MK^{\Ic}_n(\Ac)  \ar[r] \ar[d]^{\rho_!}  & \Kt_n(\Ic\hotimes\Ac)  \ar[d]^{\rho_!}  \\
\Kt_{n+1-p}(\Ic\hotimes\Bc) \ar[r]  & HC_{n-1-p}(\Bc) \ar[r]  & MK^{\Ic}_{n-p}(\Bc)  \ar[r]  & \Kt_{n-p}(\Ic\hotimes\Bc) }}  
$$
The vertical arrows are invariant under conjugation of quasihomomorphisms; the arrow in topological $K$-theory $\Kt_n(\Ic\hotimes\Ac)\to\Kt_{n-p}(\Ic\hotimes\Bc)$ is also invariant under homotopy of quasihomomorphisms. 
\end{theorem}

Let us now describe the simple trick allowing to obtain \emph{local formulas} representing the cyclic cohomology class of the bivariant Chern character. We suppose that the quasihomomorphism $\rho: \Ac\to \Ec^s \triangleright \Ic^s\hotimes\Bc $ is of parity $p\mod 2$ and summability degree $(p+1)$. We collect all the components of the eta-cochain $\etah^{n+1}$, $n\geq p$, into a single infinite sum:
\be
\eta^{>p} = \sum_{k\geq 0} \etah^{p+2k+1}\quad \in\ \hom(\Om \Mch^s_+,X(\Rch))\ .
\ee
Observe that $\eta^{>p}$ does not vanish on the space of $m$-forms $\Om^m\Mch^s_+$ for $m$ sufficiently large, hence it is only defined on the direct sum $\Om \Mch^s_+=\bigoplus_{m\geq 0}\Om^m\Mch^s_+$ and not on its completion $\Omh \Mch^s_+=\prod_{m\geq 0}\Om^m\Mch^s_+$. Following the terminology of Higson \cite{Hi}, we may call it an \emph{improper} cochain. Its parity is $p+1$ mod 2. Taking the boundary of $\eta^{>p}$ in the $\zz_2$-graded complex $\hom(\Om \Mch^s_+,X(\Rch))$, we get, using the transgression formulas (\ref{tra}):
\be
\chih^p= [\partial, \eta^{>p} ]
\ee 
and it is again a (non-trivial) proper cocycle in $\hom(\Omh \Mch^s_+,X(\Rch))$ representing the bivariant Chern character. This may be depicted through the diagram
$$
\vcenter{\xymatrix@!0@=2.5pc{
\eta^{>p}  \ar[d] &  = &  & \etah^{p+1} \ar[dl] \ar[dr]  & + & \etah^{p+3} \ar[dl] \ar[dr] & + & \etah^{p+5} \ar[dl] \ar[dr] & + & \ldots & \\
[\d,\eta^{>p} ] & = & \chih^{p} & + & 0 & + & 0 & + & 0 & + & \ldots }}
$$
We would like to sum also the components $\etah^{n+1}$ of lower degree, but of course when $n<p$ the product $x_0[F, x_1]\ldots [F,x_{n+1}]$ is no longer in the domain of the supertrace $\tau$, and $\etah^{n+1}$ is not defined. At this point we need to introduce a \emph{renormalization procedure}: suppose we can extend $\tau$ to a linear map $\tau_R$ defined over non-trace-class elements, for example by inserting some regularizing operator. It is customary to call such $\tau_R$ a renormalized supertrace, although the trace property cannot be preserved in general \cite{P4}. In practice the construction of a renormalized supertrace requires to impose some strong conditions on the elements $x_i\in \Mch^s_+$. We shall suppose the existence of a dense subalgebra $\Ac_0\subset\Ac$, such that for any $x_i\in \rho_*\Th\Ac_0$ the \emph{renormalized components of the eta-cochain}
\beq
\lefteqn{\etah^{n+1}_{R0}(x_0\dd x_1\ldots\dd x_{n+1}) =  \frac{\Gamma(\frac{n}{2}+1)}{(n+2)!} \, \frac{1}{2}\tau_R\Big (F x_0[F,x_1]\ldots [F,x_{n+1}]+ } \non \\
&&\qquad \qquad \qquad  \sum_{i=1}^{n+1}(-)^{(n+1)i} [F,x_i]\ldots [F,x_{n+1}] Fx_0 [F,x_1]\ldots [F,x_{i-1}] \Big) \non
\eeq
\beq
\lefteqn{\etah^{n+1}_{R1}(x_0\dd x_1\ldots\dd x_{n+2}) =} \non\\
&& \frac{\Gamma(\frac{n}{2}+1)}{(n+3)!} \sum_{i=1}^{n+2}  \frac{1}{2}\tau_R\nat(ix_0 F + (n+3-i)Fx_0)[F,x_1]\ldots\dd x_i \ldots [F,x_{n+2}]\non
\eeq
are well-defined. The restriction to a dense subalgebra $\Ac_0$ is necessary because $\Ac$ appears generally as a ``big'' completion of $\Ac_0$ and looses some crucial analytic properties. For instance if $\Ac_0$ is an algebra of pseudodifferential operators on a smooth manifold, $\tau_R$ may be constructed by zeta-function renormalization (see section \ref{sren} for concrete examples). Other interesting renormalizations are possible \cite{P6}. Note that $\Ac_0$ is regarded as a discrete algebra. Hence $T\Ac_0$ denotes the \emph{algebraic} tensor algebra of $\Ac_0$, i.e. obtained by taking algebraic tensor products over $\Ac_0$ and finite direct sums, the ideal $J\Ac_0$ is the kernel of the multiplication map $T\Ac_0\to\Ac_0$ and $\Th\Ac_0=\varprojlim_k T\Ac_0/(J\Ac_0)^k$ is the $J\Ac_0$-adic completion of $T\Ac_0$.\\
Now define the full \emph{renormalized eta-cochain} as the improper cochain
\be
\eta_R = \sum_{k< 0} \etah^{p+2k+1}_R + \eta^{>p} \ .
\ee
Hence the difference $\eta_R-\eta^{>p} $ is the finite sum of renormalized components and defines a proper cochain, i.e. vanishes on $x_0\dd x_1\ldots\dd x_n$ for $n$ sufficiently large. The following lemma is obvious.

\begin{lemma}\label{lan}
Let $\Ac_0\subset\Ac$ be a dense subalgebra. For any choice of renormalized supertrace $\tau_R$ the boundary $\chi_R := [\partial , \eta_R]$ is a chain map cohomologous to $\chih^p$ once restricted to the image $\rho_*\Th\Ac_0\subset \Mch^s_+$. Hence the composite $\ch_R(\rho)=\chi_R\rho_*\gamma \in \hom(X(\Th\Ac_0),X(\Rch))$ represents the bivariant periodic Chern character:
\be
\ch_R(\rho)\equiv \ch(\rho)\in HP^{n}(\Ac_0,\Bc)
\ee
in degree $n\equiv p \mod 2$. Moreover if the composite map $(\eta_R-\eta^{>p} )\rho_*\gamma:X(\Th\Ac_0)\to X(\Rch)$ is a cochain of order $\leq n$, then $\ch_R(\rho)$ represents the bivariant Chern character of degree $n$:
\be
\ch_R(\rho)\equiv \ch^n(\rho)\in HC^n(\Ac_0,\Bc)\ .
\ee
\end{lemma}
{\it Proof:} The difference $\ch_R(\rho)-\ch^p(\rho)=[\d,(\eta_R-\eta^{>p} ) \rho_*\gamma]$ is a coboundary in $\hom(X(\Th\Ac_0),X(\Rch))$. \cqfd\\

The way the components of the cocycle $\chi_R$ assemble from the renormalized improper cochain $\eta_R$ may be depicted through the following diagram, using the transgression relation $\chih^n-\chih^{n+2}= [\partial , \etah^{n+1}]$ valid for $n\geq p$:
\be
\vcenter{\xymatrix@!0@=2.5pc{
\eta_R = \ar[d]_{[\d,\ ]} &  \ldots\  + & \etah_R^{p-3} \ar[dl] \ar[dr]  & + & \etah_R^{p-1} \ar[dl] \ar[dr] & + & \etah^{p+1} \ar[dl] \ar[dr] & + & \etah^{p+3} \ar[dl] \ar[dr] & +\  \ldots  \\
\chi_R  = &\ \ldots  \chi_R^{p-4} & + & \chi_R^{p-2} & + & \chi_R^{p} & + & 0 & + & \ 0\  \ldots }} \label{can}
\ee
The components $\chi^n_R$ thus measure to which extent the renormalized eta-cochain fails to be a cocycle. In this sense the equality $\chi_R=[\d,\eta_R]$ should be called an \emph{anomaly}. In favourable circumstances renormalization can be performed in such a way that the anomaly $\chi_R$ is given by a ``local'' formula. This will be explained in the next section. \\

In order to get a closer insight to the link with chiral anomalies let us explain how the renormalized eta-cochain can be used to calculate direct images of $\Kt_1$-classes under quasihomomorphisms of even parity. Hence let $\rho:\Ac\to \Ec^s\triangleright \Ic^s\hotimes\Bc$ be a $(p+1)$-summable quasihomomorphism of even parity, with $p$ even. We suppose that $\Ic$ is the Schatten ideal of $(p+1)$-summable operators in Hilbert space, since it is an important source of examples. The topological part of the Riemann-Roch-Grothendieck theorem \ref{trr} yields a commutative diagram 
\be
\vcenter{\xymatrix{
\Kt_1(\Ic\hotimes\Ac) \ar[r]^{\rho_!} \ar[d] \ar@{.>}[rd]^{\Delta} & \Kt_1(\Ic\hotimes\Bc) \ar[d] \\
HP_1(\Ac) \ar[r]_{\ch(\rho)} & HP_1(\Bc) }}  \label{diag} 
\ee
and one wishes to calculate the diagonal map $\Delta$. From cyclic cohomology invariants of $\Bc$ one thus gets $K$-theory invariants of $\Ac$; this is used for instance in the formulation of higher index theorems \cite{CM90}. A topological $K$-theory class $[g]\in \Kt_1(\Ic\hotimes\Ac)$ is represented by an invertible element $g\in (\Ic\hotimes\Ac)^+$ such that $g-1\in \Ic\hotimes\Ac$ (as usual $^+$ denotes unitalization). The diagonal of (\ref{diag}) thus carries $[g]$ to the periodic Chern character $\ch(\rho_!(g))\in HP_1(\Bc)$. By commutativity of (\ref{diag}), its cyclic homology class is equivalently calculated by the cup-product of $\ch(g)\in HP_1(\Ac)$ with the bivariant Chern character
$$
\ch(\rho_!(g)) \equiv \ch^{2n}(\rho)\cdot\ch(g)
$$
for any choice of degree $2n\geq p$. In order to renormalize we suppose given a dense subalgebra $\Ac_0\subset\Ac$ and take $g$ such that $g-1$ belongs to the algebra of finite size matrices $M_{\infty}(\Ac_0)\subset \Ic\hotimes\Ac$. By \cite{CQ1}, the Chern character $\ch(g)\in HP_1(\Ac_0)$ is represented by the cycle of odd degree in the complex $X(\Th\Ac_0)$
\be
\ch_1(\gh)= \frac{1}{\sqrt{2\pi i}} \, \Tr\nat \gh^{-1}\dd\gh \quad \in \Om^1(\Th\Ac_0)_{\nat}\ ,
\ee
where $\gh$ is any lifting of $g$ to the unitalized tensor algebra $M_{\infty}(\Th\Ac_0)^+$. If one chooses the canonical lifting $\gh=g$ induced by the linear inclusion $\Ac_0\hookrightarrow \Th\Ac_0$, then the image of $\ch_1(\gh)$ by the homotopy equivalence $\gamma: X(\Th\Ac_0)\to \Omh\Th\Ac_0$ is the $(b+B)$-cycle \cite{P5}
$$
\gamma \ch_1(\gh)= \frac{1}{\sqrt{2\pi i}} \sum_{n\geq 0}(-)^n n!\,  \Tr(\gh^{-1}\dd\gh (\dd\gh^{-1}\dd\gh)^n)\ .
$$
The evaluation of this differential form under the chain map $\chih^{2n}\rho_*:\Omh\Th\Ac_0\to X(\Th\Bc)$ represents the cycle $\ch^{2n}(\rho)\cdot\ch(g)$. Let us define $V= \rho_*\gh^{-1}[F,\rho_*\gh]$ as an element of the ideal $\Ic^s\hotimes\Th\Bc \subset \Mch^s$, and the Maurer-Cartan form $\om=\rho_*\uh^{-1}\dd(\rho_*\uh)$. A straightforward computation gives
\be
\ch^{2n}(\rho)\cdot \ch_1(\gh) = \frac{1}{\sqrt{2\pi i}} \frac{(n!)^2}{(2n)!} \, \frac{1}{2}\tau\nat (V^{2n+2}\om-FV^{2n}\dd V) \ \in \Om^1\Th\Bc_{\nat}
\ee
in any degree $2n\geq p$. In fact for $2n>p$ one has $\frac{1}{2}\tau\nat (V^{2n+2}\om-FV^{2n}\dd V) = \tau\nat(V^{2n}\om)$, but this simplification does not occur in the critical degree $2n=p$ because the supertrace $\tau$ is defined only when $V$ is raised to a power $\geq p+1$. Two consecutive degrees are related by the transgression relation 
$$
\ch^{2n}(\rho)\cdot \ch_1(\gh) - \ch^{2n+2}(\rho)\cdot \ch_1(\gh) = \nat\dd \big(\tch^{2n+1}(\rho)\cdot \ch_1(\gh)\big) \ . 
$$
Since the boundary map $\nat\dd: \Th\Bc \to \Om^1\Th\Bc_{\nat}$ factors through the commutator quotient space $\Th\Bc_{\nat}=\Th\Bc/[\ ,\ ]$ it is enough to compute $\tch^{2n+1}(\rho)\cdot \ch_1(\gh)$ modulo commutators and one finds 
\be
\nat\tch^{2n+1}(\rho)\cdot \ch_1(\gh) = \frac{1}{\sqrt{2\pi i}} \frac{(n!)^2}{(2n+1)!} \, \frac{1}{2}\tau\nat (FV^{2n+1})\ \in \Th\Bc_{\nat} \label{nr}
\ee
for $2n\geq p$. In lower degrees the cochains $\tch_R^{2n+1}(\rho)= \etah_R^{2n+1}\rho_*\gamma$ are \emph{defined} using a renormalized supertrace. Therefore the renormalized bivariant Chern character $\ch_R(\rho)=\chi_R\rho_*\gamma$ evaluated on $\ch_1(\gh)$ reads
\be
\ch_R(\rho)\cdot \ch_1(\gh) = \nat\dd \Big( \sum_{n=0}^{p/2-1}\tch_R^{2n+1}(\rho)\cdot \ch_1(\gh) \Big) + \ch^p(\rho)\cdot\ch_1(\gh)\ ,
\ee
and represents the image of $[g]\in \Kt_1(\Ic\hotimes\Ac)$ under the diagonal map of (\ref{diag}), i.e. the Chern character $\ch(\rho_!(g))\in HP_1(\Bc)$. It may be rewritten formally as the boundary of an infinite sum
\be
\ch_R(\rho)\cdot\ch_1(\gh) = \nat\dd \Big( \sum_{n\geq 0}\tch_R^{2n+1}(\rho)\cdot \ch_1(\gh) \Big)\ , \label{infsum}
\ee
provided we organize the cancellation of the infinite sequence of terms $n\gg 0$ as depicted in (\ref{can}). We will see in section \ref{sreg} that it can be done in an alternative way by treating the sum $\sum_n\tch_R^{2n+1}(\rho)\cdot \ch_1(\gh)$ as a formal power series with respect to a ``gauge potential'' $A$, and this series is strictly equivalent to the renormalized partition function of a non-commutative chiral gauge theory. The Chern character $\ch(\rho_!(g))$ corresponds to the associated chiral anomaly and for this reason is always given by a local formula.

\section{Zeta-function renormalization and residues}\label{sren}

Let us illustrate the anomaly formula using Dirac operators and zeta-function renormalization. Theorem \ref{tres} obtained below is a bivariant generalization of the Connes-Moscovici local index formula \cite{CM95}. Consider a $(p+1)$-summable quasihomomorphism $\rho:\Ac\to \Ec^s\triangleright \Ic^s\hotimes\Bc$ of parity $p\mod 2$. For concreteness we suppose that a $\zz_2$-graded Hilbert space $H$ is given and that $\Ec^s$ is realized as a $\zz_2$-graded subalgebra of $\Lc(H)\hotimes\Bc$, where the algebra of bounded operators $\Lc(H)$ is provided with the operator norm. Moreover we take $\Ic^s$ as the Schatten ideal $\Lc^{p+1}(H)$ gifted with its norm of Banach algebra. This situation is not the most general one but contains many interesting examples, see \ref{sequi}. The quasihomomorphism is thus described by a homomorphism $\rho:\Ac\to \Lc(H)\hotimes\Bc$ with $[F,\rho(a)]\in \Lc^{p+1}(H)\hotimes\Bc$ for any $a\in\Ac$. In order to get local formulas for the bivariant Chern character we restrict the homomorphism to a dense subalgebra $\Ac_0\subset\Ac$ and impose some strong regularity conditions on the image $\rho(\Ac_0)$. We shall introduce the notion of abstract pseudodifferential operators associated to such a situation:

\begin{definition}
Let $\rho:\Ac\to \Ec^s\triangleright \Ic^s\hotimes\Bc$ be a quasihomomorphism with $\Ec^s\subset \Lc(H)\hotimes \Bc$ and $\Ic^s=\Lc^{p+1}(H)$ for some $\zz_2$-graded Hilbert space $H$. An \emph{algebra of abstract symbols} associated to a dense subalgebra $\Ac_0\subset \Ac$ is an inductive limit of $\zz_2$-graded Fr\'echet spaces
\be
\Pc=\varinjlim_{\al\in\rr} \Pc_{\al}\ ,
\ee
with continuous linear injections $\Pc_{\al}\to \Pc_{\beta}$ for $\al\leq\beta$, subject to the following conditions: 
\begin{itemize}
\item Each space $\Pc_{\al}$ is represented by unbounded operators on $H$. The unit $1$ and the operator $F$ are both elements of $\Pc_0$, respectively even and odd.
\item $\Pc$ is a $\zz_2$-graded unital filtered algebra, in the sense that for any $\al$ and $\beta$ there is a (jointly) continuous and associative product $\Pc_{\al}\times \Pc_{\beta} \to \Pc_{\al+\beta}$, and $1\in \Pc_0$ is the unit. In particular $\Pc_{\al}$ is a Fr\'echet algebra whenever $\al\leq 0$.
\item The $\zz_2$-graded algebra $\Pc_0$ injects continuously into $\Lc(H)$, and the homomorphism $\rho:\Ac_0\to \Lc(H)\hotimes\Bc$ factors through $\Pc_0\hotimes\Bc$.
\item The commutator $[F,\rho(a)]$ lies in the Fr\'echet algebra $\Pc_{-1}\hotimes\Bc$ if $a\in\Ac_0$.
\item The (unital) filtered algebra $\Pc\hotimes\Bc^+=\varinjlim_{\al}\Pc_{\al}\hotimes\Bc^+$ has a distinguished element of even parity $|D|\in \Pc_1\hotimes\Bc^+$,  whose spectrum is contained in a real interval $[\eps,+\infty)$, $\eps>0$.
\item The resolvent $\la\mapsto (\la-|D|)^{-1}$ is a bounded and holomorphic function over any half-plane $\re(\la)\leq \eps'<\eps$ disjoint from the spectrum of $|D|$, with values in the Fr\'echet algebra $\Pc_{-1}\hotimes\Bc^+$.
\item The multiplier $F\in \Pc_0$ commutes with $|D|$. We call the odd element $D=F|D|\in \Pc_1\hotimes\Bc^+$ a \emph{Dirac operator}.
\end{itemize}
\end{definition}

For $\al\leq 0$ we do not impose that the Fr\'echet topology of $\Pc_{\al}$ is multiplicatively convex. The only serious assumption is about the resolvent.\\
Let $0\to \Jc \to \Rc \to \Bc \to 0$ be a linearly split extension of $m$-algebras. For the moment $\Rc$ is not necessarily quasi-free but we suppose that all the powers $\Jc^k$ are direct summands (as topological vector spaces) in $\Rc$. Using the linear splitting $\si:\Bc\to\Rc$, we may lift $|D|\in \Pc_1\hotimes\Bc^+$ to an element $|\Dh|\in \Pc_1\hotimes\Rc^+$, and then consider $|\Dh|$ as an element of the unital algebra
\be
\Pc \hotimes \Rch^+ = \varinjlim_{\al} \Pc_{\al}\hotimes \Rch^+=\varinjlim_{\al}  \Big( \varprojlim_k \Pc_{\al}\hotimes (\Rc^+/\Jc^k) \Big)\ .
\ee 
Define the lifted Dirac operator (of odd parity) as $\Dh =F|\Dh|\in \Pc_1\hotimes\Rch^+$. Our task is to construct the resolvent of $|\Dh|$ in $\Pc\hotimes\Rch^+$. For $\la\in\cc$ not in the spectrum of $|D|$, the liftings $\si(\la-|D|)=\la-|\Dh|\in \Pc_1\hotimes \Rc^+$ and $\si((\la-|D|)^{-1})\in \Pc_{-1}\hotimes \Rc^+$ are inverse of each other in the algebra $\Pc\hotimes\Rc^+$ modulo the ideal $\Pc\hotimes\Jc$. Hence by the continuous multiplication $\Pc_1\times \Pc_{-1}\to \Pc_0$ one gets
$$
1 - \si(\la-|D|)\si((\la-|D|)^{-1}) \in \Pc_0\hotimes\Jc\ .
$$
Because $\Jc$ is a pro-nilpotent ideal the resolvent of $|\Dh|$ is given by the series 
\be
(\la-|\Dh|)^{-1}=\sum_{k=0}^{\infty} \si((\la-|D|)^{-1})\big(1 - (\la - |\Dh|)\si((\la-|D|)^{-1})\big)^k \label{resolv}
\ee
which makes sense in the pro-algebra $\Pc_{-1}\hotimes\Rch^+$. By hypothesis the function $\la\mapsto (\la-|D|)^{-1}$ is bounded (for all seminorms on $\Pc_{-1}\hotimes\Bc^+$) and holomorphic over any half-plane $\re(\la)\leq \eps'$ disjoint from the spectrum of $|D|$. Consequently, (\ref{resolv}) shows that the function $\la\mapsto (\la-|\Dh|)^{-1}$ is also bounded and holomorphic in each Fr\'echet space $\Pc_{-1}\hotimes (\Rc^+/\Jc^k)$, $k\in \nn$. The complex powers of $|\Dh|$ are thus defined, for any $z\in\cc$ with $\re(z)<-1$, by means of the Riemann integral
\be
|\Dh|^z= \frac{1}{2\pi i}\int \frac{\la^z}{\la-|\Dh|}\, d\la\ \in \Pc_{-1}\hotimes\Rch^+\ , \label{int}
\ee
where we integrate down a vertical line separating zero from the spectrum of $|D|$ as in \cite{Hi}. The condition $\re(z)<-1$ ensures the absolute convergence of the integral. By taking enough products with $|\Dh|\in \Pc_1\hotimes\Rch^+$ we define the complex powers $|\Dh|^z\in \Pc\hotimes\Rch^+$ for any $z\in\cc$. By a classical argument \cite{Se} one has the group law $|\Dh|^{z_1}|\Dh|^{z_2}=|\Dh|^{z_1+z_2}$, which implies that $|\Dh|^z$ belongs to $\Pc_{\al}\hotimes\Rch^+$ where $\al$ is the lowest integer $\geq\re(z)$. For convenience we make the further assumption that $|\Dh|^z\in\Pc_{\re(z)}\hotimes \Rch^+$, and when $\re(z)\ll 0$ it will be used for regularizing the operator trace in order to renormalize the eta-cochain $\etah^{n+1}$ in degrees $n<p$. \\

The bivariant Chern character requires to lift the homomorphism $\rho:\Ac\to \Lc(H)\hotimes\Bc$ to a homomorphism of pro-algebras $\rho_*: \Th\Ac \to \Lc(H)\hotimes \Rch$, sending the ideal $\Jh\Ac$ to $\Lc(H)\hotimes\Jch$. By hypothesis, the restriction of $\rho$ to the dense subalgebra $\Ac_0$ yields a homomorphism $\Ac_0\to \Pc_0\hotimes\Bc$, with $[F,\rho(\Ac_0)]\subset \Pc_{-1}\hotimes\Bc$. Consequently the restriction of $\rho_*$ to the pro-algebra $\Th\Ac_0=\varprojlim_k T\Ac_0/(J\Ac_0)^k$ yields a homomorphism $\Th\Ac_0\to \Pc_0\hotimes\Rch$ with $\Jh\Ac_0\to\Pc_0\hotimes\Jch$ and $[F,\rho_*T\Ac_0]\subset \Pc_{-1}\hotimes\Rch$. Writing $\Dh=F|\Dh|$ we can only show that the commutator subspace $[\Dh,\rho_*\Th\Ac_0]$ is contained in $\Pc_1\hotimes\Rch$, while a better estimate is needed. Hence we shall adapt to the present context the notion of regularity introduced by Connes and Moscovici for spectral triples \cite{CM95}. By assumption the pro-vector space $\Jch$ is a direct summand in $\Rch$, and also the $\Jc$-adic completions $\Jch^k$ of all powers $\Jc^k$. Define the filtration of $\Pc\hotimes \Rch^+$ by the subspaces 
$$
\Fc_{\al}^k= \Pc_{\al}\hotimes \Jch^k\ ,\qquad \Fc_{\al}= \sum_{k\geq 0} \Fc_{\al+k}^k \ ,\qquad \al\in\rr\ , \ k\in\nn
$$
with the convention that $\Jch^0= \Rch^+$ for $k=0$. In $\Fc_{\al}^k$ we call $\al$ the \emph{symbol degree} and $k$ the \emph{adic degree}. One has $\Fc_{\al}^k\cdot \Fc_{\beta}^l\subset\Fc_{\al+\beta}^{k+l}$, $|\Dh|^z\in \Fc_{\re(z)}^0$, $\rho_*\Th\Ac_0\subset \Fc_0^0$, $[F,\rho_*\Th\Ac_0]\subset \Fc_{-1}^0$, $\rho_*(\Jh\Ac_0)^k\subset \Fc_0^k$ and $[F,\rho_*(\Jh\Ac_0)^k]\subset \Fc_{-1}^k$. Consider the derivation $\delta=[|\Dh|,\ ]$ on the algebra $\Pc\hotimes\Rch$. Since $|\Dh|\in \Fc_1^0$, the commutator $\delta$ a priori defines a linear map from $\Fc_{\al+k}^k$ to $\Fc_{\al+k+1}^k$ for any $\al$ and $k$. The regularity condition basically requires that all the derivatives of $\Th\Ac_0$ remain in $\Fc_0$:

\begin{definition}\label{dreg}
Let $0\to \Jc\to \Rc\to \Bc\to 0$ be an extension of $m$-algebras, and $\rho:\Ac\to \Ec^s\triangleright \Ic^s\hotimes\Bc$ a quasihomomorphism gifted with an algebra of abstract symbols $\Pc$ associated to a dense subalgebra $\Ac_0\subset \Ac$. We say that $\rho$ is \emph{$\Rc$-regular} if the powers of the derivation $\delta=[|\Dh|, \ ]$ fulfill
\beq
\delta^n(\rho_*\Th\Ac_0)+ \delta^n [\Dh,\rho_*\Th\Ac_0] &\subset& \Fc_0^0 + \Fc_1^1 + \ldots + \Fc_n^n\ \subset \ \Fc_0 \ ,\non\\
\delta^n(\rho_*(\Jh\Ac_0)^k)+ \delta^n [\Dh,\rho_*(\Jh\Ac_0)^k] &\subset& \Fc_0^k + \Fc_1^{k+1} + \ldots + \Fc_n^{k+n}\ \subset \ \Fc_{-k} \ .\non
\eeq
\end{definition}

Hence roughly speaking a quasihomomorphism is regular if the commutators with $|\Dh|$ can only increase simultaneously the symbol degree and the adic degree. In this sense the images of $\Th\Ac_0$ and $\Jh\Ac_0$ are ``smooth'' in $\Pc\hotimes\Rch$. In particular if $\Bc=\Rc=\cc$ one recovers the Connes-Moscovici regularity condition that $\delta^n\rho(\Ac_0)$ and $\delta^n[D,\rho(\Ac_0)]$ are contained in some algebra of abstract symbols $\Pc_0\subset \Lc(H)$. In the general case regularity strongly depends on the choice of extension $\Rc$; for example the tensor algebra $\Rc=T\Bc$ leads to a regular quasihomomorphism only in very special cases.

\begin{example}
\textup{ Definition \ref{dreg} is motivated by the classical example of families of pseudodifferential operators parametrized by a smooth manifold $M$. Define $\Bc$ as the commutative Fr\'echet $m$-algebra $\cinf(M)$, possibly subject to some decay conditions at infinity. Its completed projective tensor product with the algebra of abstract symbols $\Pc$ is thus an algebra of smooth families of symbols parametrized by $M$. Suppose now that the product of symbols fulfills the classical property that the commutator $[\Pc_{\al},\Pc_{\beta}]$ is included in $\Pc_{\al+\beta-1}$. In other words, the product is commutative at the level of leading symbols. Hence if $|D|\in \Pc_1\hotimes\Bc^+$ is a family of order one, the commutativity of $\Bc$ implies 
$$
[|D|,\Pc_{\al}\hotimes\Bc^+]\subset \Pc_{\al}\hotimes\Bc^+\qquad \forall \al\in\rr\ .
$$
Now let $\Om^+(M)$ be the space of differential forms of even degree over $M$. We deform the ordinary commutative product of differential forms into the associative and non-commutative product defined by Fedosov \cite{F}
$$
x\circ y = xy - dxdy\qquad \forall x,y\in \Om^+(M)\ ,
$$
and denote by $\Rc$ the non-commutative algebra thus obtained. The projection of differential forms onto the degree zero subspace $\cinf(M)$ yields an extension of Fr\'echet $m$-algebras $0\to\Jc\to\Rc\to\Bc\to 0$, where $\Jc$ is the ideal of differential forms of even degree $\geq 2$. Because $\Jc$ is nilpotent, one has $\Rch=\Rc$ and $\Jch=\Jc$. The canonical linear injection of $\Bc$ into $\Rc$ allows to lift $|D|$ to an element $|\Dh|\in \Pc_1\hotimes\Rc^+$. Now the derivation defined by the \emph{Fedosov} commutator $\delta=[|\Dh|,\ ]_{\circ}$ does no longer preserve the symbol degree since $\Rc$ is not commutative. Instead, for any element $x\in \Pc_{\al}\hotimes\Jc^k=\Fc_{\al}^k$ one has
$$
\delta(x)= [|\Dh|,x]_{\circ} = [|\Dh|,x] - d|\Dh| dx + dx d|\Dh|\ .
$$
Here $[|\Dh|,x]$ is the commutator involving the ordinary (commutative) product of forms, hence remains in $\Fc_{\al}^k$. This shows that $\delta\Fc_{\al}^k\subset \Fc_{\al}^k + \Fc_{\al+1}^{k+1}$ and by iteration one gets the filtration property
$$
\delta^n \Fc_{\al}^k\subset \Fc_{\al}^k + \Fc_{\al+1}^{k+1} + \ldots + \Fc_{\al+n}^{k+n} \subset \Fc_{\al -k}\ .
$$
This fails when $\Rc$ is replaced by the tensor algebra over $\Bc$, unless $|\Dh|\in \Pc_1\otimes 1$ is a constant symbol over $M$ so that $d|\Dh|=0$. Note that because $\Rc$ is not quasi-free the homology of $X(\Rc)$ is not isomorphic but contains the de Rham cohomology of $M$ \cite{CQ1}. For a general quasihomomorphism we only need to demand that the filtration property above holds on the image of $\Th\Ac_0$ in $\Pc\hotimes \Rch$. }
\end{example}

Now we return to the general case and assume that the quasihomomorphism is regular. Define the algebra of abstract pseudodifferential operators $\Psi\Ac_0$ as the (non-closed) subalgebra of $\Pc\hotimes\Rch$ generated by $\rho_*\Th\Ac_0$, $F$ and all the complex powers $|\Dh|^z$, and such that its elements contain at least one factor in $\rho_*\Th\Ac_0$. Any abstract pseudodifferential operator is a linear combination of elements $x\in \Psi\Ac_0$ having an asymptotic expansion
\be
x\simeq \sum_{k\geq 0}(a_k+b_k F)\, |\Dh|^{z-k}\ ,\quad z\in\cc\ ,\label{asy}
\ee
where $a_k$ and $b_k$ are in the subalgebra of $\Fc_0$ generated by all the derivatives $\delta^n(\rho_*\Th\Ac_0)$ and $\delta^n [\Dh,\rho_*\Th\Ac_0]$, and $\simeq$ means that the difference $x-\sum_{k= 0}^{N}(a_k+b_k F)\, |\Dh|^{z-k}$ lies in $\Fc_{\re(z)-N-1}$ for any $N\in \nn$. The fact that any element has an asymptotic expansion is a consequence of the formula (whose proof is identical to \cite{Hi} Lemma 4.20)
\be
a_1|\Dh|^{z_1} \cdot a_2|\Dh|^{z_2} \simeq \sum_{k\geq 0} \binom{z_1}{k} a_1 \delta^k(a_2) |\Dh|^{z_1+z_2-k}\ ,
\ee
where $\binom{z}{k}$ is the binomial coefficient $\frac{z(z-1)\ldots(z-k+1)}{k!}$. The algebra $\Psi\Ac_0$ is filtered by the subspaces $(\Psi\Ac_0)_{\al}^k=\Psi\Ac_0\cap \Fc_{\al}^k$; one has $x\in (\Psi\Ac_0)_{\al}^k$ if and only if $x|\Dh|^{-\al} \in (\Psi\Ac_0)_0^k$. In particular, for any $x\in \rho_*\Th\Ac_0$ the commutator $[F,x]$ belongs to $(\Psi\Ac_0)_{-1}^0$. More generally for any elements $x_i \in \rho_*(\Jh\Ac_0)^{k_i}$ the product $x_0[F,x_1]\ldots [F,x_n]$ is a pseudodifferential operator in $(\Psi\Ac_0)_{-n}^k$ with $k=\sum_i k_i$. \\
Similarly we would like to consider a one-form $x_0[F,x_1]\ldots \dd x_i \ldots [F,x_{n+1}]$ as a kind of pseudodifferential operator in an appropriate subspace of the universal $\Pc\hotimes\Rch^+$-bimodule
\be
\Om^1(\Pc\hotimes\Rch) := \varinjlim_{\al,\beta}\, (\Pc_{\al}\hotimes\Rch^+) \dd(\Pc_{\beta}\hotimes \Rch)\ .
\ee
Note that here $\Rch^+$ is view as a \emph{unital} algebra, so the universal derivation $\dd: \Pc\hotimes\Rch^+\to \Om^1(\Pc\hotimes\Rch)$ vanishes on the subspace $\Pc\hotimes 1$ by definition. Since $\Pc$ is a $\zz_2$-graded algebra, $\dd$ actually anticommutes with odd elements and the Leibniz rule reads $\dd(xy)=\dd x\,y +(-)^{|x|}x\dd y$. Moreover the Hochschild operator $b:\Om^1(\Pc\hotimes\Rch)\to \Pc\hotimes\Rch$ is the graded commutator $\bb( x\dd y)=(-)^{|x|}[x,y]$. Then, filter $\Om^1(\Pc\hotimes\Rch)$ by the subspaces 
$$
(\Om^1\Fc)_{\al}^k=\sum_{\substack{\beta+\gamma=\al \\ l+m=k}} \Fc_{\beta}^l\dd\Fc_{\gamma}^m\ ,\qquad  (\Om^1\Fc)_{\al}=\sum_{k\geq 0} (\Om^1\Fc)_{\al+k}^k 
$$
for any $\al\in\rr$ and $k\in\nn$. We denote improperly by $\Om^1\Psi\Ac_0\subset \Om^1(\Pc\hotimes\Rch)$ the (non-closed) subspace of abstract pseudodifferential operators, generated by products of $\rho_*\Th\Ac_0$, $F$, $|\Dh|^z$ together with their $\dd$-derivatives, and such that at least one factor comes from $\rho_*\Th\Ac_0$. Let $\nat$ be the projection of $\Om^1(\Pc\hotimes\Rc)$ onto its quotient by the commutator subspace $[\Pc\hotimes\Rc^+, \Om^1(\Pc\hotimes\Rc)]$. Any abstract pseudodifferential operator in $\nat \Om^1\Psi\Ac_0$ is a linear combination of one-forms $\om$ having an asymptotic expansion
\be
\om \simeq \sum_{k\geq 0}\nat\big( (a_k+Fb_k)|\Dh|^{z-k}\dd c_k +(a'_k+Fb'_k)|\Dh|^{z-k-1}\dd |\Dh| \big)\ ,
\ee
where as before $a_k$, $b_k$, $a'_k$, $b'_k$ and $c_k$ are in the subalgebra of $\Fc_0$ generated by the derivatives $\delta^n(\rho_*\Th\Ac)$, $\delta^n [\Dh,\rho_*\Th\Ac]$, and $\simeq$ means that the difference $\om-\sum_{k=0}^{N}\big((a_k+b_k F)\, |\Dh|^{z-k}\dd c_k+\ldots\big)$ lies in $\nat(\Om^1\Fc)_{\re(z)-N-1}$ for any $N\in \nn$. Indeed, the integral formula (\ref{int}) and the identity $\dd (\la-|\Dh|)^{-1}=(\la-|\Dh|)^{-1}\dd |\Dh|(\la-|\Dh|)^{-1}$ yield the asymptotic expansion
\be
\dd |\Dh|^z \simeq \sum_{k\geq 1} \binom{z}{k} (-)^{k-1}|\Dh|^{z-k}\delta^{k-1}(\dd|\Dh|)\ ,
\ee
hence for any $a$ one gets $\nat\dd(a|\Dh|^z) \simeq \nat|\Dh|^z\dd a+ \sum_{k\geq 1} \left( 
\begin{smallmatrix}
z \\
 k \end{smallmatrix} \right)\nat \delta^{k-1}(a) |\Dh|^{z-k}\dd|\Dh|$. This shows that any element of $\nat\Om^1\Psi\Ac_0$ has an asymptotic expansion. Since $\Psi\Ac_0$ is an algebra, the Hochschild boundary $\bb$ clearly sends $\nat\Om^1\Psi\Ac_0$ to $\Psi\Ac_0$, and $\nat\dd$ maps $\Psi\Ac_0$ to $\nat\Om^1\Psi\Ac_0$ by construction. Moreover $\Om^1\Psi\Ac_0$ is filtered by the subspaces $(\Om^1\Psi\Ac_0)_{\al}^k:=\Om^1\Psi\Ac_0\cap (\Om^1\Fc)_{\al}^k$, and the obvious compatibility of the filtrations with the boundary maps $\bb$ and $\nat\dd$ leads to a family of complexes
\be
X(\Psi\Ac_0)_{\al}^k\ :\ (\Psi\Ac_0)_{\al}^k\ \rightleftarrows\ \nat(\Om^1\Psi\Ac_0)_{\al}^k
\ee
indexed by $\al\in\rr$ and $k\in\nn$. In particular for any elements $x_i\in \rho_*(\Jh\Ac_0)^{k_i}$, the one-form $\nat x_0[F,x_1]\ldots \dd x_i \ldots [F,x_{n+1}]$ lies in $\nat(\Om^1\Psi\Ac_0)_{-n}^k$ with $k=\sum_i k_i$. \\

Before renormalizing the eta-cochain with our pseudodifferential calculus, it remains to discuss partial supertraces. By definition the operator supertrace $\tau$ is defined on the $(p+1)$-th power of the Schatten ideal $\Ic^s=\Lc^{p+1}(H)$ (see \cite{P5} for the normalization of $\tau$ depending on the parity of the quasihomomorphism). Since $\Pc$ is represented by unbounded operators on $H$ we can also regard $\tau$ as an unbounded linear map on each $\Pc_{\al}$. In general its domain cannot be extended to the entire space $\Pc_{\al}$ even for $\al\ll 0$ (for example when $\Pc_{\al}$ is the space of pseudodifferential operators of order $\leq \al$ over a non-compact manifold); however it is reasonable to impose that the unbounded partial supertraces 
\beq
\tau &:& \Fc_{\al}^k=\Pc_{\al}\hotimes\Jch^k \to \Jch^k \label{part}\\
\tau\nat &:& (\Om^1\Fc)_{\al}^k \to \sum_{l+m=k} \Pc_{\al}\hotimes \Jch^l\dd\Jch^m \to \nat(\Jch^k\dd\Rch+\Jch^{k-1}\dd\Jch) \non
\eeq
extend to a well-defined chain map from the entire subcomplex $X(\Psi\Ac_0)_{\al}^k$ to $X(\Rch)$ whenever $\al\ll 0$. The idea is that a pseudodifferential operator of the kind $x |\Dh|^{\al}$ with $x\in\rho_*T\Ac_0$ should be trace-class. Recall that $X(\Rch)$ is filtered by the subcomplexes $F^{2k-1}\Xh(\Rc,\Jc): \Jch^{k}\rightleftarrows \nat(\Jch^k\dd\Rch+\Jch^{k-1}\dd\Jch)$.

\begin{definition}\label{dsd}
A regular quasihomomorphism has \emph{finite analytic dimension} $d\in\rr$ if the subspaces of pseudodifferential operators $(\Psi\Ac_0)_{\al}^k$ and $(\Om^1\Psi\Ac_0)_{\al}^k$ are in the domain of the unbounded partial supertraces (\ref{part}) whenever $\al<-d$, and the latter extend to a chain map preserving the adic degree:
$$
 X(\Psi\Ac_0)_{\al}^k \to F^{2k-1}\Xh(\Rc,\Jc)\ .
$$
Moreover the quasihomomorphism has \emph{discrete dimension spectrum} if there exists a discrete subset $Sd\subset\cc$, such that for any $x,y\in (\Psi\Ac_0)_0^*\cup \{1\}$ not simultaneously equal to 1, the zeta-functions with values in $X(\Rch)$
$$
\tau(x|\Dh|^{z}y)\in \Rch\ ,\qquad \left. \begin{array}{c}
\tau\nat(x|\Dh|^{z}\dd y) \\
 \tau\nat(x|\Dh|^{z-1}y\dd |\Dh|) \end{array} \right\} \in \Om^1\Rch_{\nat}
$$
are holomorphic over the half-plane $\re(z) <-d$ and extend to meromorphic functions over $\cc$ with poles contained in $Sd$, preserving the adic degree.
\end{definition}

On the technical side we have to assume that whenever a pseudodifferential operator constructed from $x,y\in \Psi\Ac_0$ and the complex power $|\Dh|^z$ is in the domain of the trace, the corresponding zeta-function is holomorphic and coincides with the value given by the trace. This is to avoid pathological situations in which the zeta-function would differ from the value of the trace even when the latter is defined. Note that speaking about holomorphic functions with values in $X(\Rch)$ requires that the odd part of the $X$-complex $\Om^1\Rch_{\nat}$ can be written as an inverse system of Fr\'echet spaces; this holds true for example when $\Rc$ is quasi-free and we shall not insist on that point.\\ 
In the following we assume that the analytic dimension $d$ is exactly the degree $p$ of the quasihomomorphism (morally $\rho(a)|D|^{-1}$ sits in $\Lc^{p+1}(H)\hotimes\Bc$ for any $a\in\Ac_0$), but other situations are conceivable. Now we can define a renormalized partial (super)trace $\tau_R$ by taking the \emph{finite part} at $z=0$ of the meromorphic extension of the trace $\tau$, regularized by $|\Dh|^{-z}$ when $\re(z)\gg 0$:
\be
\tau_R(x) = \Pf_{z=0} \tau(x|\Dh|^{-z}) \ ,\qquad \tau_R\nat(x\dd y) = \Pf_{z=0} \tau\nat(x|\Dh|^{-z}\dd y)\ ,
\ee
for $x,y\in \Psi\Ac_0$. The finite part of a meromorphic function $\Pf_{z=0}\zeta(z)$ is the coefficient $\zeta_0$ of its Laurent expansion $\zeta(z)=\sum_k \zeta_k z^k$ around zero. Hence by hypothesis, when the symbol degrees of $x,y$ are low enough the zeta-functions are holomorphic around zero and their values at $z=0$ coincide with the unrenormalized traces $\tau(x)$ and $\tau\nat(x\dd y)$. However, the renormalized trace is no longer a trace: because of the presence of the operator $|\Dh|^{-z}$, the naive algebraic identity (the sign $\pm$ depends on the parity of $x$ and $\tau$)
$$
\bb \, \tau\nat(x\dd y) = \pm\tau([x,y])
$$
is broken by \emph{anomalous terms}. The latter are necessarily residues of zeta-functions and can be explicitly computed using the above pseudodifferential calculus. One has
$$
\bb \, \tau_R\nat(x\dd y)= \pm \Pf_{z=0} \tau([x|\Dh|^{-z},y]) = \pm \tau_R([x,y]) \pm \Pf_{z=0} \tau(x[|\Dh|^{-z},y]) \ .
$$
For concreteness take $y\in (\Psi\Ac_0)_{\al}^l$. The commutator $[|\Dh|^{-z},y]$ admits the asymptotic expansion
$$
[|\Dh|^{-z},y] \simeq \sum_{k\geq 1} \binom{-z}{k} \delta^k(y) |\Dh|^{-z-k} \simeq  z\sum_{k\geq 1} \frac{(-)^k}{k!} \frac{\Gamma(z+k)}{\Gamma(z+1)}\, \delta^k(y) |\Dh|^{-z-k} \ ,
$$
where $\delta^k(y) |\Dh|^{-z-k} \in (\Psi\Ac_0)_{\al-z-k}^{l} + (\Psi\Ac_0)_{\al-z-k+1}^{l+1} + \ldots + (\Psi\Ac_0)_{\al-z}^{l+k}$. It means that if we work modulo high adic degrees, the difference $x[|\Dh|^{-z},y] - \sum_{k= 1}^N \binom{-z}{k} x\delta^k(y) |\Dh|^{-z-k}$ becomes trace-class for large $N$ and its trace, being holomorphic around zero, vanishes at $z=0$. Hence from $\Pf\limits_{z=0} z\zeta(z)=\res \zeta(z)$ we can write the anomalous terms as
$$
\Pf_{z=0} \tau(x[|\Dh|^{-z},y]) = \sum_{k\geq 1} \frac{(-)^k}{k!} \res \frac{\Gamma(z+k)}{\Gamma(z+1)}\,\tau(x \delta^k(y) |\Dh|^{-z-k})\ ,
$$
where the sum over $k$ in the right-hand-side is pro-finite in $\Rch$. The components of the renormalized eta-cochain $\etah_R^{n+1}$ in low degrees $n<d=p$ are defined in exactly the same way, with insertion of $|\Dh|^{-z}$ at appropriate places:
\beq
\lefteqn{\etah^{n+1}_{R0}(x_0\dd x_1\ldots\dd x_{n+1}) =\frac{\Gamma(\frac{n}{2}+1)}{(n+2)!}\,  \frac{1}{2} \Pf_{z=0}\tau\Big( |\Dh|^{-z}F x_0[F,x_1]\ldots [F,x_{n+1}]+ } \non \\
&&\qquad\qquad    \sum_{i=1}^{n+1}(-)^{(n+1)i} [F,x_i]\ldots [F,x_{n+1}] |\Dh|^{-z}Fx_0 [F,x_1]\ldots [F,x_{i-1}] \Big) \non
\eeq
\beq
\lefteqn{\etah^{n+1}_{R1}(x_0\dd x_1\ldots\dd x_{n+2}) =\frac{\Gamma(\frac{n}{2}+1)}{(n+3)!}\frac{1}{2}\times   }\label{ren}\\
&&\sum_{i=1}^{n+2}  \Pf_{z=0} \Big( \sum_{j=1}^i \tau\nat x_0F[F,x_1]\ldots [F,x_{j-1}]|\Dh|^{-z}[F,x_j] \ldots \dd x_i \ldots [F,x_{n+2}] \non\\
&&\qquad  +\sum_{j=i}^{n+2} \tau\nat Fx_0[F,x_1]\ldots \dd x_i \ldots [F,x_j]|\Dh|^{-z}[F,x_{j+1}] \ldots [F,x_{n+2}] \Big) \ ,\non
\eeq
for any $x_i\in \rho_*\Th\Ac_0$. Note that for $n\geq d=p$ the renormalized eta-cochain $\etah^{n+1}_R$ coincides with $\etah^{n+1}$ because in this situation the zeta-functions are holomorphic at $z=0$ by the summability hypothesis $[F,x]\in (\Psi\Ac_0)_{-1}^0$. According to Lemma \ref{lan}, when $\Rc$ is quasi-free the cocycle 
\be
[\partial,\eta_R]\rho_*\gamma \ \in\hom(X(\Th\Ac_0),X(\Rch))
\ee
represents the periodic Chern character $\ch(\rho)\in HP^{p}(\Ac_0,\Bc)$. In fact a finer analysis of the adic properties of the renormalized eta-cochain performed in the theorem below shows that the cocycle $[\partial,\eta_R]\rho_*\gamma$ is of order $p$ hence represents also the non-periodic bivariant Chern character $\ch^p(\rho)\in HC^p(\Ac_0,\Bc)$, and carries all the information about secondary classes. The general discussion of section \ref{sbiv} shows that the components of $[\partial,\eta_R]$, depicted through the diagram (\ref{can}), arise as an anomalous boundary of the renormalized eta-cochain caused by the insertion of $|\Dh|^{-z}$. They are explicitly computable as residues of zeta-functions, by a straightforward generalization of the argument given above for the renormalized trace. One thus obtains a \emph{local} representative (in the sense of noncommutative geometry) of the bivariant Chern character. \\
It should be stressed that the renormalization (\ref{ren}) is by no means unique; one could indeed use the operator $|\Dh|^{-z}$ in many different ways. The simplest modification is, for example, to multiply the zeta-function appearing in the definition $\etah_R^{n+1}=\Pf_{z=0}\zeta(z)$ by a function $h(z)$ holomorphic around zero, with $h(0)=1$. Since one has $h(z)=1+z\, l(z)$ for some $l(z)$, the new renormalized cochain reads
$$
\etah_{R'}^{n+1}=\etah_R^{n+1} + \res\ l(z)\zeta(z)\ .
$$
The function $l(z)$ being holomorphic, the residue term only extracts the \emph{poles} of $\zeta(z)$ at zero. This illustrates a general fact: different renormalizations $\eta_R$ and $\eta_{R'}$ always differ by a (pro-finite) sum of \emph{local counterterms} (here, residues of zeta-functions), and their boundaries $[\partial,\eta_R]\rho_* \gamma$ and $[\partial,\eta_{R'}]\rho_* \gamma$ represent the same bivariant periodic cohomology class $\ch(\rho)\in HP^p(\Ac_0,\Bc)$. We could also modify $\eta_R$ by adding more complicated terms; they would be typically of the form
\be
\res\ h(z) \tau(\delta^{k_0}X_0\ldots \delta^{k_n}X_n\, |\Dh|^{-z-m})\ , \label{contre}
\ee
where $h(z)$ is holomorphic around zero, $\delta^k$ is the $k$-th power of the derivation $[|\Dh|,\ ]$, $m$ is a positive integer, and $X_i$ denotes either $x_i$, $\dd x_i$, $[F,x_i]$ or $[\Dh,x_i]$. It is however not true that any renormalization leads to a representative of the non-periodic bivariant class $\ch^p(\rho)\in HC^p(\Ac_0,\Bc)$; in general we can only control the periodic class. A judicious choice of counterterms makes the link between the renormalization (\ref{ren}) and the bivariant generalization of the Connes-Moscovici local index formula \cite{CM95}:

\begin{theorem}[Anomaly formula]\label{tres}
Let $\rho:\Ac\to \Ec^s\triangleright \Ic^s\hotimes\Bc$ be a $(p+1)$-summable quasihomomorphism of parity $p\mod 2$, gifted with an algebra of abstract symbols associated to a dense subalgebra $\Ac_0\subset\Ac$. Suppose that $\rho$ is regular with respect to an extension $\Rc$ of $\Bc$, has finite analytic dimension $d=p$ and discrete dimension spectrum. Then the boundary of the eta-cochain renormalized by (\ref{ren}) is a cocycle of order $p$ in the bivariant complex
$$
[\partial,\eta_R]\rho_*\gamma \ \in \hom^p(X(\Th\Ac_0),X(\Rch))
$$
Hence it represents the non-periodic Chern character $\ch^p(\rho)\in HC^p(\Ac_0,\Bc)$ when $\Rc$ is quasi-free. Moreover $[\partial,\eta_R]\rho_*$ is cohomologous, in the complex $\hom(\Omh\Th\Ac_0, X(\Rch))$, to the cocycle $\chi_R\rho_*$ whose component in any degree $n\equiv p\mod 2$ is given by a linear map $\chi^n_{R0}\rho_*:\Omh^n\Th\Ac_0\to \Rch$ 
\beq
\lefteqn{\chi_{R0}^{n}(x_0\dd x_1\ldots\dd x_{n}) =\sum_{k_0,\ldots,k_{n}\geq 0} (-)^{k+n}c(k_0,\ldots,k_{n})\, \res \Big(\frac{\Gamma(z+k+\frac{n}{2})}{\Gamma(z+1)}\times}\non\\ 
&&\qquad\qquad  \sum_{i=0}^{n} (-)^{i(n-1)}\, \tau(dx_{n-i+1}^{(k_0)}\ldots x_0^{(k_i)}dx_1^{(k_{i+1})}\ldots dx_{n-i}^{(k_{n})}|\Dh|^{-2(z+k)-n})\Big)\non\\
&& \label{chir0}
\eeq
and a linear map $\chi^n_{R1}\rho_*:\Omh^{n+1}\Th\Ac_0\to \Om^1\Rch_{\nat}$
\beq
\lefteqn{\chi_{R1}^{n}(x_0\dd x_1\ldots\dd x_{n+1}) =\sum_{k_0,\ldots,k_{n}\geq 0} (-)^{k+n}c(k_0,\ldots,k_{n})\, \res \Big(\frac{\Gamma(z+k+\frac{n}{2})}{\Gamma(z+1)}\times}\non\\
&&\qquad \sum_{i=0}^{n}(-)^{in} \tau\nat (dx_{n-i+2}^{(k_0)}\ldots x_0^{(k_i)}dx_1^{(k_{i+1})}\ldots dx_{n-i}^{(k_{n})}|\Dh|^{-2(z+k)-n}\dd x_{n-i+1})\Big)\non\\
&&\qquad\qquad - \sum_{k_0,\ldots,k_{n+1}\geq 0} (-)^{k+n}c(k_0,\ldots,k_{n+1})\, \res \Big(\frac{\Gamma(z+k+\frac{n}{2}+1)}{\Gamma(z+1)}\times\non\\
&&\qquad \sum_{i=0}^{n+1}(-)^{in} \tau\nat (dx_{n-i+2}^{(k_0)}\ldots x_0^{(k_i)}dx_1^{(k_{i+1})}\ldots dx_{n-i+1}^{(k_{n+1})}|\Dh|^{-2(z+k+1)-n}\dd \Dh)\Big)\non\\
&& \label{chir1}
\eeq
where $k$ denotes the sum $\sum_i k_i$, $dx$ is the commutator $[\Dh,x]$, the superscript $x^{(k_i)}$ denotes the $k_i$-th power of the derivation $[\Dh^2,\ ]$ acting on $x$, and the constant $c(k_0,\ldots,k_{n})$ is defined by
$$
c(k_0,\ldots,k_{n})^{-1}=k_0!\ldots k_n! (k_0+1)(k_0+k_1+2)\ldots(k_0+\ldots+k_n+n+1)\ .
$$
Hence the composite $\chi_R\rho_*\gamma \in \hom(X(\Th\Ac_0),X(\Rch))$ represents the periodic Chern character $\ch(\rho)\in HP^p(\Ac_0,\Bc)$ when $\Rc$ is quasi-free. 
\end{theorem}
Before giving the proof let us make some comments about the cocycle $\chi_{R}$. Firstly, if $k+\frac{n}{2} > 0$ the quotient $\frac{\Gamma(z+k+\frac{n}{2})}{\Gamma(z+1)}$ is holomorphic at $z=0$ (it is moreover a polynomial when $n$ is even). Hence the residue term extracts only the poles of the corresponding zeta-function. Only when $n=0$, $k=0$ the above quotient of gamma-functions is singular and equals $1/z$. In this case the residue actually selects the finite part of the zeta-function. This happens when the quasihomomorphism has even degree, and concerns only the components of degree zero
\beq
\chi_{R0}^0(x_0) &=& \Pf_{z=0}\tau(x_0|\Dh|^{-2z}) + \ \mbox{poles (terms $k>0$)}\ ,\non\\
\chi_{R1}^0(x_0\dd x_1) &=& \Pf_{z=0}\tau\nat(x_0|\Dh|^{-2z}\dd x_1) + \ \mbox{poles (terms $k>0$)}\ .\non
\eeq
Secondly, denote by $(\Psi\Ac_0)_{\al}$ the sum $\sum_{l\geq 0}(\Psi\Ac_0)_{\al+l}^l$. For $x\in (\Psi\Ac_0)_0$ the commutator $x^{(1)}=[|\Dh|^2,x]=\delta^2 (x)+2\delta (x)\, |\Dh|$ is a pseudodifferential operator in $(\Psi\Ac_0)_{1}$, and inductively one has $x^{(k)}\in (\Psi\Ac_0)_{k}$ and similarly for $dx^{(k)}$. Using Definition \ref{dsd} of the analytic dimension, it follows that for any $x_i\in \rho_*\Th\Ac_0$ the pseudodifferential operator
$$
x_0^{(k_0)} dx_1^{(k_1)}\ldots dx_n^{(k_n)}|\Dh|^{-2k-n} \in (\Psi\Ac_0)_{-k-n}\ ,\qquad k_0+\ldots+ k_n=k
$$
becomes trace-class for $n+k\gg 0$, provided we mod out the large adic degrees. It follows that the evaluation of the chain map $\chi_R\rho_*\gamma$ on an arbitrary element of the $X$-complex $X(\Th\Ac_0)$ is actually given by a pro-finite sum of residues in $X(\Rch)$. \\

\noindent {\it Proof of \ref{tres}:} According to Lemma \ref{lan} we have to show that $(\eta_R-\eta^{>p})\rho_*\gamma$ is a cochain of order $p$ in the complex $\hom(X(\Th\Ac_0),X(\Rch))$. Since the difference $\eta_R-\eta^{>p}$ is the finite sum of the renormalized components $\etah_R^{n+1}$ for $n=p-2$, $p-4$, etc... it is enough to establish that the latter are cochains of order $\leq p$. We follow the proof of \cite{P5} Proposition 3.9 and introduce the coarse filtration of $\Om^nT\Ac_0$ by the subspaces
$$
H^k\Om^n T\Ac_0 = \sum_{k_0+\ldots + k_n\geq k} (J\Ac_0)^{k_0} \dd(J\Ac_0)^{k_1} \ldots \dd (J\Ac_0)^{k_n}\ ,
$$
with the convention $(J\Ac_0)^0=(T\Ac_0)^+$. Morally $H^k\Om^nT\Ac_0$ contains at least $k$ powers of the ideal $J\Ac_0$. The components of the renormalized eta-cochain are products of commutators $[F,x]$ with $|\Dh|^{-z}$ evaluated on the partial trace $\tau$, which preserves the adic degree $k$ by virtue of Definition \ref{dsd}, hence (we do no longer mention the homomorphism $\rho_*$ when it is not necessary)
\beq
\etah^{n+1}_{R0}(H^k\Om^{n+1}T\Ac_0) &\subset& \Jch^k\ , \non\\
\etah^{n+1}_{R1}(H^k\Om^{n+2}T\Ac_0) &\subset& \nat(\Jch^k\dd\Rch + \Jch^{k-1}\dd\Jch)\ . \non
\eeq
These estimates are sufficient to prove that the composite $\etah_R^{n+1}\rho_*\gamma$ is a cochain of order $n+3$. This is what we want except for the component of highest degree $\etah_R^{p-1}$, for which we need a better estimate. As explained in \cite{P5} Proposition 3.9, an optimal estimate requires to refine the $H$-filtration by introducing a new one,
$$
G^k\Om^nT\Ac_0 = \mspace{-10mu} \sum_{k_0+\ldots +k_n\geq k} \mspace{-10mu} (J\Ac_0)^{k_0} \dd T\Ac_0\ldots (J\Ac_0)^{k_{n-1}}\dd T\Ac_0 (J\Ac_0)^{k_n} + H^{k+1}\Om^nT\Ac_0
$$
and relate the Hochschild coboundary of $\etah_{R1}^{p-1}$ to the component $\chih^p_1$ of the Chern character in critical degree $p$. The crucial properties of $\chih_1^p$ which make the proof of \cite{P5} work are: i) invariance under the Karoubi operator $\kappa$ and ii) compatibility with the $G$-filtration as follows:
$$
\chih_1^p(G^k\Om^{p+1}T\Ac_0)\subset \nat \Jc^k\dd\Rc\ .
$$
In the present case one must be careful because the eta-cochain of degree $p-1$ is renormalized, and the naive relation $\etah_1^{p-1}b=(-)^{p-1}\chih_1^p$ is broken by an anomalous term. By direct computation one finds $\etah_{R1}^{p-1}b=(-)^{p-1}\chih_1^p +\psi_1^p$, where $\psi_1^p$ involves commutators with $ |\Dh|^{-z} $:
\beq
\lefteqn{\psi_1^p( x_0\dd x_1\ldots \dd x_{p+1} ) = } \non\\
&& \frac{\Gamma(p/2)}{(p+1)!} \frac{1}{2} \sum_{i=1}^p \Pf_{z=0}\tau\nat\Big( \sum_{j=1}^i (-)^j x_0F[F,x_1] \ldots [|\Dh|^{-z},x_j] \ldots \dd x_{i+1} \ldots [F,x_{p+1}] \non\\
&& \qquad\qquad\qquad\quad + \sum_{j=i+1}^{p+1} (-)^j Fx_0[F,x_1] \ldots \dd x_i \ldots [|\Dh|^{-z},x_j] \ldots [F,x_{p+1}] \Big) \ .\non
\eeq
Problematically $\psi_1^p$ is not $\kappa$-invariant and its compatibility with the $G$-filtration is not obvious. Using the abstract pseudodifferential calculus we expand the commutator 
$$
[|\Dh|^{-z},x_j] \simeq z\sum_{l\geq 1} \frac{(-)^l}{l!} \frac{\Gamma(z+l)}{\Gamma(z+1)}\, \delta^l(x_j) |\Dh|^{-z-l}\ , 
$$  
and $\psi_1^p$ can be rewritten as a pro-finite sum of residues, typically
\beq
 && \res \tau\nat(x_0F[F,x_1]\ldots \delta^l(x_j) |\Dh|^{-z-l} \ldots \dd x_{i+1} \ldots [F,x_{p+1}]) \non\\
&\mbox{or}& \res \tau\nat(Fx_0[F,x_1] \ldots \dd x_i \ldots [|\Dh|^{-z},x_j] \ldots [F,x_{p+1}] ) \non
\eeq
for some $l\geq 1$. Suppose $x_i\in (J\Ac_0)^{k_i}$, so that $x_0\dd x_1\ldots \dd x_{p+1}\in H^k\Om^{p+1}T\Ac_0$ with $k=k_1+\ldots+k_{p+1}$. Then one has $[F,x_i]\in (\Psi\Ac_0)^{k_i}_{-1}$ and 
$$
\delta^l(x_j) |\Dh|^{-z-l}\in (\Psi\Ac_0)^{k_j}_{-z-l} + (\Psi\Ac_0)^{k_j+1}_{-z-l+1}+ \ldots + (\Psi\Ac_0)^{k_j+l}_{-z}\ .
$$
It means that modulo the adic degrees $\geq k+1$, only the residue corresponding to $l=1$ remains because for $l> 1$ the trace is finite. Hence for $x_0\dd x_1\ldots \dd x_{p+1}\in H^k\Om^{p+1}T\Ac_0$ we have
\beq
\lefteqn{\psi_1^p( x_0\dd x_1\ldots \dd x_{p+1} ) \equiv } \non\\
&& \mspace{-30mu}  - \frac{\Gamma(p/2)}{(p+1)!} \frac{1}{2} \sum_{i=1}^p \res \tau\nat\Big( \sum_{j=1}^i (-)^j x_0F[F,x_1] \ldots \delta x_j |\Dh|^{-z-1} \ldots \dd x_{i+1} \ldots [F,x_{p+1}] \non\\
&& \qquad\qquad\qquad\quad + \sum_{j=i+1}^{p+1} (-)^j Fx_0[F,x_1] \ldots \dd x_i \ldots \delta x_j |\Dh|^{-z-1} \ldots [F,x_{p+1}] \Big) \non
\eeq
modulo $\nat(\Jch^{k+1}\dd\Rch + \Jch^k\dd\Jch)$. The operator under the residue has symbol degree zero. Therefore we can move $|\Dh|^{-z-1}$ and $F$ in first position because their commutators with the $x_i$'s have lower symbol degree unless they raise the adic degree. One gets
\beq
\lefteqn{\psi_1^p( x_0\dd x_1\ldots \dd x_{p+1} ) \equiv } \non\\
&& \mspace{-30mu} -\frac{\Gamma(p/2)}{(p+1)!} \frac{1}{2} \sum_{i=1}^p \res \tau\nat \, F|\Dh|^{-z-1}\Big( \sum_{j=1}^i (-)^j x_0[F,x_1] \ldots \delta x_j \ldots \dd x_{i+1} \ldots [F,x_{p+1}] \non\\
&& \mspace{200mu} + \sum_{j=i+1}^{p+1} (-)^j x_0[F,x_1] \ldots \dd x_i \ldots \delta x_j \ldots [F,x_{p+1}] \Big) \non
\eeq
modulo $\nat(\Jch^{k+1}\dd\Rch + \Jch^k\dd\Jch)$. Now it is easy to show the $\kappa$-invariance of the r.h.s. as well as the needed compatibility with the $G$-filtration, still modulo $\nat(\Jch^{k+1}\dd\Rch + \Jch^k\dd\Jch)$. This is allows to apply verbatim the proof of \cite{P5} Proposition 3.9 and show that $\etah_R^{p-1}\rho_*\gamma$ is a cochain of order $p$ as wanted.\\

\noindent The formula for $\chi_R$ is an involved computation which relies on the techniques developed in \cite{P1, P2} in connection with a JLO approach to the bivariant Chern character. We use Quillen's formalism of algebra cochains \cite{Q2} and let $\Cc$ be the bar coalgebra of the unitalized completed tensor algebra $ (\Th\Ac_0)^+ $. Thus $\Cc$ is the graded direct sum of $n$-fold tensor products 
$$
\Cc_0=\cc\ ,\qquad \Cc_n = (\Th\Ac_0)^+ \hotimes \ldots \hotimes (\Th\Ac_0)^+ \quad \mbox{($n$ times)}\ . 
$$
The coproduct on $\Cc$ comes from the linear maps $\Cc_n\to \bigoplus_{k=0}^n \Cc_k\hotimes\Cc_{n-k}$ induced by all possible decompositions of $n$-chains into a tensor product. Moreover the multiplication on $(\Th\Ac_0)^+$ leads to a Hochschild differential of odd degree $b':\Cc_n\to \Cc_{n-1}$
$$
b'(x_1\otimes\ldots \otimes x_n)=\sum_{i=1}^{n-1}(-)^{i+1} x_1\otimes \ldots \otimes x_{i}x_{i+1}\otimes\ldots \otimes x_n\ .
$$
Then $b'^2=0$ and $\Cc$ is a graded differential coalgebra. Now consider the (unital) $\zz_2$-graded algebra of non-commutative differential forms of degree $\leq 1$ over $\Pc\hotimes\Rch$
$$
\Om^*(\Pc\hotimes\Rch) := \Pc\hotimes\Rch^+ \oplus \Om^1(\Pc\hotimes\Rch)
$$ 
gifted with the derivation $\dd$. The space of linear maps $ \hom(\Cc, \Om^*(\Pc\hotimes\Rch)) $ with convolution product is a $\zz_2$-graded algebra. $\dd$ and the transposed $\boldsymbol{\delta}$ of $b'$ act on this algebra as (anticommuting) differentials of odd degree. The unital homomorphism $\rho_*:(\Th\Ac_0)^+\cong\Cc_1 \to \Pc_0\hotimes\Rch^+$ extended by zero over $\Cc_n$, $n\neq 1$, defines an element of odd degree in this DG algebra. In the same way, we regard the Dirac operator $\Dh\in \Pc_1\hotimes\Rch^+$ as a linear map $\cc\cong\Cc_0\to \Pc_1\hotimes\Rch^+$, $1\mapsto \Dh$, hence as an element of odd degree $\Dh \in \hom(\Cc, \Om^*(\Pc\hotimes\Rch))$. For any real parameter $t> 0$, the sum $\rho_*+t\Dh$ thus defines a Quillen superconnection $\nabla_t=\boldsymbol{\delta} - \dd +\rho_* + t\Dh$ on the DG algebra $\hom(\Cc, \Om^*(\Pc\hotimes\Rch))$, with curvature (see \cite{P1})
$$
\nabla_t^2 = -\dd\rho_*-t\dd\Dh +t[\Dh,\rho_*] +t^2\Dh^2\ \in \hom(\Cc, \Om^*(\Pc\hotimes\Rch) ) \ .
$$
The exponential $\exp(-\nabla_t^2)$ is an even-degree element of the convolution algebra, defined via a Duhamel-type expansion (see below). It involves the heat operator related to the complex powers of the Dirac operator through the Mellin transform:
$$
e^{-t^2\Dh^2} = \frac{1}{2\pi i} \int \frac{e^{-t^2\la^2}}{\la-|\Dh|}\, d\la\ ,\qquad |\Dh|^{-2z} = \frac{2}{\Gamma(z)}\int_0^{\infty} \frac{dt}{t} t^{2z} e^{-t^2\Dh^2}\ .
$$
The integral over the complex parameter $\la$ is performed along a vertical line separating zero from the spectrum of $|D|$ as in Eq. (\ref{int}). Next we let $\Om_1\Cc=\Cc\hotimes(\Th\Ac_0)^+\hotimes\Cc$ be the universal bicomodule (\cite{Q2}) over the coalgebra $\Cc$. The left comodule map $\Om_1\Cc \to \Cc\hotimes\Om_1\Cc$ and the right comodule map $\Om_1\Cc \to \Om_1\Cc\hotimes\Cc$ are defined via the coproduct on the left and right factor $\Cc$ respectively. Hence by duality the space of linear maps $\hom(\Om_1\Cc, \Om^*(\Pc\hotimes\Rch) )$ is a bimodule over the convolution algebra $\hom(\Cc, \Om^*(\Pc\hotimes\Rch) )$. Moreover the transposed of the universal coderivation $\Om_1\Cc \to \Cc$, which identifies a tensor product $\Cc_k\hotimes (\Th\Ac_0)^+\hotimes\Cc_{n-k}$ with $\Cc_{n+1}$, induces a derivation
$$
\d: \hom(\Cc, \Om^*(\Pc\hotimes\Rch) ) \to \hom(\Om_1\Cc, \Om^*(\Pc\hotimes\Rch) ) \ .
$$
In particular $\d\rho_*$ defines an element of this bimodule. The core of the bivariant JLO Chern character constructed in \cite{P1} is the element 
\be
\mu(t) = \int_0^1 ds\, e^{-s\nabla_t^2} \d\rho_* e^{(s-1)\nabla_t^2} \ \in  \hom(\Om_1\Cc, \Om^*(\Pc\hotimes\Rch)) \ . \label{core}
\ee
Alternatively $\mu(t)$ is the term proportional to $\d\rho_*$ in the Duhamel expansion of the exponential $\exp(\d\rho_* - \nabla_t^2)$. We compose the linear map $\mu(t)$ on the right with the \emph{cotrace} $\nat: \Om \Th\Ac_0\to \Om_1\Cc$ defined on a $n$-form by
$$
\nat(x_0\dd x_1\ldots \dd x_n)= \sum_{i=0}^n (-)^{n(i+1)} (x_{i+1}\otimes \ldots \otimes x_n)\otimes x_0\otimes (x_1\otimes\ldots\otimes x_i)
$$
whence a linear map $\mu(t)\nat: \Om\Th\Ac_0\to \Om^*(\Pc\hotimes\Rch)$. Let $p_X$ be the projection of $\Om^*(\Pc\hotimes\Rch)$ onto the $X$-complex $X(\Pc\hotimes\Rch^+): \Pc\hotimes\Rch^+ \rightleftarrows \Om^1(\Pc\hotimes\Rch)_{\nat}$, and compose $\mu(t)\nat$ on the left with $p_X$ to get a linear map for any $t>0$:
$$
p_X\mu(t)\nat\in \hom(\Om\Th\Ac_0, X(\Pc\hotimes\Rch))\ .
$$
Proposition 6.5 of \cite{P1} shows  that $p_X\mu(t)\nat$ is a \emph{chain map} from the $(b+B)$-complex of differential forms over $\Th\Ac_0$ to the $X$-complex:
$$
(\nat\dd\oplus \bb)\circ p_X\mu(t)\nat - p_X\mu(t)\nat \circ (b+B)=0\ .
$$
This is a direct consequence of the Bianchi identity $[\nabla_t,\exp(-\nabla_t^2)]=0$. Then for $z\in \cc$ consider the cocycle
$$
\chi(z) =   \frac{2}{\Gamma(z)} \int_0^{\infty} \frac{dt}{t} t^{2z}\, p_X\mu_t\nat\ \in \hom(\Om\Th\Ac_0, X(\Pc\hotimes\Rch))\ .
$$
It splits into a set of components $\chi_0^n(z):\Om^n\Th\Ac_0 \to \Pc\hotimes\Rch$ and $\chi_1^n(z):\Om^{n+1}\Th\Ac_0 \to \nat \Om^1(\Pc\hotimes\Rch)$ for any integer $n$. In order to establish explicit formulas let $\mu_0(t)$ be the component of the map $\mu(t)$ with range contained in $\Pc\hotimes\Rch^+$. It is defined via the Duhamel series
$$
\mu_0(t) = \sum_{n\geq 0} (-t)^{n} \int_{\Delta_{n+1}}ds\,  e^{-s_0t^2\Dh^2} A_0 e^{-s_1t^2\Dh^2}A_1 \ldots e^{-s_nt^2\Dh^2}A_n e^{-s_{n+1}t^2\Dh^2} \ ,
$$
where $\Delta_{n+1}$ is the standard simplex with coordinates $s_0,\ldots,s_{n+1}$, one of the $A_i$'s is equal to $\d\rho_*$ while the others are equal to $[\Dh,\rho_*]$, and we sum over all such possibilities. The $n$-th term of this series contributes to the component $\chi^n_0(z)$:
$$
\chi^n_0(z) = \frac{2(-)^n}{\Gamma(z)} \int_0^{\infty} \frac{dt}{t} t^{2z+n} \int_{\Delta_{n+1}}ds\,  e^{-s_0t^2\Dh^2} A_0 \ldots A_n e^{-s_{n+1}t^2\Dh^2}\nat \ .
$$
The range of this map is not contained in the subspace of pseudodifferential operators $\Psi\Ac_0\subset \Pc\hotimes\Rch$. However we can write an asymtotic expansion of $\chi^n_0(z)$ by moving all the heat operators $\exp(-s_it^2\Dh^2)$ to the right (see \cite{CM95}),
\beq
\lefteqn{\int_0^{\infty} \frac{dt}{t} t^{2z+n} \int_{\Delta_{n+1}}ds\, e^{-s_0t^2\Dh^2} A_0 e^{-s_1t^2\Dh^2} A_1 \ldots e^{-s_nt^2\Dh^2} A_n e^{-s_{n+1}t^2\Dh^2}   } \non\\
&\simeq& \sum_{k_0,\ldots,k_n \geq 0} \int_0^{\infty} \frac{dt}{t} t^{2z+n}(-t^2)^{k} \times \non\\
&& \qquad \qquad  \int_{\Delta_{n+1}}ds\,  \frac{s_0^{k_0}}{k_0!}  \ldots \frac{(s_0+\ldots +s_n)^{k_n}}{k_n!}\,  A_0^{(k_0)} A_1^{(k_1)} \ldots A_n^{(k_n)} e^{-t^2\Dh^2} \non\\
&\simeq& \sum_{k_0,\ldots,k_n\geq 0} (-)^k\int_0^{\infty} \frac{dt}{t} t^{2z+n+2k}    c(k_0,\ldots, k_n) \, A_0^{(k_0)} A_1^{(k_1)} \ldots A_n^{(k_n)} e^{-t^2\Dh^2} \non\ ,
\eeq
where $k=k_0+\ldots +k_n$, $c(k_0,\ldots,k_n)$ is the constant defined above, $A_i^{(k_i)}$ is the $k_i$-th power of the derivation $[\Dh^2,\ ]$ on $A_i$, and $\simeq$ means equality modulo the space of linear maps $\hom(\Om_1\Cc, \Fc_{\al})$ for any $\al\in \rr$. Perform the integral over $t$ using the Mellin transform
$$
\frac{2}{\Gamma(z)}\int_0^{\infty} \frac{dt}{t} t^{2z+2k+n} e^{-t^2\Dh^2} = \frac{\Gamma(z+k+n/2)}{\Gamma(z)} \, |\Dh|^{-2z-2k-n}\ ,
$$
so that $\chi^n_0(z)\in \hom(\Om^n\Th\Ac_0, \Pc\hotimes\Rch)$ admits the asymptotic expansion
\beq
\lefteqn{\chi^n_0(z)\simeq} \non\\
&& \mspace{-40mu} \sum_{k_0,\ldots,k_n\geq 0}\mspace{-20mu} (-)^{n+k} c(k_0,\ldots, k_n)  \frac{\Gamma(z+k+\frac{n}{2})}{\Gamma(z)}\, A_0^{(k_0)} A_1^{(k_1)} \ldots A_n^{(k_n)} |\Dh|^{-2z-2k-n}\nat\ . \non
\eeq
Observe that the $(k_0,\ldots,k_n)$-th term of the expansion is a linear map from $\Om^n\Th\Ac_0$ to the subspace of pseudodifferential operators $\Psi\Ac_0 \cap \Fc_{\al}$ with $\al=-2\re(z)-k-n$. Hence we can apply the partial trace $\tau:\Psi\Ac_0\to \Rch$ for $z\gg 0$ and define the renormalized components $\chi_{R0}^n\in\hom(\Om^n\Th\Ac_0,\Rch)$ by taking the finite part at $z=0$:
\beq
\lefteqn{\chi^n_{R0}=} \non\\
&& \mspace{-40mu} \sum_{k_0,\ldots,k_n\geq 0}\mspace{-20mu} (-)^{n+k} c(k_0,\ldots, k_n)  \Pf_{z=0} \frac{\Gamma(z+k+\frac{n}{2})}{\Gamma(z)} \tau (A_0^{(k_0)} A_1^{(k_1)} \ldots A_n^{(k_n)} |\Dh|^{-2z-2k-n})\nat\ . \non
\eeq
Since $1/\Gamma(z) = z / \Gamma(z+1) $ the latter is a sum of residues:
$$
\sum_{k_0,\ldots,k_n\geq 0}\mspace{-20mu} (-)^{n+k} c(k_0,\ldots, k_n)  \res \frac{\Gamma(z+k+\frac{n}{2})}{\Gamma(z+1)} \tau (A_0^{(k_0)} A_1^{(k_1)} \ldots A_n^{(k_n)} |\Dh|^{-2z-2k-n})\nat
$$
The explicit evaluation of this cochain on a $n$-form $x_0\dd x_1\ldots\dd x_n \in \Om^n\Th\Ac_0$ yields formula (\ref{chir0}). For the reason explained above the sum over $k_0,\ldots,k_n$ is actually pro-finite, and the sum of all components $\sum_{n}\chi^n_{R0}$ really extends to a cochain in $\hom(\Omh\Th\Ac_0,\Rch)$. \\
The degree-one component $\mu_1(t)\in \hom(\Om_1\Cc,\Om^1(\Pc\hotimes\Rch))$ of $\mu(t)$ is also defined via a Duhamel series,
$$
\mu_1(t) = \sum_{n\geq 0} t^{d} \int_{\Delta_{n+1}}ds\,  e^{-s_0t^2\Dh^2} B_0 e^{-s_1t^2\Dh^2}B_1 \ldots e^{-s_nt^2\Dh^2}B_n e^{-s_{n+1}t^2\Dh^2} \ ,
$$
where now the $B_i$'s are chosen among $\d\rho_*$, $-[\Dh,\rho_*]$, $\dd\rho_*$, $\dd\Dh$, and the power $d$ is equal to the number of operators $\Dh$ involved. Applying the universal trace and cotrace on both sides, the one-forms $\dd\rho_*$ and $\dd\Dh$ can be pushed to the right so that the composite $\nat\mu_1(t)\nat \in \hom(\Om\Th\Ac_0,\Om^1(\Pc\hotimes\Rch)_{\nat})$ is
\beq
\lefteqn{\nat\mu_1(t)\nat = \sum_{n\geq 0} (-t)^{n} \int_{\Delta_{n+1}}ds \, \nat e^{-s_0t^2\Dh^2} A_0  \ldots e^{-s_nt^2\Dh^2}A_n e^{-s_{n+1}t^2\Dh^2} \dd\rho_* \nat } \non\\
&+& \sum_{n\geq 0} (-t)^{n}t \int_{\Delta_{n+1}}ds \, \nat e^{-s_0t^2\Dh^2} A_0  \ldots e^{-s_nt^2\Dh^2}A_n e^{-s_{n+1}t^2\Dh^2} \dd\Dh\nat \ ,\non
\eeq
where one of the $A_i$'s is $\d\rho_*$ and the others are equal to $[\Dh,\rho_*]$. Moving the heat operators to the right without crossing $\dd$, one obtains as before an asymptotic expansion for the cochain $\chi^n_1(z)\in \hom(\Om^{n+1}\Th\Ac_0,\Om^1(\Pc\hotimes\Rch)_{\nat})$
\beq
\lefteqn{\chi^n_1(z)\simeq} \non\\
&& \mspace{-40mu} \sum_{k_0,\ldots,k_n\geq 0}\mspace{-20mu} (-)^{n+k} c(k_0,\ldots, k_n)  \frac{\Gamma(z+k+\frac{n}{2})}{\Gamma(z)}\nat A_0^{(k_0)} A_1^{(k_1)} \ldots A_n^{(k_n)} |\Dh|^{-2(z+k)-n}\dd\rho_*\nat \non\\
&-& \mspace{-40mu} \sum_{k_0,\ldots,k_{n+1}\geq 0}\mspace{-20mu} (-)^{n+k} c(k_0,\ldots, k_{n+1})  \frac{\Gamma(z+k+\frac{n}{2}+1)}{\Gamma(z)}\nat A_0^{(k_0)} \ldots A_{n+1}^{(k_{n+1})} |\Dh|^{-2(z+k+1)-n}\dd\Dh\nat \non
\eeq
where $\simeq$ means equality modulo $\hom(\Om^{n+1}\Th\Ac_0,\nat(\Om^1\Fc)_{\al})$ for any $\al\in\rr$. The $(k_0,\ldots,k_n)$-th term of the expansion ranges in the subspace of pseudodifferential operators $\nat\Om^1(\Psi\Ac_0)\cap \nat(\Om^1\Fc)_{\al}$ with $\al= -2\re(z)-k-n$. Hence we can apply the trace $\tau: \nat\Om^1(\Psi\Ac_0)\to \Om^1\Rch_{\nat}$ for $z\gg 0$ and define the renormalized components $\chi_{R1}^n\in \hom(\Om^{n+1}\Th\Ac_0, \Om^1\Rch_{\nat})$ by taking the finite part at $z=0$. Alternatively it is the (pro-finite) sum of residues
\beq
\lefteqn{\chi^n_{R1}=  \sum_{k_0,\ldots,k_n\geq 0} (-)^{n+k} c(k_0,\ldots, k_n) \times } \non\\
&& \qquad \qquad \res \Big( \frac{\Gamma(z+k+\frac{n}{2})}{\Gamma(z+1)}\nat A_0^{(k_0)} A_1^{(k_1)} \ldots A_n^{(k_n)} |\Dh|^{-2(z+k)-n}\dd\rho_*\nat \Big) \non\\
&& - \sum_{k_0,\ldots,k_{n+1}\geq 0}  (-)^{n+k} c(k_0,\ldots, k_{n+1}) \times \non\\
&& \qquad\qquad  \res \Big( \frac{\Gamma(z+k+\frac{n}{2}+1)}{\Gamma(z+1)}\nat A_0^{(k_0)} \ldots A_{n+1}^{(k_{n+1})} |\Dh|^{-2(z+k+1)-n}\dd\Dh\nat\Big)\ . \non
\eeq
The explicit evaluation of this cochain on a $(n+1)$-form $x_0\dd x_1\ldots\dd x_{n+1} \in \Om^{n+1}\Th\Ac_0$ yields formula (\ref{chir1}). Since the residues vanish for large degrees the sum $\sum_{n}\chi^n_{R1}$ extends to a cochain in $\hom(\Omh\Th\Ac_0,\Om^1\Rch_{\nat})$. Collecting the even and odd part we thus obtain a cochain in the bivariant periodic complex
$$
\chi_R \in \hom(\Omh\Th\Ac_0, X(\Rch))\ .
$$
It is a cocycle as a consequence of the fact that the cocycle property for $\chi(z)$ is algebraic and thus passes to asymptotic expansions in terms of pseudodifferential operators.\\ 
In order to relate $\chi_R$ to the renormalized eta-cochain we use the same tricks and introduce a parameter $u\in \rr_+$. Let $(\Om^*\rr_+=\cinf(\rr_+)\oplus \Om^1\rr_+, d_u)$ be the de Rham complex of ordinary differential forms with respect to $u$. On the convolution DG algebra $\hom(\Cc,\Om^*(\Pc\hotimes\Rch)\hotimes\Om^*\rr_+)$ we take the sum $\tilde{\nabla}_t$ of the connections $-d_u+uF$ and $\nabla_t=\boldsymbol{\delta}-\dd +\rho_*+t\Dh$; the new curvature reads
$$
(\tilde{\nabla}_t)^2 = -du F -\dd\rho_* -t\dd\Dh +[uF+t\Dh,\rho_*] +(u+t|\Dh|)^2\ .
$$
Then define $\tilde{\mu}(t) \in \hom(\Om_1\Cc,\Om^*(\Pc\hotimes\Rch)\hotimes\Om^*\rr_+)$ by formula (\ref{core}) with $\nabla_t$ replaced by $\tilde{\nabla}_t$. Write $\tilde{\mu}(t)=\mu(t)+ du \, \nu(t)$ its expansion in powers of the one-form $du$. Hence $\mu(t)$ and $\nu(t)$ are both elements of $\hom(\Om_1\Cc,\Om^*(\Pc\hotimes\Rch)\hotimes\cinf(\rr_+))$. Remark that $\mu(t,u=0)$ coincides with the old expression (\ref{core}) for any $t>0$. By \cite{P1} Proposition 6.7, the Bianchi identity for the total connection $\tilde{\nabla}_t$ leads to the transgression relation
$$
\frac{\d}{\d u} (p_X\mu(t,u)\nat) = (\nat\dd\oplus \bb)\circ (p_X\nu(t,u)\nat) + (p_X\nu(t,u)\nat) \circ (b+B)
$$
in the complex $\hom(\Om\Th\Ac_0,X(\Pc\hotimes\Rch))$. Then consider the cochain
$$
\eta(z)= \frac{-2}{\Gamma(z)} \int_0^{\infty} \frac{dt}{t} t^{2z}  \int_{0}^{\infty} du\, p_X\nu(t,u)\nat \ \in \hom(\Om\Th\Ac_0,X(\Pc\hotimes\Rch))\ . 
$$
Its boundary is related to the cocycle $\chi(z)$ by means of the above transgression:
\be
(\nat\dd\oplus \bb)\circ \eta(z) + \eta(z)\circ (b+B) = \frac{2}{\Gamma(z)} \int_0^{\infty} \frac{dt}{t} t^{2z}\,  p_X\mu(t,0)\nat = \chi(z) \label{trans}
\ee
Let us calculate the asymptotic expansion of $\eta(z)$. The component of $\nu(t,u)$ landing in degree zero $\nu_0(t,u)\in \hom(\Om_1\Cc,\Pc\hotimes\Rch^+)$ is defined via the Duhamel series
$$
\nu_0(t,u)= \sum_{n\geq 0} (-t)^d (-u)^f \int_{\Delta_{n+2}}ds\,  e^{-s_0(u+t|\Dh|)^2} A_0  \ldots A_{n+1} e^{-s_{n+2}(u+t|\Dh|)^2} 
$$
where two of the $A_i$'s are equal to $\d\rho_*$ and $F$ respectively, the others are equal to $[\Dh,\rho_*]$ or $[F,\rho_*]$, and we sum over all such possibilities. $d$ and $f$ are the numbers of commutators $[\Dh,\rho_*]$ and $[F,\rho_*]$ respectively ($d+f=n$). Only the term of order $n$ contributes to the component $\eta^n_0(z)\in \hom(\Om^n\Th\Ac_0,\Pc\hotimes\Rch)$. Moving all the heat operators $\exp(-s_i(u+t|\Dh|)^2)$ to the right we obtain an asymptotic expansion
\beq
\lefteqn{\eta^n_0(z)\simeq   \sum_{k_0,\ldots,k_{n+1}\geq 0} (-)^{n+k} c(k_0,\ldots, k_{n+1}) \times } \non\\
&& \qquad \qquad \frac{-2}{\Gamma(z)} \int_0^{\infty} \frac{dt}{t} t^{2z+d}\int_0^{\infty}du\, u^f\, A_0^{[k_0]} \ldots A_{n+1}^{[k_{n+1}]} e^{-(u+t|\Dh|)^2} \nat \non
\eeq
where $k=k_1+\ldots +k_{n+1}$, $A_i^{[k_i]}$ is the $k_i$-th power of the derivation $[(u+t|\Dh|)^2,\ ]$ on $A_i$. Each term of the expansion lies in $\hom(\Om^n\Th\Ac_0, \Psi\Ac_0)$. One has $A_i^{[1]}= 2tu\, \delta A_i + t^2A_i^{(1)}$ with the derivation $\delta=[|\Dh|, \ ]$, hence recursively
$$
A_i^{[k_i]} = \sum_{p_i=0}^{k_i} \binom{p_i}{k_i}  t^{2k_i-p_i}(2u)^{p_i}\,  \delta^{p_i} A^{(k_i-p_i)}\ .
$$
Hence we can perform the integrals over $t,u$ using the identity
\beq
\lefteqn{\frac{2}{\Gamma(z)} \int_0^{\infty} \frac{dt}{t} t^{2z+2k-p+d}  \int_{0}^{\infty}du\, u^{p+f} \, e^{-(u+t|\Dh|)^2} = } \non\\
&& \qquad \qquad \frac{1}{\Gamma(z)} \Big( \sum_{l=0}^{p+f} \binom{l}{p+f} \frac{(-)^l\, \Gamma(z+k+\frac{n+1}{2})}{2z+2k-p+d+l} \Big)\, |\Dh|^{-2z-2k+p-d} \non
\eeq
with $k=k_1+\ldots +k_{n+1}$ and $p=p_1+\ldots + p_{n+1}$, so that
\beq
\lefteqn{\eta^n_0(z)\simeq  - \sum_{\substack{ k_0,\ldots,k_{n+1}\geq 0 \\ 0\leq p_i \leq k_i}}\mspace{-20mu} (-)^{n+k} c(k_0,\ldots, k_{n+1}) 2^p \prod_{i=0}^{n+1}\binom{p_i}{k_i} \times } \non\\
&& \mspace{-30mu} \frac{1}{\Gamma(z)}\Big( \sum_{l=0}^{p+f} \binom{l}{p+f} \frac{(-)^l\, \Gamma(z+k+\frac{n+1}{2})}{2z+2k-p+d+l} \Big)\, \Big(\prod_{i=0}^{n+1}\delta^{p_i} A_i^{(k_i-p_i)}\Big) \, |\Dh|^{-2z-2k+p-d} \nat \non
\eeq
The renormalized component $\eta^n_{R0}\in \hom(\Om^n\Th\Ac_0, \Rch)$ is defined by applying the partial trace $\pm\tau:\Psi\Ac_0\to\Rch$ to this asymptotic expansion for $z\gg 0$, and then taking the finite part at $z=0$ (the sign $\pm$ is dictated by the parity of the quasihomomorphism). We want to show that, modulo counterterms, it coincides with the renormalized eta-cochain $\etah^n_{R0}$ given by Equation (\ref{ren}). $\eta^n_{R0}$ is a sum of terms proportional to
$$
\Pf_{z=0} \frac{1}{\Gamma(z)}\Big( \sum_{l=0}^{p+f} \binom{l}{p+f} \frac{(-)^l\, \Gamma(z+k+\frac{n+1}{2})}{2z+2k-p+d+l} \Big)\, \tau\Big(\prod_{i=0}^{n+1}\delta^{p_i} A_i^{(k_i-p_i)} \, |\Dh|^{-2z-2k+p-d}\Big) \nat 
$$
with $d+f=n$ and $p\leq k$. Since $1/\Gamma(z)=z/\Gamma(z+1)$, the meromorphic function
$$
\frac{1}{\Gamma(z)}\Big( \sum_{l=0}^{p+f} \binom{l}{p+f} \frac{(-)^l\, \Gamma(z+k+\frac{n+1}{2})}{2z+2k-p+d+l} \Big)
$$
vanishes at $z=0$, unless $d=k=p=0$. It follows that the series $\eta^n_{R0}$ is the sum of its first term corresponding to $d=k=p=0$, plus counterterms
\beq
\eta^n_{R0} &=& \mp \frac{(-)^n}{(n+2)!} \Pf_{z=0} \Big( \sum_{l=0}^{n} \binom{l}{n} \frac{(-)^l}{2z+l} \Big)\frac{\Gamma(z+\frac{n+1}{2})}{\Gamma(z)} \, \tau(A_0\ldots A_{n+1} |\Dh|^{-2z})\nat\non\\
&& +\ \mbox{counterterms}\ , \non
\eeq
where the counterterms are residues taking into account only the poles of the zeta-functions
$$
\tau (\delta^{p_0} A_0^{(k_0-p_0)} \ldots \delta^{p_{n+1}} A_{n+1}^{(k_{n+1}-p_{n+1})}  |\Dh|^{-2z-2k+p-d}) \nat \ .
$$
Observe that $\delta^{p_0} A_0^{(k_0-p_0)} \ldots \delta^{p_{n+1}} A_{n+1}^{(k_{n+1}-p_{n+1})}  |\Dh|^{-2z-2k+p-d}$ ranges in the subspace of pseudodifferential operators $\Psi\Ac_0 \cap \Fc_{\al}$ with $\al=-2\re(z)-k-n$. Hence the residues vanish for large $n+k$ after projection $\Rch\to \Rc/\Jc^m$, and sum over the counterterms is actually pro-finite. We can simplify further the first term of the series by noting that the meromorphic function in front of $\tau(A_0\ldots A_{n+1} |\Dh|^{-2z})\nat$ has an analytic expansion $\Gamma(\frac{n+1}{2})/2 + zh(z)$ with $h(z)$ holomorphic around zero. Hence we may absorb the contribution of $h$ inside the counterterms:
$$
\eta^n_{R0} = \mp(-)^n \frac{\Gamma(\frac{n+1}{2})}{(n+2)!}\, \frac{1}{2} \Pf_{z=0} \tau(A_0\ldots A_{n+1} |\Dh|^{-2z})\nat +\ \mbox{counterterms}\ . \non
$$
Now evaluate this map on a $n$-form $ x_0\dd x_1\ldots\dd x_n \in \Om^n\Th\Ac_0$ and recall that the parity of $n$ is opposite to the parity of the quasihomomorphism:
\beq
\lefteqn{\eta^n_{R0}( x_0\dd x_1\ldots\dd x_n )= \frac{\Gamma(\frac{n+1}{2})}{(n+1)!}\, \frac{1}{2} \Pf_{z=0} \tau\Big( F x_0[F,x_1]\ldots [F,x_n] |\Dh|^{-2z}} \non\\
&+&  \sum_{i=1}^n (-)^{ni} [F,x_i]\ldots [F,x_n]Fx_0[F,x_1]\ldots [F,x_{i-1}] |\Dh|^{-2z} \Big) +\ \mbox{counterterms} \non
\eeq
This almost corresponds to the renormalized eta-cochain $\etah_{R0}^n$ given by (\ref{ren}), except for the fact that the regularizing operator $|\Dh|^{-2z}$ is not inserted at the same places. Using commutators with $|\Dh|^{-2z}$ we can nevertheless conclude that modulo counterterms of the form (\ref{contre}), the linear maps $\eta_{R0}^n$ and $\etah_{R0}^n$ coincide in $\hom(\Om^n\Th\Ac,\Rch)$. Moreover, since the counterterms vanish for large $n$ after projection $\Rch\to\Rc/\Jc^m$, the \emph{difference} $\sum_{n}(\eta^n_{R0}-\etah^n_{R0})$ extends to a well-defined linear map in $\hom(\Omh\Th\Ac_0,\Rch)$. One proceeds in exactly the same way with the other components $\eta^n_{R1}$: they coincide with $\etah_{R1}^n$ of Eq. (\ref{ren}) modulo counterterms, and the difference $\sum_{n}(\eta^n_{R1}-\etah^n_{R1})$ extends to a well-defined linear map in $\hom(\Omh\Th\Ac_0,\Om^1\Rch_{\nat})$. Hence we are simply dealing with a choice of renormalization for the eta-cochain $\eta_R\in \hom(\Om\Th\Ac_0,X(\Rch))$, which differs from the renormalization (\ref{ren}) by a cochain in the complex of bivariant cyclic cohomology $\hom(\Omh\Th\Ac_0,X(\Rch))$. Lemma \ref{lan} thus implies that its boundary represents the bivariant Chern character. Moreover
$$
(\nat\dd\oplus \bb)\circ \eta_R + \eta_R \circ (b+B) = \chi_R\ \in \hom(\Omh\Th\Ac_0,X(\Rch))
$$
because the transgression relation (\ref{trans}) is algebraic and passes to asymptotic expansions. Hence $\chi_R$ represents the bivariant Chern character in periodic cyclic cohomology. \cqfd\\

\begin{example}\label{sequi}
\textup{ As an illustration of the anomaly formula we shall derive a simplified version of the equivariant index theorem stated in \cite{P3}. Let $G$ be a countable discrete group acting smoothly and properly by isometries on a complete $p$-dimensional Riemannian manifold $M$, without boundary. The quotient space $G\backslash M$ is supposed compact. Consider a $G$-equivariant complex vector bundle $E\to M$ endowed with a $G$-invariant hermitian structure, together with a $G$-invariant selfadjoint elliptic differential operator of order one $D_0:\cinf(E)\to\cinf(E)$ acting on the smooth sections of $E$. When $p$ is even, we suppose that $E$ is $\zz_2$-graded and $D_0$ is of odd degree, whereas if $p$ is odd everything is trivially graded. $D_0$ extends to an unbounded selfadjoint operator with dense domain on the Hilbert space $H_0=L^2(E)$ of square-integrable sections of $E$ with respect to the Riemannian metric and hermitian structure. If $p$ is even, form the Hilbert space $H=H_0\oplus H_0$. If $p$ is odd form the Hilbert space $H=(H_0\oplus H_0)\hotimes C_1$, where $C_1=\cc\oplus\eps\cc$ is the $\zz_2$-graded Clifford algebra with odd generator $\eps$. In both cases $H$ has a $\zz_2$-grading: in the even case it comes from the $\zz_2$-grading on each summand $H_0$, whereas in the odd case it comes from the Clifford algebra. Define the selfadjoint unbounded operator of odd degree on $H$
$$
D=\left(\begin{matrix}
D_0 & 1 \\
1 & -D_0 \end{matrix} \right) \quad \mbox{$p$ even}, \qquad D=\eps\left(\begin{matrix}
D_0 & 1 \\
1 & -D_0 \end{matrix} \right) \quad \mbox{$p$ odd}.
$$
Since $\eps^2=1$, one gets $D^2=1+D_0^2>0$ in both cases, hence $D$ is always invertible. Let $\cinfc(M)$ denote the algebra of smooth $\cc$-valued functions with compact support on $M$. Any function $f\in \cinfc(M)$ acts on the space of sections $H_0=L^2(E)$ by pointwise multiplication, thus is represented by a bounded operator say $f_0$ on $H_0$. We extend this representation to a bounded operator of even degree on  $H$ by the matrix $f=\bigl(\begin{smallmatrix} f_0 & 0 \\ 0 & 0 \end{smallmatrix} \bigr) \in \Lc(H)$. Put $F=D/|D|$, $F^2=1$, and decompose the Hilbert space $H$ into the direct sum $H_+\oplus H_-$ (even case) or $(H_+\oplus H_-)\hotimes C_1$ (odd case) according to which the matricial form of the operator $F$ reads
$$
F=\left(\begin{matrix}
0 & 1 \\
1 & 0 \end{matrix} \right) \quad \mbox{$p$ even}, \qquad F=\eps\left(\begin{matrix}
1 & 0 \\
0 & -1 \end{matrix} \right) \quad \mbox{$p$ odd}.
$$
In this new decomposition a function $f\in \cinfc(M)$ is represented by a diagonal matrix $\bigl(\begin{smallmatrix} f_+ & 0 \\ 0 & f_- \end{smallmatrix} \bigr)$ in the even case and by a matrix $\bigl(\begin{smallmatrix} f_{++} & f_{+-} \\ f_{-+} & f_{--} \end{smallmatrix} \bigr)$ in the odd case. The commutator $[F,f]$ is always an element of the Schatten ideal $\Lc^{p+1}(H)$. Hence we may endow $\cinfc(M)$ with the norm
$$
\|f\|_{\infty,p+1}:= \|f\|_{\infty} + \|[F,f]\|_{p+1}
$$
where $\|f\|_{\infty}$ is the operator norm on $\Lc(H)$ and $\|x\|_{p+1}=(\Tr|x|^{p+1})^{\frac{1}{p+1}}$ is the Schatten norm on $\Lc^{p+1}(H)$. The norm $\|\cdot \|_{\infty,p+1}$ is submultiplicative and $G$-invariant for the natural action of $G$ on $\cinfc(M)$ by pullback. \\
Now let $\cc G$ be the group ring of $G$, i.e. the algebra generated by symbols $U^*_g$ associated to each element $g\in G$, with product $U^*_{g_1}U^*_{g_2}=U^*_{g_2g_1}$. We use a contravariant notation for convenience. Hence any element $b\in\cc G$ is a $\cc$-valued function with finite support on $G$. Suppose that a right-invariant distance $d:G\times G\to \rr_+$ is given. Then for any $\al\geq 0$ define the following norm on $\cc G$: 
$$
\|b\|_{\al}=\sum_{g\in G} \si_{\al}(g) |b(g)|\ ,\quad \forall b\in \cc G\ ,
$$
where the function $\si_{\al}(g) := (1+d(g,1))^{\al}$ fulfills the property $\si_{\al}(g_1g_2)\leq \si_{\al}(g_1)\si_{\al}(g_2)$. Then each norm $\|\cdot \|_{\al}$ is submultiplicative, and the completion of $\cc G$ with respect to this family yields a Fr\'echet $m$-algebra $\Bc$. Its elements are functions with rapid decay over $G$. Next, let $\Ac_0$ be the algebraic crossed product $\cinfc(M)\cp G$. As a vector space it corresponds to the algebraic tensor product $\cinfc(M)\otimes \cc G$, and the multiplication is defined by convolution:
$$
(f_1U^*_{g_1})(f_2U^*_{g_2}) = f_1(f_2)^{g_1}\, U^*_{g_2g_1}\ ,\qquad \forall f_i\in \cinfc(M)\ ,\ g_i\in G\ ,
$$
$(f_2)^{g_1}$ being the pullback of the function $f_2$ by the diffeomorphism $g_1$. Thus any element $a\in \Ac_0$ is a $\cinfc(M)$-valued function with finite support on $G$. For any $\al\geq 0$, endow the algebraic crossed product with the norm
$$
\| a \|_{\al} = \sum_{g\in G} \si_{\al}(g) \|a(g)\|_{\infty,p+1}\ ,\quad \forall a\in \Ac_0=\cinfc(M)\cp G\ .
$$
Using the fact that the norm $\| \cdot \|_{\infty,p+1}$ on $\cinfc(M)$ is submultiplicative and $G$-invariant, one shows that $\|\cdot\|_{\al}$ is submultiplicative for any $\al$. Hence the completion of $\Ac_0$ with respect to this family of norms is a Fr\'echet $m$-algebra $\Ac$. Now observe that because $G$ acts by isometries on $M$ and by pullback on the sections of $E$, one gets a unitary representation $r:G \to \Lc(H)$.  Let us define a homomorphism $\rho:\Ac_0 \to \Lc(H)\hotimes\Bc$ (note that $\Lc(H)$ is viewed as a Banach algebra with the operator norm) by
\be
\rho(fU^*_g)=fr(g) \otimes U^*_g\ .
\ee
It is clearly continuous with respect to the topology induced by the norms $\|\cdot\|_{\al}$ on $\Ac_0$, hence extends to a continuous homomorphism $\rho:\Ac\to \Lc(H)\hotimes\Bc$. Moreover, since $[F,f]\in \Lc^{p+1}(H)$ for any function $f\in \cinfc(M)$, one has $[F,\rho(a)]\in \Lc^{p+1}(H)\hotimes\Bc$ for any $a\in \Ac_0$ and the linear map $a\mapsto [F,\rho(a)]$ extends to a continuous linear map from $\Ac$ to $\Lc^{p+1}(H)\hotimes\Bc$. Here $\Lc^{p+1}(H)$ is provided with the Schatten norm. Let $\Ic$ be the $(p+1)$-summable Banach algebra $\Lc^{p+1}(H_+)$, and $\Ec$ be the Fr\'echet $m$-algebra $\Lc(H_+)\hotimes\Bc$. Using the isomorphism $H_+\cong H_-$ and the matrix decomposition of $\Lc(H)$ and $\Lc^{p+1}(H)$, $\rho$ defines a $(p+1)$-summable quasihomomorphism
\be
\rho: \Ac \to \Ec^s\triangleright \Ic^s\hotimes\Bc
\ee
of parity $p$ mod 2. In particular, $\rho$ induces a morphism $\rho_!:\Kt_n(\Ic\hotimes\Ac)\to \Kt_{n-p}(\Ic\hotimes\Bc)$ on topological $K$-theory (section \ref{sbiv} and \cite{P5}). The crossed product $\Ac$ is provided with a canonical class $[e]\in \Kt_0(\Ac)$ represented by the idempotent $e\in \cinfc(M)\cp G$ (see \cite{P3}): 
$$
e(g)(x)=c(x)c(g\cdot x) \qquad \forall x\in M\ ,\ g\in G\ ,
$$
where $c\in\cinfc(M)$ is any cutoff function with property $\sum_{g\in G} c(g\cdot x)^2=1$, $\forall x\in M$. Embedding the algebra $\Ac$ into $\Ic\hotimes\Ac$ via a rank-one idempotent in $\Ic$, we may view $[e] \in \Kt_0(\Ic\hotimes\Ac)$ and define the \emph{equivariant index} of the elliptic operator $D_0$ as the $K$-theory class of the convolution algebra over $G$
\be
\Ind_G(D_0)=\rho_!([e])\ \in \Kt_{-p}(\Bc)\ .
\ee
We want to establish an explicit formula for the Chern character of the index in periodic cyclic homology $HP_{-p}(\Bc)$. To this end, choose the universal free extension $0\to J\Bc\to T\Bc \to \Bc\to 0$ for $\Bc$. Because $\Ec$ is a tensor product $\Lc(H_+)\hotimes\Bc$, one can find a $T\Bc$-admissible extension with $\Mc=\Lc(H_+)\hotimes T\Bc$ and $\Nc=\Lc(H_+)\hotimes J\Bc$:
$$
\vcenter{\xymatrix{
0\ar[r] & \Nc \ar[r] & \Mc \ar[r] & \Ec \ar[r] & 0 \\
0 \ar[r] & \Ic\hotimes J\Bc \ar[r] \ar[u] & \Ic\hotimes T\Bc \ar[r] \ar[u] & \Ic\hotimes\Bc \ar[r] \ar[u] & 0}}
$$
The bivariant Chern character $\ch(\rho)\in HP^{p}(\Ac,\Bc)$ is therefore well-defined and the topological part of the Riemann-Roch-Grothendieck theorem \ref{trr} relates in a commutative diagram the morphism $\rho_!$ to the map induced in periodic cyclic homology:
\be
\vcenter{\xymatrix{
\Kt_0(\Ic\hotimes\Ac) \ar[r]^{\rho_!} \ar[d] & \Kt_{-p}(\Ic\hotimes\Bc) \ar[d] \\
HP_0(\Ac) \ar[r]^{\ch(\rho)} & HP_{-p}(\Bc) }}
\ee
Hence one gets the index formula $\ch(\Ind_G(D_0))=\ch(\rho)\cdot\ch(e)$. Since $[e]$ is actually represented by an element of the dense subalgebra $\Ac_0\subset\Ac$, we shall give a local representative of the periodic cyclic cohomology class of the bivariant Chern character $\ch_R(\rho)\equiv \ch(\rho)\in HP^p(\Ac_0,\Bc)$ using zeta-function renormalization. So let us introduce the algebra of abstract symbols in this context. For any $\beta\geq 0$, the Sobolev space $H^{\beta}=\dom (|D|^{\beta})$ is a Hilbert space for the norm $\|\xi\|^{\beta}=\||D|^{\beta}\xi\|$, $\forall\xi\in H^{\beta}$. One has $\|\xi\|\leq \|\xi\|^{\beta} \leq \|\xi\|^{\beta'}$ whenever $\beta\leq\beta'$, whence a continuous inclusion $H^{\beta'}\to H^{\beta}$. The intersection $H^{\infty}=\cap_{\beta}H^{\beta}$ is a dense subspace of $H$. Then for any $\al\in\rr$, the space of abstract symbols $\Pc_{\al}$ is defined as the set of endomorphisms $x:H^{\infty}\to H^{\infty}$ which extend to continuous linear maps $H^{\al+\beta} \to H^{\beta}$ for any $\beta\geq 0$. The family of norms inherited from the operator norm on the Banach spaces $\Lc(H^{\al+\beta},H^{\beta})$,
$$
\|x\|_{\al}^{\beta}= \sup_{\xi\in H^{\al+\beta}} \frac{\|x\cdot\xi\|^{\beta}}{\|\xi\|^{\al+\beta}}\ ,\quad \forall x\in \Pc_{\al}\ ,\ \forall\beta\geq 0\ ,
$$
turns $\Pc_{\al}$ into a Fr\'echet space. Furthermore the product $\Pc_{\al}\times\Pc_{\al'}\to \Pc_{\al+\al'}$ given by composition is continuous. One thus gets the algebra of abstract symbols by taking the inductive limit $\Pc=\varinjlim_{\al} \Pc_{\al}$ with respect to the natural maps $\Pc_{\al}\to\Pc_{\al'}$ whenever $\al\leq \al'$. In particular the Dirac operator $D$ and its modulus $|D|$ are in $\Pc_1$. We now focus on $\Pc_0$. It is a subalgebra of $\Lc(H)$ and the injection $\Pc_0\to \Lc(H)$ is continuous. Let $f\in \cinfc(M)$ be represented by a bounded operator on $H$ as above. Because $f$ is smooth, any power of the derivation $\delta=[|D|,\ ]$ on $f$ still defines a bounded operator on $H$, and this allows to show that $f$ also defines a bounded operator on each Sobolev space $H^{\beta}$. One thus gets a representation $\cinfc(M)\to \Pc_0$. Moreover, $G$ acts by isometries and commutes with $|D|$ hence is unitarily represented on each $H^{\beta}$, whence a representation $r:G\to \Pc_0\subset \Lc(H)$. It is then clear that the homomorphism $\rho_*:\Ac_0\to \Lc(H)\hotimes\Bc$ factors through the algebra $\Pc_0\hotimes\Bc$. Also, the identity $[F,f]=([D,f]-F\delta f)|D|^{-1}$ shows that $[F,f]\in \Pc_{-1}$ for any $f\in\cinfc(M)$ hence $[F,\rho(a)]\in \Pc_{-1}\hotimes\Bc$ for any $a\in \Ac_0$. The lifted homomorphism $\rho_*:T\Ac_0 \to \Pc_0\hotimes T\Bc$ reads 
$$
\rho_*(a_1\otimes \ldots\otimes a_{n}) = f_1f_2^{g_1}\ldots f_{n}^{g_{n-1}\ldots g_1} r(g_n\ldots g_1) \otimes (U^*_{g_1}\otimes \ldots\otimes U^*_{g_n})
$$
on a tensor product of elements $a_i=f_iU^*_{g_i}\in \Ac_0$, and $\rho_*(J\Ac_0)\subset \Pc_0\hotimes J\Bc$. The Dirac operator $D\in\Pc_1$, considered as an element of $\Pc_1\hotimes\Bc^+$, has an obvious lifting $\Dh=D$ in $\Pc_1\hotimes\Th\Bc^+$. This simplifies considerably the resolvent formula $(\la-|\Dh|)^{-1},$ and the complex powers $|\Dh|^z=|D|^z$ simply lie in the algebra of abstract symbols $\Pc_{\re(z)}\subset \Pc_{\re(z)}\hotimes \Th\Bc^+$. Applying any power of the derivation $\delta=[|D|,\ ]$, the $G$-invariance of $D$ shows that $\delta^n\rho_*$ and $\delta^n[D,\rho_*]$ map $\Th\Ac_0$ to the subspace $\Pc_0\hotimes\Th\Bc = \Fc_0^0$, and $\Jh\Ac_0$ to $\Pc_0\hotimes\Jh\Bc = \Fc_0^1$. Hence the quasihomomorphism is $T\Bc$-regular. Moreover, using classical estimates for the heat kernel \cite{BGV} one shows that given $f_0,\ldots f_n\in \cinfc(M)$ and $g\in G$ the zeta-function $z\mapsto \Tr(\delta^{k_0}f_0\ldots \delta^{k_n}f_nr(g) |D|^{-z})$ is holomorphic in the domain $\re(z)>p$ and extends to a meromorphic function with simple poles at integers $\leq p$ along the real axis. Hence the quasihomomorphism has finite analytic dimension $p=\dim M$ and discrete dimension spectrum. It follows that the anomaly formula of Theorem \ref{tres} applies. If $\eh\in \Th\Ac_0$ denotes any lifting of the idempotent $e\in \Ac_0$ representing the canonical $K$-theory class in $\Kt_0(\Ac)$, the Chern character $\ch(\Ind_G(D_0))$ is represented by the image of the $(b+B)$-cycle over $\Th\Ac_0$
$$
\gamma(\eh)=\eh + \sum_{n\geq 1}(-)^n\frac{(2n)!}{(n!)^2} (\eh-\frac{1}{2})(\dd \eh\dd\eh)^n
$$
under the chain map $\chi_R\rho_*:\Omh\Th\Ac_0 \to X(\Th\Bc)$. We are thus led to compute residues of zeta-functions of the kind $z^q \,\Tr(f_0df_1^{(k_1)}\ldots df_n^{(k_n)}r(g) |D|^{-2z-m})$, for some smooth functions $f_i\in \cinfc(M)$ and $g\in G$, with $m\geq 0$ and $q\geq 0$ or $m=0$ and $q=-1$. As in Theorem \ref{tres} the differential $df$ denotes the commutator $[D,f]$ and the superscript $^{(k)}$ is the $k$-th power of the derivation $[D^2,\ ]$. Since $D^2=1+D_0^2$ the complex powers of $|D|$ can be expressed in terms of the heat operator $\exp(-tD_0^2)$ via a Mellin transform $|D|^{-2z}=\frac{1}{\Gamma(z)}\int_0^{\infty}\frac{dt}{t}\, t^z e^{-t}e^{-tD_0^2}$. For $t$ goes to zero one has a short-time Laurent series expansion $\Tr(f_0df_1^{(k_1)}\ldots df_n^{(k_n)}r(g) e^{-tD_0^2})\sim \sum_i a_it^{i/2}$, with $i$ bounded below. An easy computation yields
$$
\res\, z^q\, \Tr(f_0df_1^{(k_1)}\ldots df_n^{(k_n)}r(g) |D|^{-2z-m}) = \sum_{i} a_i \res\,z^q\,\frac{\Gamma(z+\frac{m+i}{2})}{\Gamma(z+\frac{m}{2})}\ .
$$
The sum over $i$ is actually finite because the quotient of gamma-functions in the right-hand-side is holomorphic around zero for large values of $i$. The coefficients $a_i$ of the heat expansion can be explicitly computed at least when $D_0$ is a generalized Dirac operator, $E$ is a Clifford module and $G$ acts by orientation-preserving diffeomorphisms on $M$. They are given by integrals over the submanifolds $M_g\subset M$ of fixed points for $g\in G$, see for example \cite{BGV} chapter 6. Using the identification of $X(\Th\Bc)$ with the space $\Omh\Bc$ of noncommutative differential forms (\cite{P3}), the Chern character $\ch(\Ind_G(D_0))$ can be decomposed into components $\ch_n(\Ind_G(D_0)) \in \Om^n\Bc$ with $n$ of parity $p\mod 2$. But $\Om^n\Bc=\Bc^{\hotimes n}\oplus \Bc^{\hotimes (n+1)}$ may be identified with a space of functions with rapid decay over the set $G^n\cup G^{n+1}$. Assuming for simplicity that each manifold of fixed points $M_g$ has a spin structure, one finds that the evaluation of $\ch_n(\Ind_G(D_0))$ on a point $\gt\in G^n\cup G^{n+1}$ reads
\be
\ch_n(\Ind_G(D_0))(\gt) = \sum_{M_g}\frac{(-)^{q/2}}{(2\pi i)^{d/2}}\int_{M_g} \widehat{A}(M_g)\, \frac{\ch(E/S,g)}{\ch(S_N,g)}\, \ch_n(e)(\gt) \ , \label{equivariant}
\ee
where the sum runs over all the submanifolds $M_g$ of dimension $d$ and codimension $q=p-d$, fixed by the element $g=g_n\ldots g_1\in G$ when $\gt=(g_1,\ldots,g_n)$ or $g=g_n\ldots g_0$ when $\gt=(g_0,\ldots,g_n)$. The integral contains the usual characteristic classes of the Atiyah-Segal-Singer equivariant index theorem \cite{AS}, namely the $\widehat{A}$-genus and the Chern characters of some equivariant vector bundles $E/S$ and $S_N$ over $M_g$. The characteristic class $\ch_n(e)$ is the most interesting part of the formula: it is a function over $G^n\cup G^{n+1}$ with values in differential forms over $M$, associated to the idempotent $e\in\Ac_0$. We refer to \cite{P3} for more details. (\ref{equivariant}) provides a generalization of the Connes-Moscovici index theorem for coverings \cite{CM90} to the case of proper $G$-actions admitting fixed points. It is worth mentioning that in \cite{P3} we state a stronger result, valid for any locally compact group $G$ and with periodic cyclic homology replaced by entire cyclic homology. The proof is based on bornological algebras and a JLO approach of the bivariant Chern character \cite{P1, P2}. With some more work it should be possible to deduce also this general form of the equivariant index theorem from the residue formula \ref{tres}.  }
\end{example}

\begin{example}\textup{The previous situation with $G=\{1\}$ and $M$ compact is well-known. A Dirac-type operator leads to a quasihomomorphism from the completion $\Ac$ of $\Ac_0=\cinf(M)$ to $\Bc=\cc$. Hence we can take $\Rc=\Bc$ as quasi-free extension, and the cyclic homology of $\Bc$ is computed by the $X$-complex $X(\cc): \cc \rightleftarrows 0$. When the dimension of $M$ is even, the Chern character $\ch(\Ind(D_0))$ is just a number calculating the index of the Dirac operator. When the dimension is odd the components of the renormalized eta-cochain become relevant. They are linear maps $\etah_{R}^{n+1}\rho_*:\Om^{n+1}\Ac_0\to \cc$, in all even degrees $n+1=2k$. In particular the component of degree zero evaluated at the constant function $1\in \Ac_0$ yields the usual eta-invariant (the odd trace is $\tau=-\sqrt{2 i} \Tr$)
$$
\etah_R^0(1)= \frac{1}{2}\Gamma(1/2) \Pf_{z=0} \tau(F|D|^{-2z}) = -\frac{1}{2}\sqrt{2\pi i} \Pf_{z=0}\Tr(F|D|^{-2z})\ ,
$$
which measures the spectral asymmetry of the Dirac operator. The other components $\etah_R^{n+1}\rho_*$ assemble into an improper $(b+B)$-cocyle over $\Ac_0$, giving higher analogues of the eta-invariant. }
\end{example}

\section{Spectral triples, anomalies and regulators}\label{sreg}

The anomaly formula (Theorem \ref{tres}) is a bivariant generalization of the local index formula of Connes and Moscovici for regular spectral triples \cite{CM95}. In the latter situation one gets a quasihomomorphism from an algebra $\Ac$ to the complex numbers  $\Bc=\cc$ and the Chern character gives a local expression for the index map $\Kt_*(\Ac)\to \mathbb{Z}$ as in the classical Atiyah-Singer index theorem. In this section we propose to relate the local index formula with the chiral anomaly of an adequate noncommutative quantum field theory. In fact noncommutative index theory and anomalies are in a sense equivalent, the link being provided by Bott periodicity \cite{P4}. Considering anomalies, however, provides a conceptual explanation for the ``locality'' of the boundary of the eta-cochain.\\ 

Hence let $\Ac_0$ be an associative algebra over $\cc$, endowed with an antilinear involution. A spectral triple $(\Ac_0,H,D)$ corresponds to the following data:\\

\noindent i) An involutive representation of $\Ac_0$ into the algebra $\Lc(H)$ of bounded operators on a separable Hilbert space $H$. \\
\noindent ii) An unbounded selfadjoint operator $D$ with compact resolvent on $H$. Without loss of generality we may assume that it is invertible.\\
\noindent iii) The commutator $[D,a]$ extends to a bounded operator on $H$ for any $a\in \Ac_0$. \\

\noindent In addition, the spectral triple is called \emph{odd} when $H$ is trivially graded, and \emph{even} when $H$ is $\zz_2$-graded with grading operator $\gamma_5\in \Lc(H)$, $(\gamma_5)^2=1$, $\gamma_5D=-D\gamma_5$ and $[\gamma_5,\rho(a)]=0$ for any $a\in \Ac_0$. If the inverse modulus of the Dirac operator $|D|^{-1}$ lies in the Schatten class $\Lc^{p+1}(H)$ the spectral triple is $(p+1)$-summable. In order to use zeta-function regularization we impose that $(\Ac_0,H,D)$ is \emph{regular} with discrete dimension spectrum \cite{CM95}, i.e. it fulfills the additional properties:\\

\noindent iv) $\Ac_0$ and $[D,\Ac_0]$ belong to the domains of all powers of the derivation $\delta=[|D|,\ ]$, in the sense that $\delta^n(a)$ and $\delta^n([D,a])$ are bounded operators on $H$ for any $a\in\Ac_0$ and $n\in\nn$.\\
v) For any element $x$ of the algebra generated by $\delta^n(\Ac_0)$ and $\delta^n([D,\Ac_0])$, the zeta-function
$$
\zeta(z)=\textup{Tr}(x |D|^{-z})\ , \quad z\in \mathbb{C}
$$
is holomorphic in the half-plane $\re(z)> p$ and extends to a meromorphic function with poles contained in a discrete set $Sd\subset\mathbb{C}$ (the dimension spectrum).\\

Since we are interested in the relationship with chiral anomalies, we restrict our attention to spectral triples $(\Ac_0,H,D)$ of even degree, and suppose that $p$ is also an even integer. Decompose $H$ into the direct sum $H_+\oplus H_-$ corresponding to the eigenspaces for the eigenvalues $\pm 1$ of the chirality operator $\gamma_5$. According to this decomposition one can write in terms of $2\times 2$ matrices
$$
a=\left(\begin{matrix}
a_+ & 0 \\
0 & a_- \end{matrix} \right)\ ,\quad
D=\left(\begin{matrix}
0 & Q^* \\
Q & 0 \end{matrix} \right)\ ,\quad
F=\frac{D}{|D|}=\left(\begin{matrix}
0 & U^{-1} \\
U & 0 \end{matrix} \right)
$$
for any $a\in\Ac_0$, where $U:H_+\to H_-$ is the unitary operator $Q/(Q^*Q)^{1/2}$. The isomorphism of Hilbert spaces $H_+\cong H_-$ induced by $U$ allows to rewrite the spectral triple under the form of a $(p+1)$-summable Fredholm module:
$$
a=\left(\begin{matrix}
a_+ & 0 \\
0 & U^{-1}a_-U \end{matrix} \right)\ ,\quad
F=\left(\begin{matrix}
0 & 1 \\
1 & 0 \end{matrix} \right)\ ,\quad
|D|= \Id(\cc^2)\otimes (Q^*Q)^{1/2}\ .
$$
Let us deform the representation of $\Ac_0$ into $\Lc(H)$ using the homotopy $U_t=Q/(Q^*Q)^{t/2}$, $t\in[0,1]$. Then $U_1=U$, $U_0=Q$ and the above Fredholm module is homotopic to 
$$
\rho(a)=\left(\begin{matrix}
a_+ & 0 \\
0 & Q^{-1}a_-Q \end{matrix} \right)\ ,\quad
F=\left(\begin{matrix}
0 & 1 \\
1 & 0 \end{matrix} \right)\ ,\quad
|D|= \Id(\cc^2)\otimes (Q^*Q)^{1/2}\ .
$$
The new representation $\rho:\Ac_0\to\Lc(H)$ is no longer involutive but it is better suited for our purpose. Since $|D|^{-1}$ is in the Schatten ideal, one gets $[F,\rho(a)]\in \Lc^{p+1}(H)$ for any $a\in\Ac_0$. Hence the completion of $\Ac_0$ with respect to the norm $\|\rho(a)\|_{\infty}+\|[F,\rho(a)]\|_{p+1}$ is a Banach algebra $\Ac$, and the representation in Hilbert space extends to a continuous homomorphism $\rho:\Ac\to\Lc(H)$ such that $[F,\rho(a)]\in \Lc^{p+1}(H)$ for any $a\in\Ac$. In other words one obtains a $(p+1)$-summable quasihomomorphism of even parity from $\Ac$ to $\cc$
\be
\rho: \Ac \to \Ec^s\triangleright \Ic^s \label{qua}
\ee
with $\Ec=\Lc(H_+)$ and $\Ic=\Lc^{p+1}(H_+)$. The Riemann-Roch-Grothendieck theorem \ref{trr} gives a commutative diagram for the $K$-theory exact sequences 
\be
\vcenter{\xymatrix{
\Kt_{n+1}(\Ic\hotimes\Ac) \ar[r] \ar[d]^{\rho_!} & HC_{n-1}(\Ac) \ar[r] \ar[d]^{\ch^p(\rho)} & MK^{\Ic}_n(\Ac)  \ar[r] \ar[d]^{\rho_!}  & \Kt_n(\Ic\hotimes\Ac)  \ar[d]^{\rho_!}  \\
\Kt_{n+1-p}(\Ic) \ar[r]  & HC_{n-1-p}(\cc) \ar[r]  & MK^{\Ic}_{n-p}(\cc)  \ar[r]  & \Kt_{n-p}(\Ic)  }}\label{reg}
\ee
The topological $K$-theory of the Schatten ideals is known: $\Kt_0(\Ic)=\zz$ and $\Kt_1(\Ic)=0$. On the other hand, the cyclic homology of $\cc$ is $HC_{2k}(\cc)=\cc$ if $k\geq 0$, and vanishes in all other cases. Now two distinct cases are interesting. First take $n=p$ (even), then (\ref{reg}) yields
$$
\vcenter{\xymatrix{
\Kt_{1}(\Ic\hotimes\Ac) \ar[r] \ar[d]^{\rho_!} & HC_{p-1}(\Ac) \ar[r] \ar[d]^{\ch^p(\rho)} & MK^{\Ic}_p(\Ac)  \ar[r] \ar[d]^{\rho_!}  & \Kt_0(\Ic\hotimes\Ac)  \ar[d]^{\rho_!} &   \\
0 \ar[r]  & 0 \ar[r]  & \zz  \ar[r]  & \zz \ar[r] & 0 }}
$$
In this case the quasihomomorphism only induces an index map in topological $K$-theory $\rho_!:\Kt_0(\Ic\hotimes\Ac)\to\zz$, and it is a homotopy invariant for $\rho$. In particular if $e=e^*=e^2\in M_{N}(\Ac_0)$ is a projector, it determines a $K$-theory class $[e]\in \Kt_0(\Ic\hotimes\Ac)$ and the integer $\rho_!(e)$ is the Fredholm index of the compression of the Dirac operator $eDe$ on the Hilbert space $e(\cc^N\otimes H)$:
$$
\Ind(eDe)=\dim\ker (eDe)_+-\dim\ker (eDe)_-\ \in\zz\ .
$$
This index may be computed by evaluating the residue formula \ref{tres} on the Chern character of the idempotent $e$ in $HP_0(\Ac_0)$. One thus recovers the Connes-Moscovici local index theorem \cite{CM95}. A second interesting case is when $n=p+1$ (odd), so that (\ref{reg}) becomes
$$
\vcenter{\xymatrix{
 & \Kt_0(\Ic\hotimes\Ac) \ar[r] \ar[d]^{\rho_!} & HC_{p}(\Ac) \ar[r] \ar[d]^{\ch^p(\rho)} & MK^{\Ic}_{p+1}(\Ac)  \ar[r] \ar[d]^{\rho_!}  & \Kt_1(\Ic\hotimes\Ac)  \ar[d]^{\rho_!}  \\
0 \ar[r] & \zz \ar[r]  & \cc \ar[r]  & \cc^{\times}  \ar[r]  & 0  }}
$$
where $MK^{\Ic}_1(\cc)=\cc^{\times}$ is the multiplicative group of non-zero complex numbers and $\cc\to\cc^{\times}$ is the exponential map. In this case we obtain a regulator $\rho_!:MK^{\Ic}_{p+1}(\Ac)\to \cc^{\times}$. It is not invariant under homotopy but only under conjugation of quasihomomorphisms. A related construction was introduced by Connes and Karoubi in the context of algebraic $K$-theory \cite{CK}. We shall now establish the link with chiral anomalies \cite{P4}.\\

A non-commutative chiral gauge theory from the spectral triple $(\Ac_0,H,D)$ is constructed as follows. Consider the unbounded operator $Q:\dom(Q)\subset H_+\to H_-$ and let $\psi\in \dom(Q)$ and $\psib\in H_-^*$ represent some ``classical'' matter fields.  The equations of motion for $\psi,\psib$ are governed by the action functional
\be
S(\psi,\psib)=\langle \psib,Q\psi \rangle\ \in\cc\ .
\ee
Now let $g\in\Ac^+_0$ be a unitary element such that $g-1\in \Ac_0$. Recall that $\Ac_0^+$ denotes the algebra $\Ac_0$ with unit adjoined (eventually $g$ will be taken in the algebra $M_{\infty}(\Ac_0)^+$ of finite size matrices but for the moment we want to keep the notations simple). Its representation in the Hilbert space $H=H_+\oplus H_-$ reads $g=\bigl(\begin{smallmatrix}
g_+ & 0 \\
0 & g_- \end{smallmatrix} \bigr)$. We let $g$ act on $\psi,\psib$ by chiral gauge transformations:
$$
\psi\to g_+^{-1}\psi\ ,\qquad \psib\to \psib g_-\ .
$$
The action functional $S(\psi,\psib)$ is not invariant under gauge transformations; however the conjugate operator $g_-^{-1}Qg_+$ differs from $Q$ only by a bounded operator. We may absorb this change by adding to $Q$ another bounded operator $A:H_+\to H_-$ (gauge potential), which transforms according to the rule
$$
A\to A^g= ( g_-^{-1}Qg_+-Q) + g_-^{-1}Ag_+\ ,
$$
so that $(Q+A)\to g_-^{-1}(Q+A)g_+$ transforms according to the adjoint representation, and the full action
\be
S(\psi,\psib,A):= \langle \psib,(Q+A)\psi \rangle
\ee
is gauge-invariant. In practice the gauge potential $A$ is a finite linear combination of elements of the form $a_-(Qb_+-b_-Q)$, where $a,b\in\Ac_0$ and $a_{\pm},b_{\pm}$ are their representations on $H_{\pm}$, so that the product $Q^{-1}A$ lies in the Schatten ideal $\Ic=\Lc^{p+1}(H_+)$. Passing from classical to quantum field theory amounts to calculate, among other things, the partition function of the theory given by the functional integral \cite{IZ}
\be
Z(A)=\int d\psi\, d\psib\, e^{-S(\psi,\psib,A)}\ \in \cc\ .\label{Z}
\ee
If we decide to quantize $\psi$ and $\psib$ as fermions this purely formal object is defined as a Berezin integration over the Grassmann variables $\psi\in H_+$, $\psib\in H_-^*$ and normalized by $Z(0)=1$. Hence, the partition function (\ref{Z}) describes quantum fermionic fields $\psi,\psib$ moving in a non-quantized (classical) gauge potential $A$. The formal properties of Berezin integration \cite{IZ} show that $Z(A)$ may be defined as a renormalized determinant
$$
Z(A) = \Det(Q^{-1}(Q+A))\ .
$$
Renormalized determinants can be constructed in many different ways, one of the most popular being probably zeta-function renormalization \cite{S}. In any case it amounts to define a logarithm $W(A)=\ln Z(A)$, at least in a neighborhood of a fixed background gauge potential. Following physicists we shall expand $W(A)$ as a formal power series in $A$ represented by Feynman graphs, and then renormalize each graph. The advantage of this method is to show clearly the ambiguity inherent to renormalization, and the origin of anomalies. Moreover, the logarithm always makes sense as a formal power series. Hence using the naive relation $\ln\Det = \Tr\ln$, we try to define
$$
W(A)= \Tr \ln(1+Q^{-1}A) = \sum_{n\geq 1}\frac{(-)^{n+1}}{n}\, \Tr((Q^{-1}A)^n)\ .
$$
This expansion is depicted in terms of Feynman graphs as follows. We represent the propagator of the fields $\psi,\psib$ by an arrow $Q^{-1}=\xymatrix{*{}\ar@{-}[r] |-{\SelectTips{cm}{}\object@{>}}  & *{}}$ and the insertion of the gauge potential by a point $A=\bullet$. The products $Q^{-1}AQ^{-1}A\ldots$ are represented by chains, closed by taking the trace. $W(A)$ is thus given by a formal series of fermionic loops:\\
$$
W(A)  =\ \xymatrix{ *=0{\bullet} \ar@(ur,dr)@{-}[] |-{\SelectTips{cm}{}\object@{>}} }   
\ - \frac{1}{2}\  
\vcenter{\xymatrix@R=3pc{ *=0{\bullet} \ar@/^/@{-}[d] |-{\SelectTips{cm}{}\object@{>}} \\*=0{\bullet} \ar@/^/@{-}[u] |-{\SelectTips{cm}{}\object@{>}} } }
\ +\ \frac{1}{3}\ 
\vcenter{\xymatrix@C=1pc{ & *=0{\bullet} \ar@{-}[dr] |-{\SelectTips{cm}{}\object@{>}} &  \\
*=0{\bullet} \ar@{-}[ur] |-{\SelectTips{cm}{}\object@{>}} & & *=0{\bullet} \ar@{-}[ll] |-{\SelectTips{cm}{}\object@{>}} }}
\ -\frac{1}{4}\ 
\vcenter{\xymatrix@R=3pc@C=2.8pc{ *=0{\bullet} \ar@{-}[r] |-{\SelectTips{cm}{}\object@{>}} & *=0{\bullet} \ar@{-}[d] |-{\SelectTips{cm}{}\object@{>}} \\
*=0{\bullet} \ar@{-}[u] |-{\SelectTips{cm}{}\object@{>}} & *=0{\bullet} \ar@{-}[l] |-{\SelectTips{cm}{}\object@{>}} }}
\ +\frac{1}{5}\ 
\vcenter{\xymatrix@R=0.25pc@C=0.5pc{ & & *=0{\bullet} \ar@{-}[ddrr] |-{\SelectTips{cm}{}\object@{>}} & \\
 & & & & \\
*=0{\bullet} \ar@{-}[uurr] |-{\SelectTips{cm}{}\object@{>}} &  &  &  & *=0{\bullet} \ar@{-}[dddl] |-{\SelectTips{cm}{}\object@{>}} \\
 & & & &  \\
 & & & & \\
 & *=0{\bullet} \ar@{-}[uuul] |-{\SelectTips{cm}{}\object@{>}} & & *=0{\bullet} \ar@{-}[ll] |-{\SelectTips{cm}{}\object@{>}} & }}
\ +\ldots
$$
Recall that the operator $Q^{-1}A$ lies in $\Lc^{p+1}(H_+)$. Hence two remarks are in order:\\

\noindent 1) The $n$-th term $W^n(A)=\frac{(-)^{n+1}}{n}\, \Tr((Q^{-1}A)^n)$ is finite only when $n\geq p+1$. It follows that the first $p$ terms are divergent and need to be renormalized, for example by inserting some regularizing operator under the trace. We proceed as in section \ref{sren} and define the renormalized term $W^n_R(A)$ as the finite part of a zeta-function at zero: 
\be
W_R^n(A)=\frac{(-)^{n+1}}{n} \Pf_{z=0} \Tr((Q^{-1}A)^n|D|^{-2z}_+)\ ,
\ee
where $|D|^{-2z}_+=|Q^*Q|^{-z}$ is trace-class for $\re(z)$ sufficiently large, and the zeta-function $\Tr((Q^{-1}A)^n|D|^{-2z}_+)$ extends to a meromorphic function over the entire plane by hypothesis on the spectral triple. Obviously $W_R^n(A)=W^n(A)$ whenever $n\geq p+1$.\\

\noindent 2) Once renormalized, the series $W_R(A)=\sum_{n=0}^{\infty}W_R^n(A)$ converges if the operator norm verifies $\|Q^{-1}A\|_{\infty} <1$. In general this is not true and the series definitely diverges. It was of course expected since $W(A)$ is the expansion of the complex logarithm $\ln Z(A)$ around $A=0$. From now on we treat $W_R(A)$ exactly as a formal power series of the gauge potential $A$.\\

\noindent We are interested in the behaviour of the series $W_R(A)$ under chiral gauge transformations. In general it is not invariant:
\be
W_R(A^g)\neq W_R(A)\ ,
\ee
and the lack of invariance only comes from the renormalized terms $W_R^n(A)$ for $n\leq p$. In fact we are interested in ``smooth'' one-parameter families of gauge transformations, formalized as follows. Let $\Ac_0(0,1)$ be the algebraic tensor product of $\Ac_0$ with the algebra $\cinf(0,1)$ of smooth functions over $[0,1]$ whose values and all derivatives vanish at the endpoints \cite{P5}. Consider an element $g\in \Ac_0(0,1)^+$ in the unitalization of this algebra, with $g-1\in \Ac_0(0,1)$.  $g$ may be interpreted as a smooth loop of unitary elements in $\Ac^+_0$ (or $M_{\infty}(\Ac_0)^+$ if we work with matrix algebras), with basepoint $g_0=g_1=1$. Important examples are provided by idempotent loops: if $e\in\Ac_0$ is a projector, its associated loop $g=1+(\beta-1)e\in \Ac_0(0,1)^+$ is constructed from any invertible function $\beta\in\cinf(0,1)^+$ with winding number 1 (the Bott generator of the circle). Now given any unitary loop $g$, fix a background gauge potential $\underline{A}$ and consider the family of gauge-transforms
$$
A=\underline{A}^g=( g_-^{-1}Qg_+-Q) + g_-^{-1}\underline{A}\, g_+\ .
$$
$A$ is therefore a smooth loop in the space of bounded operators $H_+\to H_-$, and $Q^{-1}A$ is a loop in $\Lc^{p+1}(H_+)$. Then $W_R(A)$ is regarded as a formal series of smooth functions over the circle:
$$
W^n_R(A)\in \cinf(S^1)\quad \forall n\in\nn\ . 
$$
The chiral anomaly is the image of $W_R(A)$ under the differential $s:\cinf(S^1)\to \Om^1(S^1)$. At this point we make some sign conventions which differ from \cite{P4} but are more adapted to the physics literature: since $A$ and $Q$ are linar operators between $H_+$ and $H_-$, they have \emph{odd} degree. It means that $Q$ anticommutes with $s$ (as a constant odd-valued function over the circle), and the differential $sA$ is expressed by means of the Maurer-Cartan form $\om_{\pm}=g_{\pm}^{-1}s g_{\pm}$ and the structure equation (compare with BRS transformations \cite{MSZ}):
\be
s A= -\om_-(Q+A)-(Q+A)\om_+\ .
\ee
Write $q\om=\om_+ + Q^{-1}\om_-Q$. Then in any degree $n\geq p+1$, using the cyclicity of the trace the boundary of $W_R^n(A)=W^n(A)$ reads
\beq
s W^n(A) &=& -(-)^{n+1}\Tr (Q^{-1}s A (Q^{-1}A)^{n-1})\non\\
&=& (-)^{n+1}\Big(\Tr \big(q\om (Q^{-1}A)^{n-1}\big) + \Tr \big(q\om (Q^{-1}A)^{n}\big)\Big)\non
\eeq
Hence, in the sense of formal power series in $A$ the terms $\Tr (q\om (Q^{-1}A)^{n})$ cancel each other in the infinite sum $\sum_{n>p}sW^n(A)$ due to the alternate signs. For $n\leq p$ however, the regularizing operator $|D|^{-2z}_+$ destroys the cyclicity of the trace and cancellation does no longer holds. Using the abstract pseudodifferential calculus, one sees that the renormalized trace applied to a commutator yields residues of zeta-functions:
$$
\Pf_{z=0}\Tr([x,y]|D|^{-2z}_+)=- \sum_{k\geq 1} \frac{(-)^k}{k!} \res \frac{\Gamma(z+k)}{\Gamma(z+1)}\, \Tr(x\, y^{(k)}|D|^{-2(z+k)}_+)\ ,
$$
where as usual $y^{(k)}$ denotes the $k$-th power the derivation $[|D|_+^2,\ ]$ on $y$. This is actually a finite sum because the zeta-functions are holomorphic at $z=0$ for $k$ large enough. Hence the boundary of the term $W_R^n(A)$ for $n\leq p$ reads
\beq
s W_R^n(A) &=& (-)^{n+1}\Pf_{z=0}\Big(\Tr \big(q\om (Q^{-1}A)^{n-1}|D|^{-2z}_+\big) + \Tr \big(q\om (Q^{-1}A)^{n}|D|^{-2z}_+\big)\Big)\non\\
&& +   \frac{(-)^{n+1}}{n}\sum_{i=0}^{n-1}\Pf_{z=0}\Big(\Tr\big([(Q^{-1}A)^i,q\om(Q^{-1}A)^{n-i-1}]|D|^{-2z}_+\big)\non\\
&& \qquad\qquad\qquad + \Tr\big([(Q^{-1}A)^{i+1},\om_+(Q^{-1}A)^{n-i-1}]|D|^{-2z}_+\big)\non\\
&& \qquad\qquad\qquad + \Tr\big([(Q^{-1}A)^i,Q^{-1}\om_-Q(Q^{-1}A)^{n-i}]|D|^{-2z}_+\big)\Big)\ ,\non
\eeq
the last three terms being sums of residues. We view $s W_R^n(A)$ as a polynomial in $A$ and $\om$, so that the infinite sum $\sum_{n\geq 0} s W_R^n(A)$ makes sense as a formal power series. Consequently, the terms $\Tr (q\om (Q^{-1}A)^{n}|D|^{-2z}_+)$ still cancel each other and only the first term $\Tr(q\om|D|^{-2z}_+)$ plus the residues remain. The anomaly $\Delta(\om,A)=s W_R(A)$ is thus given by the \emph{finite sum} of one-forms
\be
\Delta(\om,A)=\Pf_{z=0} \Tr\big(q\om|D|^{-2z}_+\big) + \ \mbox{residues}\ \in \Om^1(S^1)\ ,
\ee
and it is a polynomial in $A$ of degree at most $p$. This may also be depicted diagramatically, for example in the case $p=2$, as follows:
$$
\vcenter{\xymatrix@!0@=2.5pc{
W_R(A) \ar[d]_s & = &  & W^1_R(A) \ar[dl] \ar[dr]  & + & W^2_R(A) \ar[dl] \ar[dr] & + & W^3(A) \ar[dl] \ar[dr] & + & W^4(A) \ar[dl] \ar[dr] & +  \ldots  \\
\Delta(\om,A) & = & \Delta^0(\om,A) & + & \Delta^1(\om,A) & + & \Delta^2(\om,A) & + & 0 & + & \ 0  \ldots }}
$$
Each $\Delta^i(\om,A)$ is a homogeneous polynomial of degree $i$ in $A$. Taking the quotient by the image of $s$, we may view the anomaly as a de Rham cohomology class of degree one $\Delta(\om,A)\in H^1(S^1)$. It is non-trivial in general because it is the boundary of a formal power series (and not a convergent sum) in $\cinf(S^1)$. Intuitively, the lack of gauge-invariance of $W_R(A)$ originates from the overlapping of the vertices $\bullet$ in the renormalized Feynman graphs, whence the ``locality'' of the anomaly. It must also be stressed that renormalization is not unique: there are different ways to extract finite quantities $W_R^n(A)$ from the regularizing operator $|D|^{-2z}_+$, and this may add finitely many local counterterms (residues) to the series $W_R(A)$. Consequently the cocycle $s W_R(A)$ will only change by a finite sum of boundaries, so that the cohomology class of the anomaly remains unchanged.\\ 
The formula $\Delta(\om,A)=sW_R(A)$ is completely analogous to the construction of the residue Chern character $\chi_R=[\d,\eta_R]$ from the boundary of the renormalized eta-cochain (section \ref{sren}). The formal power series $W_R(A)$ is simply replaced by the improper cochain $\eta_R$. This explains why Theorem \ref{tres} is called an anomaly formula. It is moreover possible to exhibit an explicit correspondence between $W_R(A)$ and the renormalized eta-cochain associated to the quasihomomorphism (\ref{qua}). The latter has even parity, hence the renormalized eta-cochain (\ref{ren}) is odd. Since the $X$-complex of the quasi-free algebra $\Rc=\cc$ reduces to $X(\cc):\cc\rightleftarrows 0$, the only non-zero components $\etah_{R0}^{2n+1}\rho_*:\Om^{2n+1}\Th\Ac_0\to \cc$ are explicitly given by
$$
\etah^{2n+1}_{R0}(x_0\dd x_1\ldots\dd x_{2n+1}) =  \frac{\Gamma(n+1)}{(2n+1)!}  \frac{1}{2}\Pf_{z=0}\tau(F x_0[F,x_1]\ldots [F,x_{2n+1}]|D|^{-2z})
$$
where $\tau$ is the supertrace of operators on $H$, and $F$ admits the matrix representation $F=\bigl(\begin{smallmatrix}
0 & 1 \\
1 & 0 \end{smallmatrix} \bigr)$. The transgressed Chern character $\tch^{2n+1}_R(\rho)=\etah_R^{2n+1}\rho_*\gamma$ is the composite with the Goodwillie equivalence, and yields a linear map $X(\Th\Ac_0)\to\cc$ of odd degree. As in section \ref{sbiv} we evaluate it on the Chern character of a unitary element $g\in M_{\infty}(\Ac_0)^+$ such that $g-1\in M_{\infty}(\Ac_0)$. Using the linear inclusion $\Ac_0\hookrightarrow \Th\Ac_0$ lift $g$ to an invertible element $\gh\in M_{\infty}(\Th\Ac_0)^+$. Its Chern character $\ch_1(\gh)=\frac{1}{\sqrt{2\pi i}}\Tr\nat\gh^{-1}\dd\gh$ defines an odd cycle of the complex $X(\Th\Ac_0)$.

\begin{lemma}\label{lV}
The transgressed Chern character $\tch^{2n+1}_R(\rho): X(\Th\Ac_0)\to \cc$ evaluated on $\ch_1(\gh)=\frac{1}{\sqrt{2\pi i}}\Tr\nat\gh^{-1}\dd\gh$ reads
$$
\tch^{2n+1}_R(\rho)\cdot \ch_1(\gh)=\frac{1}{\sqrt{2\pi i}}\frac{(n!)^2}{(2n+1)!} \frac{1}{2} \Pf_{z=0}\tau(FV^{2n+1}|D|^{-2z}) \ ,
$$
where $V$ is the compact operator $\rho(g)^{-1}[F,\rho(g)]$ on $H$.
\end{lemma}
{\it Proof:} By definition $\tch_R^{2n+1}(\rho)=\etah_R^{2n+1}\rho_*\gamma$. In \cite{P5} we calculated the image of the Chern character $\ch_1(\gh)=\frac{1}{\sqrt{2\pi i}}\Tr\nat\gh^{-1}\dd\gh$ under the Goodwillie equivalence:
$$
\gamma\ch_1(\gh)=\frac{1}{\sqrt{2\pi i}}\sum_{n\geq 0}(-)^n n!\, \Tr(\gh^{-1}\dd\gh(\dd\gh^{-1}\dd\gh)^n)\ .
$$
Since $\rho_*(\gh)=\rho(g)\in\Lc(H)$, the evaluation of the above non-commutative form on $\etah_{R0}^{2n+1}\rho_*$ is straightforward.  \cqfd\\

For $2n\geq p$ the operator $V^{2n+1}$ is trace-class and one recovers the non-renormalized formula (\ref{nr})
$$
\tch^{2n+1}(\rho)\cdot \ch_1(\gh) = \frac{1}{\sqrt{2\pi i}} \frac{(n!)^2}{(2n+1)!} \, \frac{1}{2}\tau (FV^{2n+1})\ .
$$
Now fix the background connection $\underline{A}=0$; for any $g\in M_{\infty}(\Ac_0)^+$ the gauge transform $A=\underline{A}^g$ is given by
$$
A=g_-^{-1}Qg_+-Q\ ,
$$
and $1+Q^{-1}A$ is invertible. The operator $V$ can be expanded in a formal power series in $A$ with no term of degree zero. Indeed one has $F= \bigl( \begin{smallmatrix} 0 & 1 \\ 1 & 0 \end{smallmatrix} \bigr)$ and $\rho(g)= \bigl( \begin{smallmatrix} g_+ & 0 \\ 0 & Q^{-1}g_- Q \end{smallmatrix} \bigr)$, so that
$$
V=\left(\begin{matrix}
0 & (1+Q^{-1}A)^{-1}-1 \\
Q^{-1}A & 0 \end{matrix} \right)
$$
and $(1+Q^{-1}A)^{-1}-1=\sum_{k\geq 1} (-Q^{-1}A)^k$. The renormalized eta-cochain reads
\beq
\lefteqn{\tch^{2n+1}_R(\rho)\cdot \ch_1(\gh) = }\non\\
&& \qquad \frac{(-)^n}{\sqrt{2\pi i}} \frac{(n!)^2}{(2n+1)!} \Pf_{z=0} \Tr \left( \left(\frac{Q^{-1}A}{1+Q^{-1}A}\right)^{2n+1}(1+Q^{-1}A/2) |D|_+^{-2z} \right) \ , \non
\eeq
and its expansion as a formal power series in $A$ contains only terms of degree $\geq 2n+1$. Consequently, the infinite sum 
\be
\tch_R(\rho)\cdot \ch_1(\gh) = \sum_{n\geq 0} \tch^{2n+1}_R(\rho)\cdot \ch_1(\gh)
\ee
makes sense as a formal power series in $A$. The following lemma establishes the link with the partition function of the non-commutative gauge theory.

\begin{lemma}
Let $g\in M_{\infty}(\Ac_0)^+$ be a unitary such that $g-1\in M_{\infty}(\Ac_0)$ and take the gauge potential $A=g_-^{-1}Qg_+-Q$. The renormalized partition function $W_R(A)$ coincides with the formal power series $\tch_R(\rho)\cdot \ch_1(\gh)$ modulo a counterterm:
\be
W_R(A) + P(A)=\sqrt{2\pi i}\, \tch_R(\rho)\cdot \ch_1(\gh) \ ,
\ee
where $P(A)$ is a polynomial in $A$ given by a finite sum of residues, with $V=\bigl(\begin{smallmatrix}
0 & (1+Q^{-1}A)^{-1}-1 \\
Q^{-1}A & 0 \end{smallmatrix} \bigr)$:
$$
P(A)=\frac{1}{8}\sum_{n\geq 0}\frac{(n!)^2}{(2n+1)!}\sum_{k\geq 1} \frac{(-)^k}{k!} \res \frac{\Gamma(z+k)}{\Gamma(z+1)}\, \tau (V^{2n+1}V^{(k)}|D|^{-2(z+k)})\ .
$$
\end{lemma}
{\it Proof:} Note that the operator $V=\rho(g)^{-1}[F,\rho(g)]$ fulfills the flat curvature condition $[F,V]+V^2=0$. Our first task will be the computation of $\ln(1+FV)$ in the following terms,
$$
\ln(1+FV)=\sum_{n\geq 0}\frac{(n!)^2}{(2n+1)!}\, (FV^{2n+1}+\frac{1}{2} V^{2n+2})\ ,
$$
understood as an equality of formal power series in $A$. To this end, for any $t\in [0,1]$ define $F_{tV}=F+tV$. One has $F_{tV}^2=1+(t^2-t)V^2$, hence $F_{tV}$ is invertible in the sense of formal power series. Consider the (formal) integral
$$
I=\int_0^1 dt\int_0^{\infty}du\, F_{tV}e^{-uF_{tV}^2}V\ .
$$
We shall compute $I$ by two different manners. First, integrate over $u$:
$$
I=\int_0^1dt\, F_{tV}(F_{tV}^2)^{-1}V=\int_0^1dt\, F_{tV}^{-1}V\ .
$$
The inverse of $F_{tV}$ is $(F+tV)^{-1}=(F(1+tFV))^{-1}=(1+tFV)^{-1}F$ because $F^2=1$. Hence
$$
I=\int_0^1dt\, (1+tFV)^{-1}FV = \sum_{n\geq 1}\frac{(-)^{n+1}}{n}\, (FV)^n =\ln(1+FV)
$$
by expanding the series $(1+tFV)^{-1}$ in powers of $t$. On the other hand, we could first expand the exponential of $F_{tV}^2=1+(t^2-t)V^2$ and then integrate over $u$:
$$
I=\int_0^1 dt\int_0^{\infty}du\, F_{tV}\sum_{n\geq 0}\frac{u^n}{n!}e^{-u}(t^2-t)^nV^{2n+1}=\sum_{n\geq 0}\int_0^1dt\, (t^2-t)^nF_{tV}V^{2n+1}\ .
$$
Finally, it remains to use the integrals
$$
\int_0^1 (t^2-t)^ndt= \frac{(n!)^2}{(2n+1)!}\ ,\qquad \int_0^1 t(t^2-t)^ndt=\frac{1}{2}\frac{(n!)^2}{(2n+1)!}
$$
to express $I$ as wanted:
$$
I=\sum_{n\geq 0}\frac{(n!)^2}{(2n+1)!}\, (FV^{2n+1}+\frac{1}{2} V^{2n+2})\ .
$$
One has $1+FV=\bigl(\begin{smallmatrix}
1+Q^{-1}A & 0 \\
0 & (1+Q^{-1}A)^{-1} \end{smallmatrix} \bigr)$. Let us compute the supertrace
\beq
\lefteqn{\frac{1}{2} \tau \ln\big((1+FV)|D|^{-2z}\big) =  \frac{1}{2} \tau \ln \left(\begin{smallmatrix}
1+Q^{-1}A & 0 \\
0 & (1+Q^{-1}A)^{-1} \end{smallmatrix} \right)|D|^{-2z}}\non\\
&& =\frac{1}{2} \Tr\, \ln(1+Q^{-1}A)|D|_+^{-2z} - \frac{1}{2} \Tr\, \ln(1+Q^{-1}A)^{-1}|D|_+^{-2z}\non\\
&&=  \Tr\, \ln(1+Q^{-1}A)|D|^{-2z}_+\ ,\non
\eeq
whose finite part at $z=0$ is precisely the renormalized logarithm $W_R(A)$. Hence we deduce the equality of formal power series in $A$:
$$
W_R(A)=\frac{1}{2}  \sum_{n\geq 0}\frac{(n!)^2}{(2n+1)!} \Pf_{z=0} \tau\big((FV^{2n+1}+\frac{1}{2} V^{2n+2})|D|^{-2z}\big)\ .
$$
The terms involving the even powers $V^{2n+2}$ are in fact residues. Indeed, $V$ is of odd parity so that $\tau(V^{2n+2}|D|^{-2z})=-\tau(V^{2n+1}|D|^{-2z}V)$ by cyclicity of the supertrace, and we can write
$$
\tau(V^{2n+2}|D|^{-2z})= -\frac{1}{2} \tau(V^{2n+1}[|D|^{-2z},V])\ .
$$
Using the asymptotic expansion of the commutator, one thus has
$$
\Pf_{z=0}\tau(V^{2n+2}|D|^{-2z})=-\frac{1}{2}\sum_{k\geq 1} \frac{(-)^k}{k!} \res\frac{\Gamma(z+k)}{\Gamma(z+1)}\, \tau(V^{2n+1}V^{(k)}|D|^{-2(z+k)})\ .
$$
This is a finite sum because the residues vanish whenever $k+2n\geq p$. Therefore, Lemma \ref{lV} implies
\beq
\lefteqn{W_R(A)=\sqrt{2\pi i}\, \tch_R(\rho)\cdot \ch_1(\gh) -\frac{1}{8}\sum_{n\geq 0}\frac{(n!)^2}{(2n+1)!}\sum_{k\geq 1} \frac{(-)^k}{k!}\times}\non\\
&&\qquad\qquad\qquad \res \frac{\Gamma(z+k)}{\Gamma(z+1)}\, \tau (V^{2n+1}V^{(k)}|D|^{-2(z+k)})\ ,\non
\eeq
where the sum of residues is finite and yields a polynomial in $A$.  \cqfd\\

The above lemma means that the formal series $\sqrt{2\pi i}\, \tch_R(\rho)\cdot \ch_1(\gh)$ is nothing else but another choice of renormalization for $W(A)$. The polynomial $P(A)$ is simply a counterterm. Now, if $g$ is a loop of unitaries in the matrix algebra $M_{\infty}(\Ac_0)^+$ with basepoint the identity, the cohomology class of the anomaly is represented by 
\be
\Delta(\om,A) \equiv \sqrt{2\pi i}\, s(\tch_R(\rho)\cdot\ch_1(\gh)) \ \in H^1(S^1)\ ,
\ee
where the invertible loop $\gh$ in $M_{\infty}(\Th\Ac_0)^+$ is the canonical lift of $g$ obtained by the linear inclusion $\Ac_0\hookrightarrow \Th\Ac_0$. Replacing the universal derivative $\dd$ by $s$ and the supertrace $\tau\nat$ by $\tau$, equation (\ref{infsum}) reinterpreted as an equality of formal power series in $A$ amounts to
$$
\ch_R(\rho)\cdot \frac{\Tr(\gh^{-1}s\gh )}{\sqrt{2\pi i}} = s(\tch_R(\rho)\cdot\ch_1(\gh))\ ,
$$
whence the identification of cohomology classes $\Delta(\om,A)\equiv \ch_R(\rho)\cdot \Tr(\gh^{-1}s\gh)$. Moreover the left-hand-side is a polynomial in $A$ and depends linearly on $\om$. Using Bott periodicity one thus gets

\begin{corollary}\label{cta}
Let $(\Ac_0,H,D)$ be a regular finitely summable spectral triple of even parity, and let $W_R(A)$ be any renormalization of the partition function for the associated chiral gauge theory. Consider a projector $e\in M_{\infty}(\Ac_0)$, and the corresponding idempotent loop $g=1+(\beta-1)e \in M_{\infty}(\Ac_0)(0,1)^+$. Then the anomaly $\Delta(\om,A)=s W_R(A)$ integrated over the loop of gauge potentials  $A=g_-^{-1}Qg_+-Q$ computes the index
\be
\Ind(eDe)= \frac{1}{2\pi i}\oint \Delta(\omega,A) \ \in\zz\ .
\ee
The latter does not depend on the choice of renormalization and is given by a local formula (e.g. a sum of residues of zeta-functions). 
\end{corollary} 
{\it Proof:} Modulo a coboundary, one has $\Delta(\om,A)\equiv \ch_R(\rho)\cdot\Tr(\gh^{-1}s\gh )$ with the canonical lifting $\gh\in M_{\infty}(\Th\Ac_0)(0,1)^+$ of $g$. Hence integrating the anomaly over the circle one gets
$$
\frac{1}{2\pi i}\oint \Delta(\omega,A)= \frac{1}{2\pi i}\,\ch_R(\rho)\cdot \int_0^1\Tr(\gh^{-1}s\gh):=\frac{1}{\sqrt{2\pi i}}\, \ch_R(\rho)\cdot\cs_0(\gh)\ .
$$
Now Bott periodicity (Lemma 4.4 of \cite{P5}) implies the equivalence of homology classes $\cs_0(\gh)\equiv \sqrt{2\pi i}\, \ch_0(\eh)$ in the complex $X(\Th\Ac_0)$, where $\eh\in M_{\infty}(\Th\Ac_0)$ is an idempotent lifting of $e$. Hence one finds
$$
\frac{1}{2\pi i}\oint \Delta(\omega,A)=\ch_R(\rho)\cdot\ch_0(\eh)=\rho_!(e)
$$
which is nothing else but the index $\Ind(eDe)\in\zz$. \cqfd\\

So far the formal power series $W_R(A)$ allowed us to interpret the anomaly formula, section \ref{sren}, but the renormalized determinant $Z(A)=\Det(1+Q^{-1}A)$ is still not defined. Morally, $Z(A)$ should be the exponential of $W_R(A)$ but we want a complex number, not a formal series in $A$.  We shall construct a renormalized determinant $Z_R(A)$ by integrating the anomaly along a smooth path of gauge potentials $A=g_-^{-1}Qg_+-Q$, starting from the background $\underline{A}=0$:
$$
\Delta(\om,A)=\frac{sZ_R(A)}{Z_R(A)}\ .
$$
Note that it is not unique, because the anomaly itself $\Delta(\om,A)=s W_R(A)$ is defined modulo the boundary of a finite sum of counterterms. We could as well take the anomaly as the boundary of the formal series $\sqrt{2\pi i}\,\tch_R(\rho)\cdot \ch_1(\gh)$, the resulting new determinant would differ from the former only by a phase factor involving the counterterms. Let us define $\Ac_0[0,1]$ as the algebraic tensor product of $\Ac_0$ with the space $\cinf[0,1]$ of smooth functions over $[0,1]$ such that all derivatives vanish at the endpoints (while the functions themselves take arbitrary values at $0$ and $1$). We say that a unitary element $g\in M_{\infty}(\Ac_0)^+$ is homotopic to $1$ if there exists a unitary path $u\in M_{\infty}(\Ac_0)[0,1]^+$ such that $u_0=1$ and $u_1=g$. The renormalized determinant of $g$ is the complex number
\be
\Det_R(g) := \exp(\sqrt{2\pi i}\,\ch_R(\rho)\cdot \cs_0(\uh))\  \in \cc^{\times}\ ,
\ee
where the Chern-Simons form of the canonical lifting $\uh\in M_{\infty}(\Th\Ac_0)[0,1]^+$ is defined as in \cite{P5}:
$$
\cs_0(\uh)=\frac{1}{\sqrt{2\pi i}}\int_0^1 \Tr(\uh^{-1}s\uh) \ \in \Th\Ac_0\ .
$$

\begin{lemma}
The renormalized determinant $\Det_R(g)$ does not depend on the unitary path $u\in M_{\infty}(\Ac_0)[0,1]^+$ joining $1$ to $g$, and its logarithmic derivative coincides with the anomaly
$$
s(\sqrt{2\pi i}\,\tch_R(\rho)\cdot \ch_1(\gh))=\frac{s\Det_R(g)}{\Det_R(g)}\ .
$$
\end{lemma}
{\it Proof:} First we show the independence of $\Det_R$ with respect to the path $u$. If $u'\in M_{\infty}(\Ac_0)[0,1]^+$ is such another path, we form a loop $u''$ with basepoint $1$, by first going along $u$ up to $g$ and then backward along $u'$ up to $1$. Hence one has the equality of cycles $\cs_0(\uh'')=\cs_0(\uh)-\cs_0(\uh')$ in $X(\Th\Ac_0)$, and the evaluation of $\cs_0(\uh'')$ on the renormalized Chern character $\ch_R(\rho):X(\Th\Ac_0)\to\cc$ coincides with its evaluation on the $p$-th dimensional Chern character $\ch^p(\rho)$. But $\ch^p(\rho)$ extends to a cyclic cocycle over the completion $\Ac$ of $\Ac_0$. Now if we regard $u''$ as a loop of invertible matrices over the Banach algebra $\Ac^+$, by Bott periodicity it defines a $K$-theory class $[e]-[p_0]\in \Kt_0(\Ac)$ for some idempotent $e\in M_{\infty}(\Ac^+)$, and by Lemma 4.4 of \cite{P5} one gets
$$
\cs_0(\uh'')=\cs_0(\uh)-\cs_0(\uh')\equiv \sqrt{2\pi i}\, \ch_0(\eh)\ \in HP_0(\Ac)\ .
$$
Therefore evaluation of these cyclic homology classes on $\ch^p(\rho)$ yields the numbers
$$
\sqrt{2\pi i}\, \ch^p(\rho)\cdot(\cs_0(\uh)-\cs_0(\uh')) = 2\pi i\, \ch^p(\rho)\cdot\ch_0(\eh)\ \in 2\pi i\,\zz\ ,
$$
so that the exponential of $\sqrt{2\pi i}\, \ch_R(\rho)\cdot\cs_0(\uh)$ is independent of $u$.\\
Now compute the logarithmic derivative of the determinant:
$$
\frac{s\Det_R(g)}{\Det_R(g)} = \sqrt{2\pi i}\, s(\ch_R(\rho)\cdot\cs_0(\uh)) = \ch_R(\rho)\cdot \Tr(\gh^{-1}s\gh)\ .$$
We conclude from $\ch_R(\rho)\cdot \Tr(\gh^{-1}s\gh)= \sqrt{2\pi i}\, s(\tch_R(\rho)\cdot\ch_1(\gh))$. \cqfd\\

Once we take $\Det_R(g)$ as a definition of the renormalized determinant, the integral of the anomaly
$$
\Ind(eDe)= \frac{1}{2\pi i}\oint \Delta(\omega,A)
$$
just counts the winding number of the partition function $Z_R(A)\in \cc^{\times}$ around the idempotent loop. This proves the result announced in \cite{P4} and interprets the topological index $\rho_!:\Kt_0(\Ac)\to\zz$ as a topological anomaly. But the renormalized determinant allows to interpret the regulator $\rho_!:MK^{\Ic}_{p+1}(\Ac)\to\cc^{\times}$ as well. Recall from \cite{P5} that a multiplicative $K$-theory class of degree $p+1$ over $\Ac$ is represented by a pair $(\gh,\te)$ where $\gh$ is an invertible element of $(\Ic\hotimes\Th\Ac)^+$ and $\te$ is a chain of even degree in the quotient complex $X_p(T\Ac,J\Ac)= X(T\Ac)/F^p_{J\Ac}X(T\Ac)$ satisfying the transgression relation $\ch_1(\gh)=\nat\dd\te$. Since we work with a dense subalgebra $\Ac_0\subset\Ac$ we shall restrict to pairs $(\gh,\te)$ with $\gh\in M_{\infty}(\Th\Ac_0)^+$ and $\te\in X_p(T\Ac,J\Ac)$ in the image of $\Th\Ac_0$:

\begin{corollary}\label{cdet}
Let $(\Ac_0,H,D)$ be a regular spectral triple of even parity and analytic dimension $p$ (even). Let the pair $(\gh,\te)\in MK^{\Ic}_{p+1}(\Ac)$ represent a multiplicative $K$-theory class, where $\gh\in M_{\infty}(\Th\Ac_0)^+$ is the canonical lift of a unitary element $g\in M_{\infty}(\Ac_0)^+$ homotopic to $1$ and $\te$ is in the image of $\Th\Ac_0$. Then its evaluation under the regulator map $\rho_!:MK^{\Ic}_{p+1}(\Ac)\to\cc^{\times}$ is the complex number
\be
\rho_!(\gh,\te)= \exp(\sqrt{2\pi i}\,\ch_R(\rho)\cdot \te)\,\Det_R^{-1}(g) \ \in \cc^{\times}\ . \label{phase}
\ee
\end{corollary}
{\it Proof:} By hypothesis, there exists a unitary path $u\in M_{\infty}(\Ac_0)[0,1]^+$ joining $u_0=1$ to $u_1=g$. Then the pair $(\gh,\te)$ is equivalent to $(1,\te-\cs_0(\uh))$ in $MK^{\Ic}_{p+1}(\Ac)$, with $\uh\in M_{\infty}(\Th\Ac_0)[0,1]^+$ the canonical lifting of $u$. By definition (see \cite{P5})
$$
\rho_!(1,\te-\cs_0(\uh))= \big(1,\ch^p(\rho)\cdot (\te-\cs_0(\uh))\big)\ \in MK^{\Ic}_1(\cc)\ .
$$
Since the analytic dimension is $p$ we have the equality of non-periodic cyclic cohomology classes $\ch_R(\rho)\equiv \ch^p(\rho)$ in $HC^p(\Ac_0)$. It means that the $p$-th dimensional Chern character $\ch^p(\rho)$, once restricted to the image of $\Th\Ac_0$ in the complex $X_p(T\Ac,J\Ac)$, is cohomologous to the renormalized Chern character $\ch_R(\rho)$:
$$
\ch_R(\rho)=\ch^p(\rho)+  \sum_{n=0}^{p/2-1}\tch_R^{2n+1}(\rho) \circ\d
$$
with the boundary map $\d=\nat\dd$ on $X_p(T\Ac,J\Ac)$. Therefore
$$
\ch^p(\rho)\cdot (\te-\cs_0(\uh))= \ch_R(\rho)\cdot ( \te-\cs_0(\uh)) - \sum_{n=0}^{p/2-1}\tch_R^{2n+1}(\rho)\cdot \nat\dd(\te-\cs_0(\uh))\ .
$$
By hypothesis $\nat\dd\te=\ch_1(\gh)$ in $X_p(T\Ac,J\Ac)$, and Lemma 4.3 of \cite{P5} implies $\nat\dd\cs_0(\uh)=\ch_1(\gh)$, whence
$$
\rho_!(\gh,\te)= \big(1,\ch_R(\rho)\cdot(\te- \cs_0(\uh))\big)
$$
at the level of multiplicative $K$-theory classes. Finally, the group isomorphism $MK^{\Ic}_1(\cc)\cong \cc^{\times}$ identifies a pair $(1,\al)$ with the complex number $\exp(\sqrt{2\pi i}\, \al)$, for any $\al\in \cc$. Thus
\beq
\rho_!(\gh,\te) &=& \exp\big(\sqrt{2\pi i}\, \ch_R(\rho)\cdot(\te- \cs_0(\uh))\big) \non\\
&=& \exp(\sqrt{2\pi i}\, \ch_R(\rho)\cdot\te)\,\Det_R^{-1}(g) \non
\eeq
as claimed.  \cqfd\\
 
Remark, as it should be, that the complex number $\rho_!(\gh,\te)$ does not depend on the renormalization chosen for the Chern character $\ch_R(\rho)$, because the phase factor $\exp(\sqrt{2\pi i}\,\ch_R(\rho)\cdot \te)$ automatically absorbs any change in the renormalized determinant. This was expected since the regulator map $MK^{\Ic}_{p+1}(\Ac)\to \cc^{\times}$ is canonically defined. The fact that it is a morphism of abelian groups allows also to interpret the \emph{multiplicative anomaly}: for two gauge transformations $g,h \in M_{\infty}(\Ac_0)^+$ the renormalized determinant fails to be multiplicative:
\be
\Det_R(gh) \neq \Det_R(g) \Det_R(h)\ .
\ee
This lack of multiplicativity is compensated by the phase factor in (\ref{phase}). Using the equivalence relation on multiplicative $K$-theory yields an explicit formula. Let $\gh$, $\hh$ and $\widehat{gh} \in M_{\infty}(\Th\Ac_0)^+$ be the canonical liftings of $g$, $h$ and $gh$ respectively. The difference $\widehat{gh} - \gh\hh$ thus lies in the pro-nilpotent ideal $M_{\infty}(\Jh\Ac_0) \subset M_{\infty}(\Th\Ac_0)$, the product $ \hh^{-1}\gh^{-1}\widehat{gh} $ lies in $1 + M_{\infty}(\Jh\Ac_0)$ and its logarithm given by the usual power series defines an element
$$
\ln (\hh^{-1}\gh^{-1}\widehat{gh} ) \in M_{\infty}(\Jh\Ac_0)\ .
$$

\begin{corollary}\label{cma}
Let $(\Ac_0,H,D)$ be a regular spectral triple of even parity and analytic dimension $p$, with renormalized Chern character $\ch_R(\rho):\Th\Ac_0\to\cc$. Let $g,h\in M_{\infty}(\Ac_0)^+$ be two unitary elements homotopic to $1$. The multiplicative anomaly is given by the formula
\be
\Det_R(gh) = \exp\big( \ch_R(\rho)\cdot \Tr\ln (\hh^{-1}\gh^{-1} \widehat{gh} ) \big) \, \Det_R(g) \Det_R(h)\ .
\ee
\end{corollary}
{\it Proof:} Because $g$ and $h$ are homotopic to $1$ there exist chains $\te$ and $\te'$ in the complex $X_p(T\Ac,J\Ac)$ such that $(\gh,\te)$ and $(\hh,\te')$ represent multiplicative $K$-theory classes in $MK^{\Ic}_{p+1}(\Ac)$. Their sum is represented by the pair 
$$
(\gh,\te)+(\hh,\te') = (\gh\oplus\hh,\te+\te') \ .
$$
Using rotation $2\times 2$ matrices, one shows as in \cite{P5} that $(\gh\oplus\hh,\te+\te')$ and the pair $(\gh \hh,\te+\te')$ represent the same multiplicative $K$-theory class. Now consider the linear homotopy
$$
\uh_t = t\, \widehat{gh} + (1-t) \gh\hh = \gh\hh +t (\widehat{gh} - \gh\hh) 
$$
for $t\in [0,1]$. Since $\widehat{gh} - \gh\hh$ lies in a pro-nilpotent ideal, $\uh$ defines an invertible path in $M_{\infty}(\Th\Ac_0)^+$ joining $\uh_0=\gh\hh$ to $\uh_1= \widehat{gh} $. Its inverse is given by the series
$$
\uh^{-1}_t = \sum_{n\geq 0} (-t)^n (\hh^{-1}\gh^{-1}\widehat{gh} - 1 )^n \hh^{-1}\gh^{-1} \ .
$$
The Chern-Simons form associated to $\uh$ thus reads (to be correct one should transform the variable $t$ by a smooth bijective function $[0,1]\to [0,1]$ with all derivatives vanishing at the endpoints)
\beq
\sqrt{2\pi i} \, \cs_0(\uh) &=& \int_0^1 \Tr(\uh^{-1}s\uh)\ =\ \sum_{n\geq 0}\int_0^1 dt (-t)^n \Tr((\hh^{-1}\gh^{-1}\widehat{gh} - 1 )^{n+1})  \non\\
&=& \Tr\ln(\hh^{-1}\gh^{-1}\widehat{gh}) \non
\eeq
and the pair $(\gh\hh,\te+\te')$ is equivalent to $(\widehat{gh}, \cs_0(\uh) +\te+\te')$ in $MK^{\Ic}_{p+1}(\Ac)$. One has
\beq
&&\rho_!(\gh\oplus\hh, \te+\te') = \exp\big(\sqrt{2\pi i}\,\ch_R(\rho)\cdot (\te+\te')\big)\,\Det_R^{-1}(g) \Det_R^{-1}(h)\ , \non\\
&& \rho_!(\widehat{gh}, \cs_0(\uh)+\te+\te') = \non\\
&&\qquad \qquad\qquad  \exp\big(\ch_R(\rho)\cdot  (\Tr\ln(\hh^{-1}\gh^{-1}\widehat{gh}) + \sqrt{2\pi i}(\te+\te'))\big)\,\Det_R^{-1}(gh) \ , \non
\eeq
hence equating these numbers yields the desired result. \cqfd\\

\begin{example} \textup{We look at a simple two-dimensional example. Let $M$ be a closed spin manifold of dimension $p=2$, and $E\to M$ be a hermitian vector bundle with unitary connection $\nabla:\cinf(E)\to \cinf(T^*M\otimes E)$. We take $D_0:\cinf(S\otimes E)\to \cinf(S\otimes E)$ as the twisted Dirac operator acting on the sections of spinors with coefficients in $E$, i.e. $D_0$ is built from the spin connection on $S$ and $\nabla$ on $E$. Proceeding as in Example \ref{sequi} we extend $D_0$ to an invertible selfadjoint operator $D$ on the Hilbert space $H=H_0\oplus H_0$, with $H_0=L^2(S\otimes E)$, and define $\Ac_0=\cinf(\End(E))$ as the algebra of smooth endomorphisms of $E$. In particular $\Ac_0$ is Morita equivalent to the commutative algebra $\cinf(M)$ and is represented on $H$ by acting only on the first summand $H_0$. The spectral triple $(\Ac_0,H,D)$ has analytic dimension $p=2$ and the renormalized Chern character $\ch_R(\rho): \Th\Ac_0\to\cc$ is a 2-dimensional trace which reads, through the isomorphism of pro-vector spaces $\Th\Ac_0\cong \Omh^+\Ac_0$,
\beq
\ch_R(\rho)(a_0) &=& \frac{1}{2\pi i} \int_M \tr(a_0\nabla^2)\ , \non\\
\ch_R(\rho)(a_0da_1da_2) &=& \frac{1}{2\pi i} \int_M \frac{1}{2}\tr(a_0[\nabla, a_1][\nabla, a_2])\ \non
\eeq
for any $a_i\in \Ac_0$. Here $\nabla^2\in \cinf(\Lambda^2 T^*M\otimes \End(E))$ is the curvature two-form of the connection, and $\tr$ is the trace of endomorphisms. Let $g\in \Ac_0$ be a unitary endomorphism of $E$ (gauge transformation) homotopic to the identity, and choose a path $u\in \Ac_0[0,1]$ of unitary endomorphisms from the identity to $g$ (here we don't need to pass to unitalized algebras since $\Ac_0$ is already unital). The renormalized determinant of $g$ is given in terms of the Maurer-Cartan form $\om= u^{-1}su \in \Ac_0 \otimes\Om^1[0,1]$:
$$
\Det_R(g)=\exp\Big( \frac{1}{2\pi i}\int_0^1\!\! \int_M\tr(\om\nabla^2 + \frac{1}{2}\, g^{-1}[\nabla, g][\nabla,\om])\Big)\ .
$$
The integral of the chiral anomaly under the exponential is known in the physics literature as a Wess-Zumino action \cite{WZ}, while $\om$ is the Faddeev-Popov ghost \cite{MSZ}. For two gauge transformations $g,h$ homotopic to the identity the multiplicative anomaly given by Corollary \ref{cma} reads}
$$
\Det_R(gh) = \exp \Big( \frac{1}{2\pi i} \int_M \frac{1}{2} \tr(h^{-1}g^{-1} [\nabla,g][\nabla,h]) \Big) \, \Det_R(g) \Det_R(h)\ .
$$
\end{example}

\end{document}